\documentclass[12pt]{amsart}
\usepackage{epsfig}
\usepackage{amsmath}
\usepackage{amssymb}
\usepackage{amscd}
\usepackage{latexsym}
\usepackage{tabularx}
\usepackage{a4wide}
\usepackage[usenames]{color}
\usepackage{enumerate}
\usepackage{subfigure}
\usepackage{url}
\usepackage{verbatim}

\usepackage{circuitikz}
\ctikzset{bipoles/resistor/height=0.15}
\ctikzset{bipoles/resistor/width=0.4}

\usepackage[
colorlinks, citecolor=blue,
pdfauthor={Omid Amini, Eduardo Esteves}, 
pdftitle={Voronoi Tilings and Degenerations of Line Bundles III},
pdfstartview ={FitV},
]{hyperref}

\usepackage{tikz, float} 
\usetikzlibrary {positioning}

\usetikzlibrary{calc,decorations.markings}
\usetikzlibrary{shapes,snakes}

\numberwithin{equation}{section}

\tikzstyle{Cwhite}=[scale = .8,circle, fill = white, minimum size=3mm] 
\tikzstyle{Cgray}=[scale = .4,circle, fill = gray, minimum size=3mm] 
\tikzstyle{Cblack2}=[scale = .4,circle, fill = black, minimum size=5mm] 
\tikzstyle{Cblack}=[scale = .7,circle, fill = black, minimum size=3mm]
\tikzstyle{C0}=[scale = .9,circle, fill = black!0, inner sep = 0pt, minimum size=3mm]
\tikzstyle{C1}=[scale = .7,circle, fill = black!0, inner sep = 0pt, minimum size=3mm]
\tikzstyle{Cred}=[scale = .4,circle, fill = red, minimum size=3mm] 

\usepackage[matrix,arrow,tips,curve]{xy}
\usepackage{pb-diagram,pb-xy}
\usepackage{verbatim}
\usepackage{mathrsfs}
\usepackage{color}
\usepackage[leqno]{amsmath}

\newcommand{\E}{\mathbb E}

\newcommand{\m}{\mathfrak m}

\small\normalsize
\newcommand{\M}{\mathcal{M}}
\newcommand{\wt}{\widetilde}
\renewcommand{\:}{\colon}
\renewcommand{\deg}{\mathrm{deg}\,}
\newcommand{\T}{\mathcal T}
\newcommand{\I}{\mathcal I}

\newcommand{\lra}{\longrightarrow}
\newcommand{\ol}{\overline}

\newcommand{\IP}{{\mathbf P}}
\newcommand{\IR}{{\mathbf R}}
\newcommand{\IS}{{\mathbf S}}

\newcommand{\Quot}{{\mathrm{Quot}}}

\renewcommand{\IJ}{\mathbf J}

\newcommand{\bfd}{\mathbf b}
\newcommand{\bfb}{\mathbf b}
\newcommand{\bfc}{\mathbf c}

\newtheorem{thm}{Theorem}[section]

\newtheorem{prop}[thm]{Proposition}

\newtheorem{claim}[thm]{Claim}
\newtheorem{cor}[thm]{Corollary}

\newtheorem{remark}[thm]{Remark}

\newtheorem{defi}[thm]{Definition}

\newcommand{\dl}{\mathfrak d}

\def\L{\mathcal L}
\def\O{\mathcal O}
\newcommand{\Q}{\mathbb{Q}}

\def\X{\mathcal X}

\newcommand{\Z}{\mathbb{Z}}

\newcommand{\R}{\mathbb{R}}

\newcommand{\N}{\mathbb{N}}

\renewcommand{\div}{\mathrm{div}}

\renewcommand{\k}{\kappa}

\newcommand{\g}{\mathfrak g}

\newcommand{\F}{\mathcal F}

\newcommand{\he}{\mathrm{h}}  
\newcommand{\te}{\mathrm{t}}  
\newcommand{\comp}[1]{\overbar[.7]{#1}} 
\newcommand{\mg}{\mathscr M} 
\newcommand{\mgbar}{\comp\mg} 

\usepackage{scalerel}

\newcommand{\ind}[1]{{_{\scaleto{#1}{4.8pt}}}}
\newcommand{\indm}[1]{{_{\scaleto{#1}{3.2pt}}}} 
\newcommand{\indmbar}[1]{{_{\scaleto{#1}{4.6pt}}}} 
\newcommand{\indmbi}[2]{_{{\scaleto{#1}{3.2pt}}_{\hspace{-.03cm}#2}}} 

\makeatletter
\newsavebox\myboxA
\newsavebox\myboxB
\newlength\mylenA

\newcommand*\overbar[2][0.75]{%
    \sbox{\myboxA}{$\m@th#2$}%
    \setbox\myboxB\null
    \ht\myboxB=\ht\myboxA%
    \dp\myboxB=\dp\myboxA%
    \wd\myboxB=#1\wd\myboxA
    \sbox\myboxB{$\m@th\overline{\copy\myboxB}$}
    \setlength\mylenA{\the\wd\myboxA}
    \addtolength\mylenA{-\the\wd\myboxB}%
    \ifdim\wd\myboxB<\wd\myboxA%
       \rlap{\hskip 1\mylenA\usebox\myboxB}{\usebox\myboxA}%
    \else
        \hskip -0.5\mylenA\rlap{\usebox\myboxA}{\hskip 0.5\mylenA\usebox\myboxB}%
    \fi}
\makeatother

\usepackage{marginnote}
\newcommand{\hooklongrightarrow}{\lhook\joinrel\longrightarrow}

\numberwithin{equation}{section}

\newcommand\sbullet[1][.5]{\mathbin{\vcenter{\hbox{\scalebox{#1}{$\bullet$}}}}}

\renewcommand{\o}{\mathfrak o}
\renewcommand{\P}{\mathbf P}
\newcommand{\Gm}{\mathbf{G_m}}

\newcommand{\x}{\textsc{x}} 

\newcommand{\bt}{\mathrm{bt}}

\title{Voronoi Tilings, Toric Arrangements and Degenerations of Line Bundles III}
 \author{Omid Amini}
 \author{Eduardo Esteves}
 \date{\today}
  \address{CNRS - CMLS, \'Ecole Polytechnique, Palaiseau, France}
\email{omid.amini@polytechnique.edu}

\address{Instituto Nacional de Matem\'atica Pura e Aplicada, Estrada Dona Castorina 110,
22460-320 Rio de Janeiro RJ, Brazil}
\email{esteves@impa.br}
\begin{document}

\maketitle
\begin{abstract}
 We describe limits of line bundles on nodal curves
  in terms of toric arrangements associated to Voronoi tilings of Euclidean spaces. 
These tilings encode information on the relationship between the 
possibly infinitely many limits, and ultimately give rise to a new
definition of \emph{limit linear series}. This article and the first
two that preceded it are the first in a series aimed to explore this new
approach. 

In Part I, we set up the combinatorial framework and 
showed how graphs weighted with integer lengths associated to the
edges provide tilings of Euclidean spaces by certain polytopes
associated to the graph itself and to its subgraphs. 

In Part II, we described the arrangements of 
toric varieties associated to the tilings of Part I in several ways: using normal fans, as unions of orbits, by equations and as degenerations of tori. 

In the present Part III, we show how these combinatorial and 
toric frameworks allow us to describe all \emph{stable limits} 
of a family of line bundles along a degenerating family of curves. 
Our main result asserts that the collection of all these limits is 
parametrized by a connected 0-dimensional closed substack of 
the Artin stack of all torsion-free rank-one sheaves on the limit
curve. Moreover, we thoroughly describe this closed substack and all
the closed substacks that arise in this way as certain torus quotients
of the arrangements of toric varieties of Part II determined by the Voronoi
tilings of Euclidean spaces studied in Part I.
\end{abstract}

\tableofcontents

     \addtocontents{toc}{\protect\setcounter{tocdepth}{1}}

\section{Introduction}

This is a sequel to our previous works~\cite{AE1, AE2} whose aim is to achieve the description of stable limits of line bundles on nodal curves by means of graph theory and toric geometry. This is motivated by the desire to understand all the possible limits of linear series $\g^r_d$ over any sequence $X_1, X_2, \dots$ of smooth projective curves of genus $g$ whose corresponding points $x_1, x_2, \dots$ in $\mg_g$ converge to a given point $x$ on the Deligne-Mumford compactification $\mgbar_g$.

 \medskip
 
In this introduction, after providing an overview of the previous works, we focus on the contribution of this paper and its preceding companions and briefly discuss the content of our forthcoming work.

\subsection{Overview.} Nodal curves are curves (one-dimensional, reduced,
connected but not necessarily irreducible projective schemes over an
algebraically closed field) 
that fail to be smooth in the weakest possible form: 
the singularities
are normal crossings, that is, ordinary nodes. Among them, the most
important are (Deligne--Mumford) stable curves, characterized as those having ample
canonical bundle. This is in fact the property that allows for the
construction of their moduli
space, $\comp{\mg_g}$, where $g$ stands for the (arithmetic) genus. 
One of the key properties of stable curves is the Stable Reduction Theorem,
which says that a family of stable curves parameterized by 
the germ of a punctured smooth curve can be completed, after a finite base change,
in a unique way to a family over the whole germ. It implies that the
moduli of stable curves is complete, in fact projective, and is thus a
compactification for the moduli of smooth curves of genus at least
two. However, no such thing
holds, in general, for line bundles and linear series over curves.

A linear series over a curve $X$ is simply the data of a line bundle and a
linear subspace of the space of sections of that bundle. If $d$ is the
(total) degree of the line bundle and $r$ is the projective dimension
of the subspace, we say we have a $\mathfrak
g^r_d$. If $r\geq 0$, it corresponds to a rational map  
$X\dashrightarrow\P^r$, whence the importance of their study for
understanding the projective geometry of curves. Line bundles over a curve are
parametrized by the Picard scheme, and the linear series by fibrations
by Grassmannians over the Picard scheme. It is all well if $X$ is
smooth, as the Picard scheme has projective connected
components indexed by the degree $d$; the one with degree 0 is even an
algebraic group, the Jacobian. 

However, the Picard scheme of a stable curve has projective components
only if the curve is of compact type. Furthermore, whereas for smooth
curves, a $\mathfrak g^r_d$ gives rise to a rational map that can even
be uniquely extended to the whole curve, 
the map induced by a $\mathfrak g^r_d$
on a stable curve may fail to be defined on whole components of the
curve. Finally, there is nothing similar to stable reduction;
quite to the contrary, the trivial line bundle over the total space of
a family of smooth curves parameterized by 
the germ of a punctured smooth
curve can be completed in infinitely many ways over the whole germ 
if the family of curves is completed by adding a reducible
stable curve.

A compactification of the Picard scheme of an irreducible curve was 
suggested by Mumford in the sixties, carried out by D'Souza~\cite{dsouza} in the
seventies, and thoroughly studied by Altman and Kleiman~\cite{AK80, AK79} and others in the
years that followed. As it was later observed by Eisenbud and Harris~\cite{EH83},
it is possible to complete the proof of the celebrated Brill--Noether
Theorem by considering this compactification over rational irreducible
cuspidal curves. The compactification  parameterizes torsion-free,
rank-one sheaves, and the notion of linear series extends simply by just
considering spaces of sections of these sheaves. In fact, in this setup it is even possible to
extend the notion of divisors: these are the pseudo-divisors studied
by Hartshorne~\cite{Hartshorne}. Finally, there are fine proper moduli spaces
parameterizing torsion-free, rank-one sheaves of a given degree for
families of irreducible curves. 

Though even in codimension one (in their moduli), stable curves 
fail to be irreducible, fortunately, general reducible stable
curves are of compact type. Curves of compact type are those curves 
for which the connected components of the Picard scheme are 
compact, indeed projective. Eisenbud and Harris 
studied degenerations of line bundles
along families of smooth curves degenerating to curves of compact
type, and coined the term \emph{limit linear series}~\cite{EH86}. 
Remarkably, they observed that
it was actually useful that the family of line bundles could degenerate to many
limit line bundles, as just one limit would rarely carry enough information
about the degeneration. They found out that for many applications it
was important, and enough, to consider as many limit line bundles as
components of the limit stable curve, each focusing on one of the
components. Their study allowed for a deep understanding of the moduli
of stable curves up to codimension two, and found many applications 
in describing the geometry of 
a general curve and the theory of 
divisors of $\comp{\mg_g}$~\cite{EH86, EH87-WP, EH87mon}. It is thus
no surprise they wished for an extension of their theory to more
general curves; as they wrote in \cite{EHBull}: ``The known special
cases suffice for many applications, but there is probably a gold mine
awaiting a general insight.''

\vspace{.2cm}

Already in the late seventies, the study of degenerations to more complicated nodal
curves began, in the hands of Oda and
Seshadri~\cite{OS79}. To do away with the problem of many 
limits line bundles, they
introduced the notion of stability, one that had naturally arisen in
the study of vector bundles over smooth curves and the construction of
their moduli space by Mumford's Geometric Invariant Theory. Out of the
many limits a family of line bundles could have, at most one would be 
stable, and a finite nonzero number of them would be
semistable. Moreover, the semistable limits would be equivalent in the
sense they would all have the same quotients in their 
Jordan--H\"older filtrations. The convenient notion of stability has been
centerpiece in all the studies that followed, including in the
landmark construction of a compactification of the relative Picard
scheme over the whole moduli of curves by Caporaso in the early 90's~\cite{caporaso}.

However, as far as limits of linear series are considered, stability
appears to be very restrictive. The line bundles Eisenbud and Harris
considered, one for each component of the stable curve, would rarely be semistable; even if one could force one of
them to be semistable, by deforming the notion of stability, the
others would seldom be. Moreover, generally semistable line bundles in the same
equivalence class have almost completely unrelated spaces of 
sections.

On the other hand, an approach similar to that by Eisenbud and Harris,
choosing for each component of the limit curve a ``best'' limit line
bundle, was tentatively carried out by many: Ziv Ran had an early
draft on this already in the 80's, whereas the second author, 
in collaboration with
Medeiros and Salehyan~\cite{EM, ES07}, used this approach to study limits of
Weierstrass points for a wide class of stable curves in the nineties.
However, one could not carry the approach further.

\vspace{.2cm}
In the middle of the first decade of the present century, Osserman
\cite{Os06-moduli} 
introduced a new idea, that not only extended the theory to include
curves in positive characteristic, but also modified substantially the
approach by Eisenbud and Harris. Osserman's idea was remarkably
simple, even though it seemed at first to complicate further the
study: One should consider not only the ``best'' limit line bundle for each
component of the limit curve, but also all of those in between, more
precisely, all of
those with nonnegative degree on every component; the intermediate
limits carried only partial information over each component, but
perhaps crucial information
nonetheless. And, in fact, if Eisenbud's and Harris's limit
linear series was not refined, Osserman found out that the extra limits did
carry more information. For instance, Osserman and the second author
discovered later that they described limits of divisors of the family
of linear series that cannot be accounted for only by the ``best''
line bundles~\cite{EO}. 

The approach by Osserman introduced a new challenge: How to account
for, how to keep track of all the data that show up when considering
all effective limit line bundles? Osserman introduced his approach first
only for two-component curves of compact type, the simplest case, 
and recently extended the approach to more 
general curves~\cite{Os-pseudocompact, Oss16}, specifically, curves of 
pseudo-compact type where a parallel theory to that by
Eisenbud and Harris can be established.  

\vspace{.2cm}

This was not the only challenge: another was to consider one-parameter
families of curves whose total space was not smooth. For such
families, when the limit curve is not of compact type, it might as
well happen that a family of line bundles has no line bundle in the
limit. One can blow up the total space and replace the stable curve by
one of its semistable models, splitting apart each node to introduce a
chain of a variable number of rational smooth curves. This way one
obtains a limit line bundle, but on a different curve! Worse, one has
to deal with a whole set of new variables, one for each node: the
number of rational curves in each chain. This procedure, 
semistable reduction, has already appeared when dealing with curves of 
compact type; it did not lead to many difficulties in applications
because in them, a class of curves was considered, and the class did
not change by semistable reduction. For more general curves, as it
became apparent in the work by Esteves and Medeiros~\cite{EM} on 
curves with two components, the specific type of reduction 
that appeared had enormous influence on, for instance, the limits of
Weierstrass points. 

\vspace{.2cm}

The whole study has been combinatorially intensive from the start,
already in Oda--Seshadri ~\cite{OS79}. Farther than obstructing the study of the
problem, combinatorics has been one of the main tools used by all
those that ventured in the field. It has even spanned a new line of
approach, using tropical and non-Archimedean geometry, that has 
been very fruitful. In fact, the first author in a joint work with
Baker~\cite{AB15} has introduced a notion of limit linear series which
also extends, 
in a somewhat more combinatorial way originating from non-Archimedean 
analysis, the definition by Eisenbud and Harris. He has used this approach as an alternative tool in the study of reduction of Weierstrass points and their distributions~\cite{Ami-W}. Moreover, this circle of ideas has 
 been used by Jensen and Payne to prove specific cases of the Maximal Rank Conjecture~\cite{JP1, JP2, JP3}, as well as in the recent work by Farkas, Jensen and Payne on the Kodaira dimension of the moduli spaces of curves of genus 22 and 23~\cite{FJP}.  A comparison of Osserman's
approach~\cite{Os-pseudocompact} to the work of Amini and Baker~\cite{AB15} can be found in~\cite{Oss19}.
 
\smallskip

The overall approach we take in these series of works features a new interplay between the combinatorics and
the geometry. We do not do away with any data (as it is done,
sometimes harmlessly for the applications in mind), but we use the
combinatorics to organize all the information. In a nutshell we prove
in the present paper that the collection of all limit torsion-free, rank-one sheaves of a
family of line bundles along a family of smooth curves degenerating to
a nodal curve $X$ is parameterized by a connected 0-dimensional closed
substack $\mathfrak I$ of the Artin stack $\mathbf J$ of all torsion-free, rank-one sheaves on
$X$. 

More significatively, we characterize combinatorially all the
closed substacks $\mathfrak I$ that arise from degenerations. The characterization
is given by Theorems~\ref{thm:main4}~and~\ref{regeneration}, 
which are the main results in the article. In short, the $\mathfrak I$ are
certain torus quotient of a combinatorial arrangement of toric varieties. Thus
our theorems provide a far-reaching generalization to all nodal 
curves of the combinatorial-toric approach carried out by Esteves
and Medeiros in \cite{EM} for the case of stable curves with two irreducible components.

\vspace{.2cm}

In a future work, if a family of linear series is given, in addition to the substack
$\mathfrak I$, we will show the existence of 
a closed substack $\mathfrak I'$ of $\mathbf G$, the Grassmann fibration over
$\mathbf J$ parameterizing vector spaces of sections of torsion-free, rank-one sheaves on
$X$, satisfying two conditions: First, the stack $\mathfrak I'$ will
lie over $\mathfrak I$. 
Second, the induced relative torsion-free rank-one sheaf $\mathcal I$ on
$X\times\mathfrak I'$ over $\mathfrak I'$ and the induced locally free subsheaf
$\mathcal V$ of the pushforward $p_{2*}\mathcal I$ will be such that
their restrictions over the points of $\mathfrak I'$ parameterize all the limits of the family of
linear series. We will describe $\mathfrak I'$ combinatorially as
well; it is actually very similar in esprit to the description we give of $\mathfrak I$.

The connection with the work by Osserman will be elaborated in our future study. In the case of a two-component curve of compact type, this has been carried out in detail in the
doctoral thesis by Rizzo \cite{Rizzo}, who explains how the collection
of linear series considered by Osserman as a limit linear series for a
two-component curve of compact type can actually be viewed as members of a $\Gm$-equivariant
family of linear series on $X$ parameterized by a 2-punctured chain of rational
curves. The quotient of this chain by $\Gm$ is the stack
$\mathfrak I'$ truncated in nonnegative degrees.

\subsection{The present work.} 

We will now proceed to
describe the work we do in this paper, in continuation to what we did
in \cite{AE1} and \cite{AE2}. Later we will point out what lies ahead
in the path we are taking.

Fix a connected nodal curve $X$ over an algebraically closed field
$\k$. Consider its associated dual
graph $G=(V,E)$, which is a pair consisting of a vertex set $V$ in one-to-one
correspondence with the set of irreducible components of $X$, and 
an edge set $E$ in one-to-one correspondence with the set of nodes of $X$. 
An edge connects two vertices if the corresponding node lies on the two corresponding 
components. For our purposes, we discard those edges that form a
loop.

Let now $\mathbb E$ be the set of all the oriented edges (also called arcs) obtained out of $E$: for each
edge, there are two possible arcs, pointing to the two different vertices
connected by the edge. For $e\in\mathbb E$, we write $e=uv$ to mean
that $e$ is an arc connecting $u$ to $v$, even if it might not be
the only one. Also, we let $\te_e$ denote the tail and $\he_e$ the
head of $e$. In addition, $\ol e$ denotes the same edge with
the reverse orientation. 

Recall that given a commutative ring $A$, one associates to the graph
$G$ the complex of 
$A$-modules 
$$
d_A\colon C^0(G,A) \to C^1(G,A).
$$
Here, $C^0(G,A)$ is the $A$-module of functions $V\to A$, and 
$C^1(G,A)$ is the $A$-module of all functions $f\colon\mathbb E\to A$ satisfying 
$f(\ol e)=-f(e)$ for each $e\in\mathbb E$. And $d_A(f)(e)=f(v)-f(u)$ for
each $e=uv\in\mathbb E$. 

The characteristic functions $\chi_{\indm v}$, for $v$ in $V$, form a basis of
the $A$-module $C^0(G,A)$, whereas the functions $\chi_{\indm e}-\chi_{\indmbar{\ol
  e}}$, for a collection of $e\in\E$ giving an orientation to the whole
$G$, form a basis of $C^1(G,A)$. There are
bilinear forms $\langle\,,\rangle$ on $C^0(G,A)$ and $C^1(G,A)$ satisfying
\begin{align*}
&\langle \chi_{\indm v},\chi_{\indm w}\rangle=\delta_{v,w}\quad\text{for }v,w\in V;\\
&\langle\chi_{\indm e}-\chi_{\indmbar{\ol e}},\chi_{\ind f}-\chi_{\ind{\ol
  f}}\rangle=\delta_{e,f}-\delta_{\ol e,f}\quad\text{for
  }e,f\in\E.
\end{align*} 

Define the homomorphism 
$d^*_A\colon C^1(G,A)\to C^0(G,A)$ by putting $d^*_A(\chi_{\indm e}-\chi_{\indmbar{\ol e}}):=
\chi_{\indm v}-\chi_{\indm u}$ for each $e=uv\in\E$. Then $d^*_A$ is the adjoint to  $d_A$, that is,
$\langle f,d^*_A(h)\rangle=\langle d_A(f),h\rangle $ for each $f\in C^0(G,A)$ and $h\in C^1(G,A)$.
In addition, the degree map $\deg\colon C^0(G,A)\to A$, sending
$f$ to $\sum_{v\in V} f(v)$, is a cokernel for $d_A^*$. The kernel is $H^1(G,A)$.

Let $H_{0,A}:=\{f\in
C^0(G,A)\,|\,\deg(f)=0\}$ and $F_A:=\text{Im}(d_A)$. 
Let $\Delta_A:=d^*_Ad_A$, the Laplacian of $G$. The homomorphism 
$d^*_A$ induces an injection $F_A\to H_{0,A}$. For $A=\mathbb R$, it
is a bijection, and the bilinear form $\langle \,,\rangle$ on $C^1(G,\R)$ 
induces by restriction a norm on $F_{\mathbb R}$
corresponding via $d^*_{\mathbb R}$ to the quadratic form $q$
on $H_{0,\mathbb R}$ satisfying $q(f)=\langle f,\Delta_{\mathbb R}(f)\rangle$ for
each $f\in C^0(G,\mathbb R)$. 

\vspace{.2cm}

In \cite{AE1}, we described a certain family of tilings of
$H_{0,\mathbb R}$ by polytopes. Each tiling consists of a family of
polytopes covering $H_{0,\mathbb R}$ such that
\begin{itemize}
\item each face of a polytope which is in the tiling  belongs itself to the tiling; and
\item the intersection of a finite number
of polytopes in the tiling is a face of each of the polytopes.
\end{itemize}
By removing from a polytope in the tiling all
the faces of positive codimension it contains, we get the
corresponding open face. The open faces form then a stratification of 
the whole space $H_{0,\mathbb R}$. We call tiles the polytopes of
maximum dimension.

For instance, let $\Lambda_A:=\text{Im}(d^*_A)$. Then 
$\Lambda_{\mathbb R}=H_{0,\mathbb R}$ and $\Lambda_{\mathbb Z}$ is a
sublattice of $H_{0,\mathbb Z}$ of finite index equal to the number of
spanning trees of $G$, by the Kirchhoff Matrix Tree Theorem. The
\emph{standard Voronoi tiling} of $G$ is the Voronoi decomposition 
$\text{Vor}_G$ of
$H_{0,\mathbb R}$ with respect to $\Lambda_{\mathbb Z}$ and $q$. The
tiles are 
$$
\text{Vor}_G(\beta):=\bigl\{\eta\in H_{0,\mathbb R}\,|\,
q(\eta-\beta)\leq q(\eta-\alpha)\text{ for every }\alpha\in 
\Lambda_{\mathbb Z}-\{\beta\}\bigr\}
$$
for $\beta\in \Lambda_{\mathbb Z}$. 

This is one of the infinitely many  tilings we
consider. There are variants of it, that we call \emph{twisted mixed
Voronoi tilings} and denote $\text{Vor}^{\m}_{G,\ell}$. 
Though the standard Voronoi tiling is homogeneous,
meaning all tiles are translates of the tile centered at the origin,
$\text{Vor}_G(O)$, a twisted mixed Voronoi tiling is obtained by
putting together translations of the tiles $\text{Vor}_H(O)$
associated to connected spanning subgraphs $H$ of $G$. More precisely,
the twisted Voronoi
tiling $\text{Vor}^{\m}_{G,\ell}$ depends on $\m\in C^1(G,\Z)$ (the
``twisting'') and an edge length
function $\ell\: E\to\N$; its tiles are the polytopes 
$d^*_\R(\dl^\m_f)+\text{Vor}_{G^\m_f}(O)$ for
$f\in C^0(G,\Z)$ with $G^\m_f$ connected, 
where $\dl^\m_f\in C^1(G,\R)$ is a modification of $d_\Z(f)$, namely 
$$
\dl^\m_f(e):=\frac{\delta^{\m,\ell}_e(f)-\delta^{\m,\ell}_{\ol
  e}(f)}{2},\text{ where }
\delta^{\m,\ell}_e(f):=\Big\lfloor\frac{f(v)-f(u)+\m(e)}{\ell(e)}\Big\rfloor
\text{ for each $e=uv\in\E$,}
$$
and $G^\m_f$ is the spanning subgraph of $G$ obtaining by removing the edges
$e\in\E$ for which $\dl^\m_f(e)\not\in\Z$. We refer to \cite{AE1}
for a thorough presentation of these tilings. 

\vspace{.2cm}

In the present article we establish a correspondence between the stratifications of 
$H_{0,\mathbb R}$ associated to the $\text{Vor}^{\m}_{G,\ell}$
and the stable limits of line bundles in
one-parameter 
smoothings of $X$. Properties of each stratum, and
the way they fit together in the stratification of $H_{0,\mathbb R}$
are reflected in properties of the limits and the relationship between
them. 

More explicitly, let $\pi\:\mathcal X\to B$ be a (one-parameter)
smoothing of $X$. Here, $B$ is the spectrum of $\k[[t]]$ and 
$\pi$ is a projective flat morphism whose generic fiber is
smooth and special fiber is isomorphic to $X$. We fix such an
isomorphism. Let $\eta$ and $o$ be the generic and special points of
$B$. The total space $\mathcal X$ is regular except possibly at the nodes of
$X$. For $e\in E$, the completion of the local ring of $\mathcal X$ at
the corresponding node $N_e$ is $\k[[t]]$-isomorphic to
$\k[[u,v,t]]/(uv-t^{\ell_e})$ for a certain $\ell_e>0$, called the
\emph{singularity degree} of $\pi$ at $N_e$. If $\ell_e=1$, then
$\mathcal X$ is regular at $N_e$. If all $\ell_e=1$, then
$\pi$ is said to be \emph{Cartier}; it is \emph{regular} if $\mathcal
X$ is regular at \emph{all} the nodes of $X$. A finite base change is obtained by
sending $t$ to $t^n$ for some $n$. The resulting family is similar to
the original one: the special fiber is the same, the generic fiber is
a base field extension of the original one, but the
singularity degrees $\ell_e$ change to $n\ell_e$.

Let $L_\eta$ be an invertible sheaf on the generic fiber. If $\pi$ is
Cartier, it extends to an \emph{almost invertible} sheaf $\mathcal L$
on $\X$, a sheaf that is invertible at all $N_e$. It is
not unique, as $\mathcal L\otimes\mathcal O_{\mathcal X}(\sum
f(v)X_v)$ is another extension for each $f\in C^0(G,\mathbb
Z)$. (Here $X_v$ is the component of $X$ corresponding to $v\in V$,
which can and will be viewed as a Cartier divisor of $\X$ because
$\pi$ is Cartier.) For
general $\pi$, the sheaf $L_\eta$ extends to a relatively
torsion-free, rank-one sheaf $\mathcal I$ on $\X/B$, that is, a $B$-flat 
coherent sheaf on $\mathcal X$ whose fibers over $B$ are torsion-free,
rank-one. Again, it is not the unique extension: in~\cite{Esteves01} a procedure
similar to the one explicited above shows how to change $\mathcal I$
into other extensions. Furthermore, one could do a finite base change
to $\pi$, extend $L_\eta$ to the new generic fiber and consider its
extensions. Of course, they will be extensions on a different total
space. But the special fibers are the same, and thus the restrictions
of all these extensions to $X$ are torsion-free, rank-one sheaves that
we call the \emph{stable limits} of $L_\eta$.

Let $\mathbf J$ denote the Artin stack parameterizing torsion-free,
rank-one sheaves on $X$. It is the disjoint union of the closed and
open substacks $\mathbf J^d$, each parameterizing sheaves $I$ with
degree $d$, that is, with $\chi(I)-\chi(\mathcal O_X)=d$. Letting
$d:=\deg(L_\eta)$, we may consider the subset $\mathfrak I$ of
$\mathbf J^d$ parameterizing all the stable limits of $L_\eta$. We
proceed to describe $\mathfrak I$ thoroughly.

First, we give a meaningful structure to $\mathbf J$ as a quotient
stack, as follows. Fixing 
$\mathbf b\in C^0(G,\Z)$, we let 
$$
\mathbf J^{\mathbf b}:=\prod_{v\in V}\mathbf J_v^{\mathbf b(v)},
$$
where $\mathbf J_v^{\mathbf b(v)}$ parameterizes torsion-free, rank-one
sheaves of degree $\mathbf b(v)$ on the component $X_v$ for each $v\in
V$. Over $\mathbf J^{\mathbf b}$ we construct a scheme $\mathbf
R_{\mathbf J^{\mathbf b}}$ parameterizing gluings along the
nodes $N_e$ associated to $e\in E$ 
of the sheaves given by points of $\mathbf J^{\mathbf b}$, their
modifications and their degenerations. 

More precisely, given torsion-free, rank-one sheaves $K_v$ on $X_v$ for
each $v\in V$, we view their gluings along the nodes as subsheaves of
the direct sum $\oplus_vK_v$ whose quotients are supported with length
1 on each and every node $N_e$. This
allows for ``degenerate'' gluings, and thus for torsion-free sheaves
that may fail to be invertible at part or all of the $N_e$. The
parameter space for these subsheaves is a product of $\mathbf P^1$,
one for each $e\in E$, and describes an irreducible component of the
fiber of  $\mathbf R_{\mathbf J^{\mathbf b}}$ over the point on
$\mathbf J^{\mathbf b}$ parameterizing the $K_v$.

The other components are obtained by modifying the $K_v$ as
follows. Fix $\mathbf c\in C^1(G,\Z)$. For each $c\in C^1(G,\Z)$, put
\[
K_v^{\bfc -c}:=K_v\otimes\mathcal O_{X_v}\Big(
\sum_{\substack{e\in\mathbb E\\ \he_e=v}}(\bfc(e)-c(e))N_e\Big)
\]
and do the same gluing as above; we obtain another irreducible
component (if $c\neq\bfc$) 
of the fiber of  $\mathbf R_{\mathbf J^{\mathbf b}}$ over the same point on
$\mathbf J^{\mathbf b}$ parameterizing the $K_v$, and all components
are obtained this way. 

The fibers of $\mathbf R_{\mathbf J^{\mathbf b}}$ over 
$\mathbf J^{\mathbf b}$ are thus arrangements of an infinite number of
simple toric varieties. Each fiber is the same: what we call the
arrangement of toric varieties or simply toric
tiling $\mathbf R$ 
associated to the Voronoi decomposition of $C^1(G,\R)$ in hypercubes with respect to
$C^1(G,\Z)$. More precisely, each Voronoi tile is a rational polytope,
to each rational polytope we may associate its normal fan, and to the
normal fan the corresponding toric variety. The polytopes in the
Voronoi decomposition form a complex, in the sense that each two of
them intersect in a common face, when they intersect. We may thus glue the
toric varieties associated to each two polytopes by identifying the
orbit closures of the common face in each variety with the toric
variety associated to the face itself, viewed as a polytope of smaller
dimension. This gives us $\mathbf R$. 

There is nothing special about the above Voronoi decomposition with
regard to the above construction. It can be carried out for any tiling
of an Euclidean space by rational polytopes intersecting each other in
faces. As we have already mentioned, we consider other tilings in the present article.

The scheme $\mathbf R_{\mathbf J^{\mathbf b}}$ parameterizes
torsion-free, rank-one sheaves with a lot of redundancy. Two groups act
on it independently. First, there is a natural action of 
$\mathbf G_{\mathbf m}^V/\Gm$,
where we view $\Gm$, the multiplicative group of $\k$, embedded
diagonally in $\mathbf G_{\mathbf m}^V$. The action is given by
observing that $\mathbf G_{\mathbf m}^V$ is the automorphism group of 
$\oplus_v K_v^{\bfc -c}$ for each point on $\mathbf J^{\mathbf b}$
parameterizing the $K_v$ and each $c\in C^1(G,\Z)$, and that a
subsheaf is fixed under the diagonal action of $\Gm$. The action
preserves each fiber of $\mathbf R_{\mathbf J^{\mathbf b}}$ over
$\mathbf J^{\mathbf b}$ and each component of that fiber.

The second group action moves fibers around. Recall $H^1(G,\Z)$, the
group of cycles of $G$, the kernel of $d^*_\Z$. For each $\gamma\in
C^1(G,\Z)$ and torsion-free rank-one sheaves $K_v$ on $X_v$ for each
$v\in V$, put:
$$
K_v^{-\gamma}:=K_v\otimes\mathcal O_{X_v}\Big(-
\sum_{\substack{e\in\mathbb E\\ \he_e=v}}\gamma(e)N_e\Big).
$$
If $\gamma\in H^1(G,\Z)$ and the $K_v$ are parameterized by a point
$s\in\mathbf J^{\mathbf b}$, so are the $K_v^{-\gamma}$, parameterized
by a point we denote by $\tau^{\gamma}(s)$. Then $\tau^{\gamma}$ is an
automorphism of $\mathbf J^{\mathbf b}$. It lifts in a rather trivial
way to an automorphism of $\mathbf R_{\mathbf J^{\mathbf b}}$: Since 
$\bfc -(c+\gamma)=-\gamma+(\bfc -c)$, we may associate to a 
point on the component corresponding to 
$c+\gamma$ of the fiber of $\mathbf R_{\mathbf J^{\mathbf b}}$ over
$s$ the point on the component corresponding 
to $c$ of the fiber of $\mathbf R_{\mathbf J^{\mathbf b}}$ 
over $\tau^{\gamma}(s)$ parameterizing the same subsheaf of the same 
direct sum, for each $c\in C^1(G,\Z)$. The action of $H^1(G,\Z)$ may move
fibers and is free, as even if a fiber is fixed, its components are
not.

The two actions are independent and $\mathbf J^d$ can be obtained as
the quotient stack:
\begin{equation}\label{stackJ}
\mathbf J^d=\Bigg[\frac{\mathbf R_{\mathbf J^{\mathbf b}}}
{H^1(G,\Z)\times \mathbf G_{\mathbf m}^V/\Gm}\Bigg].
\end{equation}

We may thus describe the collection $\mathfrak I$ of stable limits of a
sheaf $L_\eta$ by describing its inverse image in $\mathbf R_{\mathbf
  J^{\mathbf b}}$. Our Theorem~\ref{thm:main4} claims that this
inverse image is the disjoint union of certain connected subschemes of certain
fibers of $\mathbf R_{\mathbf J^{\mathbf b}}$ over $\mathbf
J^{\mathbf b}$, each subscheme isomorphic to its image in the
quotient
$$
\mathbf S^d:= \Bigg[\frac{\mathbf R_{\mathbf J^{\mathbf b}}}
{H^1(G,\Z)}\Bigg],
$$
all the images being the same. Under chosen identifications we may
view each subscheme as a subscheme of the arrangement of toric
varieties $\mathbf R$. We prove this subscheme is
$Y^{a,b}_{\ell,\mathfrak m}$ for certain choices of $\ell$, $\m$, $a$
and $b$. Each $Y^{a,b}_{\ell,\mathfrak m}$ is itself an arrangement of toric
varieties of dimension $|V|-1$ which was thoroughly described in
\cite{AE2}.

Here we explain the parameters defining 
$Y^{a,b}_{\ell,\mathfrak m}\subseteq\mathbf R$: 
characters $a\: C^1(G,\Z)\to\Gm(\k)$ and $b\: H^1(G,\Z)\to\Gm(\k)$, an
edge length function $\ell\: E\to\N$ and an element 
$\mathfrak m\in C^1(G,\Z)$ we call a twisting. These data arise from
the smoothing $\pi\:\mathcal X\to B$ and from $L_\eta$, as follows: 
The length function is
simply the collection of singularity degrees $\ell_e$; and the character $a$ keeps
track of the infinitesimal data of the arc defined by $\pi$ on the
moduli $\comp{\mg_g}$ or on a versal deformation space of $X$. 
The character $b$ describes the gluing data of
an almost invertible stable limit, if $L_\eta$ admits one, and then we may
set $\mathfrak m=0$, no twisting is necessary. If not, then
$L_\eta$ admits an almost invertible limit on $X^\ell$, the semistable model
of $X$ obtained by splitting the branches of each $N_e$ apart, and
connecting them by a chain of $\ell_e-1$ smooth rational curves for
each $e\in E$. It even admits an admissible almost invertible limit, meaning
an almost invertible limit whose restriction to each component of each chain
has degree zero, but possibly one, where the degree is one. Then $b$ is
related to the gluing data of an admissible almost invertible limit and the
twisting 
$\mathfrak m$ keeps track of where that limit has degree one on each
added chain.

Different choices of $a$, $b$ and $\mathfrak m$ may yield the same
subscheme $Y^{a,b}_{\ell,\mathfrak m}\subseteq\mathbf R$. We have left
this analysis for a later work. But the structure of
$Y^{a,b}_{\ell,\mathfrak m}$ depends on $\ell$ and $\mathfrak m$ only,
its equations being a deformation of the equations defining
$Y^{1,1}_{\ell,\mathfrak m}$, which we denote by $Y_{\ell,
  \m}^{\bt}$ and call the basic toric tiling, as pointed out in \cite{AE2}, Prop.~4.6.

In \cite{AE2}, Section 4, we explained how $Y_{\ell,\m}^{\bt}$
is determined from
the tiling $\text{Vor}^{\mathfrak m}_{G,\ell}$ of $H_{0,\R}$. There we
remark we can naturally view $Y_{\ell,\m}^{\bt}$ and its deformations
$Y^{a,b}_{\ell,\mathfrak m}$ as closed subschemes
of $\IR$. 

Finally, we may consider the subscheme
$Y^{a,b}_{\ell,\m}\subseteq\mathbf R$ for arbitrary choices of $a$, $b$, $\ell$
and $\m$, and under chosen identifications, as a subscheme of a 
fiber of $\mathbf R_{\mathbf J^{\mathbf b}}$ over
$\mathbf J^{\mathbf b}$. If we denote by $\mathfrak I$ the image of
this subscheme in the quotient $\mathbf J^d$, our
Regeneration Theorem~\ref{regeneration} claims
that $\mathfrak I$ parameterizes the collection of stable limits of an
invertible sheaf under a smoothing of $X$ with singularity degrees $\ell_e$.

Besides the above description, $Y^{a,b}_{\ell,\m}$ is also described
in \cite{AE2}, Subsection~4.2, by giving the equations of each of its
irreducible components in the
corresponding component of $\mathbf R$, and in \cite{AE2}, Thm.~5.3, 
by means of its orbits under the action of $\mathbf
G_{\mathbf m}^V/\Gm$. These characterizations are the ones we use in the
present article. 

But we have also described $Y^{a,b}_{\ell,\m}$ in \cite{AE2},
Thm.~6.3, globally by its (infinitely many)
equations in $\mathbf R$, and in \cite{AE2}, Thm.~7.2, as an equivariant 
degeneration of the torus $\mathbf G_{\mathbf m}^V/\Gm$. The first
description will be important in understanding the moduli of the
$Y^{a,b}_{\ell,\m}$, to be worked out later. And the second yields that the connected
0-dimensional substack $\mathfrak I$ of $\mathbf J^d$ parameterizing
a collection of stable limits is the degeneration of a point! 

\subsection{Future work.} To attain our goal of giving a new definition of limit linear series
and constructing their moduli space, our series of articles will
continue. We summarize now what comes ahead.

First of all, we have constructed $\mathbf J^d$ as the quotient
expressed in \eqref{stackJ}, we have shown it parameterizes
torsion-free, rank-one, degree-$d$ sheaves on $X$, 
but we have not proved it is actually the moduli stack of these
sheaves. That $\mathbf J^d$ is bound to be the moduli stack is an
observation by Margarida Melo and Filippo Viviani, and the proof that
it actually is will appear later.

Second, given a smoothing $\pi\:\X\to B$ of $X$ and an invertible sheaf
$L_\eta$ on its generic fiber, we may consider the punctured arc in
the relative Artin stack $\mathbf J_{\X/B}$ parameterizing torsion-free
rank-one sheaves on the fibers of $\pi$. It is natural to think that
by adding the Artin substack $\mathfrak I$ of stable limits of
$L_\eta$ that the punctured arc will be completed to a $B$-flat closed
substack $A\subseteq\mathbf J_{\X/B}$. To show this, we need to construct $A$,
the degeneration of a point to $\mathfrak I$, whose existence we
proved in \cite{AE2}, Thm.~7.2, and moreover a relative torsion-free, rank-one sheaf on
$\X\times_BA$ whose restriction to the generic fiber is $L_\eta$. This
will also appear later.

Third, we will describe the moduli of all collections of
stable limits, which is tantamount to describing the various ways the 
$Y^{a,b}_{\ell,\m}$ appear as subschemes of fibers of 
$\mathbf R_{\mathbf J^{\mathbf b}}$ over
$\mathbf J^{\mathbf b}$. Already mentioned above, the global equations of
$Y^{a,b}_{\ell,\m}$ in $\mathbf R$, laid out by \cite{AE2}, Thm.~6.3, 
and a thorough description of how
$Y^{a,b}_{\ell,\m}$ actually depends on $a$, $b$, $\ell$ and $\m$ are
the fundamental pieces in this construction. However, since the
$Y^{a,b}_{\ell,\m}$ have infinitely many components, the moduli will
be formal. To get an actual moduli, we need to truncate $\mathbf
R_{\mathbf J^{\mathbf b}}$. The most obvious truncation is to restrict
to the open subscheme parameterizing sheaves obtained as gluings of
sheaves of nonnegative degrees on the components $X_v$, as for
understanding limits of effective divisors these sheaves are enough.

Fourth, the moduli of collections of stable limits can be thought of as a new
compactification of the Picard scheme of the curve $X$. To show it is
natural we will construct a relative version of it over the moduli
stack of stable curves $\comp{\mg_g}$. The construction may require a local analysis
and blowups of $\comp{\mg_g}$, in the way carried out by Main\`o in her
construction of the moduli space of enriched curves \cite{Maino}. But
further blowups will be necessary, infinitely many to get the formal
moduli, but finitely many for the truncated moduli. 

Finally, as mentioned before, our notion of limit
linear series, and the construction of their moduli space will be
similar to what is done for smooth curves: If $\mathfrak I$ is a 
closed substack of $\mathbf J^d$ parameterizing a collection of stable
limits, a ``limit linear series'' of ``sections'' of this collection
is a locally free subsheaf
$\mathcal V$ of the pushforward $p_{2*}\mathcal I$ satisfying certain
conditions, where $\mathcal
I$ is the sheaf induced on $X\times\mathfrak I$ (or certain base
extensions) by the universal
sheaf over $\mathbf J^d$.

\subsection{Organization} The layout of the paper is as follows. In Section~\ref{embsheaves} we construct the quotient stack $\mathbf J^d$
describing thoroughly its atlas and the group action. In
Section~\ref{degbdles} we consider smoothings of a nodal curve and
describe limits of line bundles. Limits of the trivial bundle are
considered in detail, and are explained in terms of the versal
deformation space of the curve, expanding on work done in \cite{Maino}
and \cite{EM}. In Section~\ref{stackydeg} we prove one of our main
theorems, Theorem~\ref{thm:main4}, describing collections of stable
limits as quotients by 
$\mathbf G_{\mathbf m}^V/\Gm$ of the
$Y^{a,b}_{\ell,\m}$. Finally, in
Section~\ref{sec:regeneration} we show that any such a quotient is a
collection of stable limits, our Theorem~\ref{regeneration}.

\subsection{Basic notations} In addition to the notations already introduced, throughout the
present article, we 
will denote by $N_e$ the node of $X$ associated to each $e\in E$ and
by $X_v$ the irreducible component of $X$ associated to each $v\in
V$. 

Given a collection $A\subseteq E$ of edges,
an orientation of $A$ is simply a section $\mathfrak o\:A\to\mathbb E$
over $A$ of the forgetful map $\E\to E$. We denote by $A^{\mathfrak
  o}$ the image of the orientation. To simplify the presentation, we fix an orientation 
$\mathfrak o\: E\to\mathbb E$ for $G$. Given $e\in\E$, we will write
$e$ as well for the (non-oriented) edge in $E$. Given $e\in E$, we denote
$e^{\mathfrak o}:=\mathfrak o(e)$. 

We will drop the subscripts
from $d_A$ and $d^*_A$ when appropriate. Given $\alpha$ in
$C^0(G,\Z)$ (resp.~$C^1(G,\Z)$), we write $\alpha_x:=\alpha(x)$ for
each vertex (resp.~oriented edge) $x$ of $G$. Similarly,
given a character $a$ of $C^0(G,\Z)$ (resp.~$C^1(G,\Z)$), we write
$a_v:=a(\chi_{\indm v})$ for each $v\in V$ (resp.~$a_e:=a(\chi_{\indm e}-\chi_{\indmbar{\ol
  e}})$ for each $e\in\E$). The multiplicative group is denoted $\Gm$.

\section{The space of embedded sheaves}\label{embsheaves}

\subsection{Gluing sheaves.} For each $v\in V$ and each $d\in\Z$, 
let $\IJ^d_v$ denote the degree-$d$ 
compactified Jacobian of $X_v$; it is a projective variety of 
dimension $g_v$,
where $g_v$ is the (arithmetic) genus of $X_v$, parameterizing
torsion-free, rank-one sheaves on $X_v$ of degree $d$. We refer to
\cite{Esteves01} 
for basic definitions and results concerning torsion-free rank-one sheaves. 

Fix an integer $d$ and let $\bfb\in C^0(G,\Z)$ be of degree 
$d$, meaning $\sum_{v\in V} \bfb_v=d$. Define
$$
\IJ^{\bfb}:=\prod_{v\in V}\IJ^{\bfb_v}_v.
$$
Note that $\IJ^{\bfb}$ is a projective variety of dimension $\sum g_v$, and comes with natural projection maps to $\IJ^{\bfb_v}_v$ for each $v\in V$.

For each vertex $v\in V$, let $\iota_v\:X_v\hookrightarrow X$ denote the
inclusion map, and $\L_v$ the pullback to $X_v\times\IJ^{\bfb}$ of a
``Poincar\'e sheaf,'' or universal sheaf on $X_v\times\IJ^{\bfb_v}_v$. 
Observe that $\L_v$ is actually not well-defined, in the
sense that a Poincar\'e sheaf is not unique. However, if we pick a 
point $P_v$ on the smooth locus of $X_v$, then 
$\L_v$ may be defined by imposing the condition that it be rigid at $P_v$, that
is, $\L_v|_{P_v\times\IJ^{\bfb}}\cong\O_{P_v\times\IJ^{\bfb}}$. At any rate, this is
not relevant to us, so we just choose one Poincar\'e sheaf
$\L_v$ for each vertex $v\in V$. For a point $s$ on $\IJ^{\bfb}$, we
denote by $\L_v(s)$ the restriction of $\L_v$ to $X_v\cong X_v \times s$.

For each oriented edge $e=uv\in\mathbb E$, let 
$$
\F_e:=\L_u|_{N_e\times\IJ^{\bfb}}\oplus\L_v|_{N_e\times\IJ^{\bfb}}.
$$
It is a rank-two vector bundle over $\IJ^{\bfb}$, under the natural isomorphism 
$N_e\times \IJ^{\bfb}\to \IJ^{\bfb}$. Let $\IP_{\IJ^{\bfb}}(\F_e)$ be the corresponding $\IP^1$-bundle over $\IJ^{\bfb}$, and define
$$
\IP_{\IJ^{\bfb}}:=\prod_{e\in E^{\mathfrak o}}\IP_{\IJ^{\bfb}}(\F_e),
$$
the product fibered over $\IJ^{\bfb}$. 
It is a $\IP_0$-bundle over 
$\IJ^{\bfb}$, where
$$
\IP_0:=\prod_{e\in E^\o}\IP^1.
$$

Fix now an element $\bfc\in C^1(G,\Z)$. 
For each $v\in V$, let $\wt\L_v$ denote the pullback of 
$(\iota_{v}\times 1_{\IJ^{\bfb}})_*\L_v$ to $X\times
\IP_{\IJ^{\bfb}}$. 
We will also denote by 
$\wt\F_e$ the pullback of $\F_e$ to $N_e\times\IP_{\IJ^{\bfb}}$ and 
by 
$$
q^{\bfc}_e\: \wt\F_e\lra \M_e
$$
the pullback of the universal quotient over 
$\IP_{\IJ^{\bfb}}(\F_e)$ to 
$N_e\times \IP_{\IJ^{\bfb}}$ 
for each $e\in E^{\mathfrak o}$. Notice that 
$$
\wt\F_e=\wt\L_{\te_e}|_{N_e\times\IP_{\IJ^{\bfb}}}\oplus
\wt\L_{\he_e}|_{N_e\times\IP_{\IJ^{\bfb}}}\quad\text{for each }e\in
E^{\mathfrak o}.
$$
We use the quotients $q^{\bfc}_e$ to construct a 
natural relative 
torsion-free, rank-one, degree-$d$ 
sheaf $\I^{\bfc}$ on  $X\times\IP_{\IJ^{\bfd}}/\IP_{\IJ^{\bfd}}$, the kernel of the
composition of surjections:
$$
\bigoplus_{v\in V}\wt\L_v \longrightarrow
\bigoplus_{v\in V}\bigoplus_{\substack{e\in E^{\mathfrak o}\\ e\ni v}}
\wt\L_v|_{N_e\times\IP_{\IJ^{\bfd}}} = 
\bigoplus_{e\in E^{\mathfrak o}}\Big(\wt\L_{\te_e}|_{N_e\times\IP_{\IJ^{\bfb}}}\oplus
\wt\L_{\he_e}|_{N_e\times\IP_{\IJ^{\bfb}}}\Big) 
\longrightarrow \bigoplus_{e\in E^{\mathfrak o}}\M_e.
$$
For a point $t$ on 
$\IP_{\IJ^{\bfd}}$, denote by $\I^{\bfc}(t)$ the induced 
torsion-free rank-one sheaf on $X \simeq X\times t$.

For each $c\in C^1(G,\Z)$, we may modify the above construction as follows: Let
$$
\L^c_v:=\L_v\Big(\sum_{\substack{e\in\mathbb E\\ \he_e=v}}(\bfc_e-c_e)N_e\times\IJ^{\bfb}\Big),
$$
and do the same construction as above but 
replacing the sheaves $\L_v$ with the sheaves $\L^c_v$. More precisely, 
for each $e=uv\in\E$, put
$$
\F^{c_e}_e:=\L^c_u|_{N_e\times\IJ^{\bfb}}\oplus
\L^c_v|_{N_e\times\IJ^{\bfb}},
$$
and set
$$
\IP^c_{\IJ^{\bfb}}:=\prod_{e\in E^{\mathfrak o}}\IP_{\IJ^{\bfb}}(\F^{c_e}_e),
$$
the fibered product over $\IJ^{\bfb}$. Instead of $\I^{\bfc}$, let $\I^{c}$ be the relative
torsion-free, rank-one, degree-$d$ sheaf on
$X\times\IP^c_{\IJ^{\bfb}}/\IP^c_{\IJ^{\bfb}}$ obtained 
from the pullbacks $q^c_e$ to $N_e\times\IP^c_{\IJ^{\bfb}}$ of the universal quotients on
$\IP_{\IJ^{\bfb}}(\F^{c_e}_e)$,  instead of the $q^{\bfc}_e$. As above, for a point $t$ on 
$\IP^{c}_{\IJ^{\bfd}}$, we denote by $\I^{c}(t)$ the induced
torsion-free rank-one sheaf on $X$. We will also say that a
torsion-free, rank-one sheaf $I$ is \emph{represented} by
$t\in\IP^{c}_{\IJ^{\bfd}}$ if $\I^{c}(t)\cong I$. 

Observe that for $t\in\IP^{c}_{\IJ^{\bfd}}$ lying over
$s\in\IJ^{\bfb}$, and for $e=uv\in E^{\mathfrak o}$, 
if $\I^c(t)$ is invertible at $N_e$, then it 
generates both $\L^c_u(s)$ and $\L^c_v(s)$ in a neighborhood
of $N_e$.  On the other hand, if $\I^c(t)$
fails to be invertible at $N_e$, then
either $q^c_e(t)(\L^c_u(s)|_{N_e})=0$ or $q^c_e(t)(\L^c_v(s)|_{N_e})=0$. In
the first case, in a neighborhood of $N_e$, 
the subsheaf of $\L^c_v(s)$ generated by $\I^c(t)$ is
$\L^c_v(s)(-N_e)$, whereas that of $\L^c_u(s)$ is $\L^c_u(s)$ itself. In the second case,
the reverse is true, that is, the same statement holds with $u$ and
$v$ exchanged. Letting $E_t$ 
denote the set of edges $e\in E$ for which $\I^c(t)$ fails to be invertible at
$N_e$, we obtain an orientation $\mathfrak o^c_t\: E_t\to\E$ by
assigning to each $e\in E_t$ the oriented edge whose head $v$ is such that 
$\I^c(t)$ generates $\L^c_v(s)(-N_e)$ in a neighborhood of $N_e$. Thus
$\I^c(t)$ generates the subsheaf
$$
\L^c_v(s)\Big(-\sum_{\substack{e\in E_t^c\\
    \he_e=v}}N_e\Big)
$$
for each $v\in V$, where, for simplicity, we put $E^c_t:=E_t^{\mathfrak o^c_t}$.

\subsection{The atlas.}\label{atlas}
We shall view the sheaves $\I^c$ as restrictions of a sheaf defined over
a larger base, containing all the schemes $\IP^{c}_{\IJ^{\bfb}}$ for all 
$c\in C^1(G,\Z)$ as closed subschemes. This works as follows.  For each vertex $v\in V$, and each natural number $n\in \mathbb N$, define the sheaf $\L^{(n)}_v$ as 
$$
\L^{(n)}_v:=\L_v\Big(n\sum_{\substack{e\in E^\o\\ e\ni v}}N_e\times\IJ^{\bfb}\Big).
$$
 By means of  the natural embeddings $\O_{X_v}\hookrightarrow\O_{X_v}(N_e)$, we may
view all the sheaves $\L^{c}_v$, for bounded $c\in C^1(G, \mathbb Z)$, more precisely for those $c$ verifying 
$|c_e-\bfc_e|\leq n$ for every $e\in\mathbb E$, as subsheaves of the sheaf  $\L^{(n)}_v$. We may thus view the schemes $\IP^{c}_{\IJ^{\bfd}}$ as
closed subschemes of the relative Quot-scheme 
$$
\Quot_{X\times\IJ^{\bfd}/\IJ^{\bfd}}\Bigl(\,
\bigoplus_{v\in V}(\iota_{v}\times 1_{\IJ^{\bfb}})_*\L^{(n)}_v\,\Bigr),
$$
more precisely, 
of the piece of the relative Quot-scheme parameterizing subsheaves of
rank one and degree $d$ or, equivalently, quotients of finite length
equal to $(2n+1)|E|$. 

Given two orientations
$\mathfrak o_1$ and $\mathfrak o_2$ of a subset $A\subseteq E$, define 
$(\mathfrak o_2,\mathfrak o_1)\in C^1(G,\Z)$ by
$$
(\mathfrak o_2,\mathfrak o_1)_e:=
\begin{cases}
1&\text{if $e\in A^{\mathfrak o_2}-A^{\mathfrak o_1}$,}\\
-1&\text{if $e\in A^{\mathfrak o_1}-A^{\mathfrak o_2}$,}\\
0&\text{otherwise.}
\end{cases}
$$

\begin{prop}\label{cappc} Notations as above, let $c,c'\in
  C^1(G,\mathbb Z)$ 
with $|c_e-\bfc_e| \leq n $ and $|c'_e-\bfc_e| \leq n$ for all $e\in \E$.
Viewing $\IP^{c}_{\IJ^{\bfb}}$
and $\IP^{c'}_{\IJ^{\bfb}}$ in the Quot-scheme, the following statements are true:
\begin{enumerate}
\item $\IP^{c}_{\IJ^{\bfb}}$
intersects $\IP^{c'}_{\IJ^{\bfb}}$ if and only if 
$|c'_e-c_e|\leq 1$ for every $e\in\mathbb E$. 
\item More precisely, given 
$t\in \IP^{c}_{\IJ^{\bfb}}$, we have that $t\in\IP^{c'}_{\IJ^{\bfb}}$
if and only if $c-c'=(\mathfrak o',\mathfrak o^c_t)$ where $\mathfrak
o'\:E_t\to\E$ is an orientation of the set $E_t$ of edges $e$ for
which $\I^c(t)$ fails to be invertible at $N_e$. In this case,
$\mathfrak o'=\mathfrak o_t^{c'}$.
\end{enumerate}
\end{prop}

\begin{proof} The first statement follows from the second. Indeed, if 
$t\in\IP^{c}_{\IJ^{\bfb}}\cap\IP^{c'}_{\IJ^{\bfb}}$, then Statement~(2)
yields that $|c'_e-c_e|=|(\mathfrak o',\mathfrak o^c_t)_e|\leq 1$ for
each $e\in\mathbb E$. Conversely, if $|c'_e-c_e|\leq 1$ for every
$e\in\mathbb E$, we
  may choose a point $t\in\IP^{c}_{\IJ^{\bfb}}$ such that $E_t$ is
  the subset of edges $e\in E$ in the support of $c-c'$ and
  $E_t^c$ is the collection of oriented edges $e$ such
  that $c'_e-c_e=1$. If $\mathfrak o'$ is the ``opposite
  orientation,'' that is, $E_t^{\mathfrak o'}$ is the collection of oriented edges $e$ such
  that $c_e-c'_e=1$, then $c-c'=(\mathfrak o',\mathfrak o_t)$. Thus, 
$t\in\IP^{c'}_{\IJ^{\bfb}}$ by Statement~(2).

As for Statement~(2), observe first that 
$\I^c(t)$ is a subsheaf of  $\bigoplus\iota_{v*}\L_v^{(n)}(s)$ such
that
$$
\I^c(t)\subseteq\bigoplus_{v\in
  V}\iota_{v*}\L^c_v(s)\Big(-\sum_{\substack{e\in E_t^c\\\he_e=v}}N_e\Big),
$$
where $s\in\IJ^{\bfb}$ lies under $t$. Thus, for each other orientation
$\mathfrak o'$ of $E_t$, 
we have
$$
\I^c(t)\subseteq\bigoplus_{v\in
  V} \iota_{v*}\L^c_v(s)\Big(\sum_{\substack{e\in E_t^{\mathfrak
      o'}\\\he_e=v}}N_e -\sum_{\substack{e\in E_t^c\\\he_e=v}}N_e\Big)=\bigoplus_{v\in
  V} \iota_{v*}\L^{c'}_v(s)
$$
if $c'=c-(\mathfrak o',\mathfrak o^{c}_t)$, and thus 
$\I^c(t)=\I^{c'}(t')$ as subsheaves of
$\bigoplus\iota_{v*}\L_v^{(n)}(s)$ for a certain
$t'\in\IP^{c'}_{\IJ^{\bfb}}$ such that
$E_{t'}=E_t$ and $\mathfrak o^{c'}_{t'}=\mathfrak o'$. So $t=t'$ in the
Quot-scheme, and thus $t\in\IP^{c'}_{\IJ^{\bfb}}$. 

Conversely, if $t=t'\in\IP^{c'}_{\IJ^{\bfb}}$ as well, then
\begin{equation}\label{ctc}
\I^c(t)\subseteq\bigoplus_{v\in
  V}\iota_{v*}\L^c_v(s)\Big(-\sum_{\substack{e\in A\\\he_e=v}}
(c'_e-c_e)N_e\Big)
\end{equation}
where $A:=\{e\in\E\,|\, c'_e>c_e\}$. Since $\I^c(t)$ generates
$\L^c_u(s)$ and $\L^c_v(s)(-N_e)$ in a neighborhood of $N_e$ for each
$e=uv\in E_t^c$, 
it follows that $|c'_e-c_e|\leq 1$ for every $e\in\mathbb
E$. Also, $c'_e-c_e=1$ only if $e\in E_t^c$,
that is $A\subseteq E_t^c$. 

Similarly, 
\begin{equation}\label{tct}
\I^{c'}(t')\subseteq\bigoplus_{v\in
  V}\iota_{v*}\L^{c'}_v(s)\Big(-\sum_{\substack{e\in A'\\\he_e=v}}
(c_e-c'_e)N_e\Big)
\end{equation}
where $A':=\{e\in\E\,|\, c_e>c'_e\}$. Then $A'\subseteq
E_{t'}^{c'}$. 

Note that $E_t=E_{t'}$, since $\I^c(t)=\I^{c'}(t')$. Thus $\mathfrak o^{c'}_{t'}$ is
another orientation for $E_t$. Also, if $e\in A$ then $\ol e\in A'$,
and hence 
$e\in E_t^c$ and $\ol e\in E_t^{c'}$. Thus $A\subseteq E_t^c-E_{t'}^{c'}$. 

On the other hand, let $e\in\E$ such that $e\in E_t^c-E_{t'}^{c'}$. As subsheaves of
$\bigoplus\iota_{v*}\L_v^{(n)}(s)$, the two sheaves in \eqref{ctc} are
equal to the corresponding ones in \eqref{tct}. Suppose by contradiction that $e\not\in
A$. Then $c'_e\leq c_e$. 
Since $e\in E_t^c$, we have that $\I^c(t)$ is contained in
$$
\bigoplus_{v\in
  V-\{\he_e\}}\iota_{v*}\L^c_v(s)\Big(-\sum_{\substack{f\in A\\\he_f=v}}
(c'_f-c_f)N_f\Big)\bigoplus
\iota_{\he_e*}\L^c_{\he_e}(s)\Big(-\sum_{\substack{f\in A\\\he_f=\he_e}}
(c'_f-c_f)N_f-N_e\Big).
$$
But, because of the equality of \eqref{ctc} and \eqref{tct}, also
$\I^{c'}(t')$ is contained in 
$$
\bigoplus_{v\in
  V-\{\he_e\}}\iota_{v*}\L^{c'}_v(s)\Big(-\sum_{\substack{f\in A'\\\he_f=v}}
(c_f-c'_f)N_f\Big)\bigoplus\iota_{\he_e*}\L^{c'}_{\he_e}(s)
\Big(-\sum_{\substack{f\in A'\\\he_f=\he_e}}
(c_f-c'_f)N_f-N_e\Big).
$$
But then $e\in E_{t'}^{c'}$, an absurd. 

It follows that  $A=E_t^c-E_{t'}^{c'}$. Then $c-c'=(\mathfrak o^{c'}_{t'},\mathfrak o^c_t)$.
\end{proof}

It follows from Proposition~\ref{cappc} that we may let 
$n$ tend to $\infty$, and consider the union of the 
$\IP^{c}_{\IJ^{\bfb}}$ for all $c\in C^1(G,\Z)$. 
We will denote this union by 
$$
\IR_{\IJ^{\bfb}}:=\bigcup_{c\in C^1(G,\Z)}\IP^{c}_{\IJ^{\bfb}}.
$$
It is a scheme locally of finite type over 
$\IJ^{\bfb}$. In fact, there is another way of describing
$\IR_{\IJ^{\bfb}}$, which shows that $\IR_{\IJ^{\bfb}}$ is naturally
fibered over $\IJ^{\bfb}$ with fibers equal to $\IR$, for the scheme
$\IR$ defined in \cite{AE2}, Subsection~3.2, and recalled below.

More precisely, for each
$e=uv\in E^{\mathfrak o}$ and 
$i\in\Z$, let
$$
\F_{e}^i:=\L_u(-(\bfc_e-i)N_e\times\IJ^{\bfb})|_{N_e\times\IJ^{\bfb}}\oplus
\L_v((\bfc_e-i)N_e\times\IJ^{\bfb})|_{N_e\times\IJ^{\bfb}}.
$$
Note that this definition is compatible with that of
$\F^{c_e}_e$, given previously, as the two sheaves coincide if 
$i=c_e$. 

As before, we may view the $\IP_{\IJ^{\bfb}}(\F_{e}^{i})$ for  integers $i$ with $|i-\bfc_e|\leq n$ 
as closed subschemes of the relative Quot-scheme
$$
\Quot_{X\times\IJ^{\bfb}/\IJ^{\bfb}}\Bigl(\,(\iota_{u}\times 1_{\IJ^{\bfb}})_*\big(\L_u\bigl(nN_e\times\IJ^{\bfb}\bigr)\big)\oplus
(\iota_{v}\times 1_{\IJ^{\bfb}})_*\big(\L_v\bigl(nN_e\times\IJ^{\bfb}\bigr)\big)\,\Bigr),
$$
more precisely, of the component of the relative Quot-scheme 
parameterizing quotients of finite length equal to $2n+1$.
Letting $n$ tend to $\infty$, we may
consider the union of all the $\IP_{\IJ^{\bfb}}(\F_{e}^{i})$ for all $i\in\Z$. We
will denote this union by 
$$
\IR^e_{\IJ^{\bfb}}:=\bigcup_{i\in\Z}\IP_{\IJ^{\bfb}}(\F_{e}^{i}).
$$
It is an $\IR_e$-bundle over $\IJ^{\bfd}$, 
where $\IR_e$ is the doubly infinite
chain of  smooth rational curves over $\k$, as in \cite{AE2},
Subsection~3.2. In addition, $\IR_{\IJ^{\bfb}}$ can be naturally
identified with the fibered product of
the $\IR^{e}_{\IJ^{\bfb}}$ for all $e\in E^{\mathfrak o}$ over $\IJ^{\bfb}$, that is,
$$
\IR_{\IJ^{\bfb}}=\prod_{e\in E^{\mathfrak o}}\IR^{e}_{\IJ^{\bfb}}.
$$
It is an $\IR$-bundle over $\IJ^{\bfb}$, where $\IR:=\prod_{e\in E^{\mathfrak o}}\IR_e$.

Yet more precisely, we identify a fiber of $\IR_{\IJ^{\bfb}}/\IJ^{\bfb}$
with $\IR$ in the following way: A point 
$s\in\IJ^{\bfb}$ corresponds to a collection of torsion-free, rank-one sheaves
$(L_v\,;\, v\in V)$ with $L_v$ of degree $\bfb_v$. We fix trivializations
$L_v|_{N_e}\cong \k$ for each $e\in E$ and each $v\in e$. We fix
as well trivializations $\O_{X_v}(N_e)|_{N_e}\cong \k$ for each $e\in
E$ and each $v\in e$, and consider the induced trivializations
$\O_{X_v}(mN_e+D)|_{N_e}\cong \k$ for each $m\in\Z$ and each Cartier
divisor $D$ of $X_v$ not containing $N_e$ in its support. These
trivializations induce trivializations 
$$
L_v\Big(\sum_{\substack{e\in\mathbb E\\ \he_e=v}}(\bfc_e-c_e)N_e\Big)|_{N_f}\cong \k
$$
for each $f\in E$, each $v\in f$ and each $c\in C^1(G,\Z)$. These
trivializations give rise to an isomorphism between the fiber of
$\IR^E_{\IJ^{\bfb}}/\IJ^{\bfb}$ over $s$ and $\IR$. 

The universal subsheaf on 
$$
X\times\Quot_{X\times\IJ^{\bfb}/\IJ^{\bfb}}\Bigl(\,
\bigoplus_{v\in V}\bigl(\iota_{v}\times
1_{\IJ^{\bfb}}\bigr)_*\L^{(n)}_v\Bigr)
$$
for 
$n\to\infty$ restricts to a relative torsion-free, rank-one, degree-$d$ sheaf on 
$X\times\IR_{\IJ^{\bfb}}/\IR_{\IJ^{\bfb}}$, which shall be denoted by $\I$. Its
restriction to the subscheme $X\times \IP^{c}_{\IJ^{\bfb}}$ is $\I^{c}$ for each $c\in C^1(G,\Z)$.

\medskip

For each $c\in C^1(G,\Z)$, an open dense subscheme of $\IP^c_{\IJ^{\bfb}}$ 
parameterizes invertible
sheaves of multidegree $\bfb+d_\Z^*(\bfc-c)$. Since $\deg$ is a
cokernel for $d^*_\Z$, 
it follows that, up to translation, we could
have changed $\bfc$ for any other 1-cochain and $\bfb$ for any other
0-cochain of degree $d$. Moreover, the orientation $\mathfrak o$
is just a convenient means of ordering the curves in the chain $\IR_e$
for each $e\in E$; a different orientation would lead to
the same fibration $\IR_{\IJ^{\bfb}}/\IJ^{\bfb}$. 

\medskip

Another way of interpreting Proposition~\ref{cappc} is through the
following definition and proposition: Given $t\in \IR_{\IJ^{\bfb}}$,
let $c\in C^1(G,\Z)$ such that $t\in \IP^c_{\IJ^{\bfb}}$. Define
$c(t)\in C^1(G,\Q)$ by:
$$
c(t)_e:=c_e+\begin{cases}
+\frac{1}{2}&\text{if $e\in E^c_t$,}\\
-\frac{1}{2}&\text{if $\ol e\in E^c_t$,}\\
0&\text{otherwise.}
\end{cases}
$$
Notice that $\E_t:=\{e\in\E\,|\,c(t)_e\not\in\Z\}$ is precisely the set
of oriented edges supported in $E_t$.

Define as well, for each $c\in C^1(G,\frac{1}{2}\Z)$,
$$
\mathcal L^c_v:=\mathcal L_v\Big(\sum_{\substack{e\in\E\\ \he_e=v}}\lfloor
\mathbf c_e-c_e\rfloor N_e\times\IJ^{\bfb}\Big)
$$
and
$$
\mathbf P^c_{\IJ^{\bfb}}:=\bigcap_{\substack{c'\in C^1(G,\Z)\\ 
|c'_e-c_e|\leq\frac{1}{2} \forall e\in\E}}\mathbf P^{c'}_{\IJ^{\bfb}}.
$$

\begin{prop}\label{cappc12} Let $t\in \IR_{\IJ^{\bfb}}$. Then the
  following statements are true:
\begin{enumerate}
\item $c(t)$ is well-defined.
\item $t\in \IP^{c}_{\IJ^{\bfb}}$ for $c\in C^1(G,\Z)$ if and only if
  $|c_e-c(t)_e|\leq\frac{1}{2}$ for every $e\in\E$. In particular, $t\in\IP^{c(t)}_{\IJ^{\bfb}}$.
\item The torsion-free sheaf $\I(t)_v$ generated by $\I(t)$ on $X_v$
  for each $v\in V$ is isomorphic to $\mathcal L^{c(t)}_v(s)$, where
  $s\in\IJ^{\bfb}$ is the point lying under $t$.
\end{enumerate}
\end{prop}

\begin{proof} Suppose $t$ lies on $\IP^{c}_{\IJ^{\bfb}}\cap
  \IP^{c'}_{\IJ^{\bfb}}$. It
  follows from Proposition~\ref{cappc} that $c'-c=(\mathfrak
  o^c_t,\mathfrak o^{c'}_t)$. Then
$$
c'_e:=c_e+\begin{cases}
+1&\text{if $e\in E^c_t-E^{c'}_t$}\\
-1&\text{if $e\in E^{c'}_t-E^{c}_t$}\\
0&\text{otherwise}
\end{cases}\,\,=\,\,
c(t)_e - \begin{cases}
-\frac{1}{2}&\text{if $\ol e\in E^{c'}_t$}\\
+\frac{1}{2}&\text{if $e\in E^{c'}_t$}\\
0&\text{otherwise,}
\end{cases}
$$
where we used the fact that $E^c_t$ and $E^{c'}_t$ are orientations of
the same set, $E_t$. Thus, the definition of $c(t)$ does not
change if $c'$ were chosen instead of $c$.

Furthermore, from the above argument it follows that
$|c_e-c(t)_e|\leq\frac{1}{2}$ for each $e\in\E$ and $c\in
C^1(G,\Z)$ such that $t\in \IP^{c}_{\IJ^{\bfb}}$. Conversely, let $c\in
C^1(G,\Z)$ such that $t\in \IP^{c}_{\IJ^{\bfb}}$. If $c'\in
C^1(G,\Z)$ is such that $|c'_e-c(t)_e|\leq\frac{1}{2}$ for every
$e\in\E$, then $c'_e=c(t)_e=c_e$ for each $e\in\E$ such that
$c(t)_e\in\Z$, that is, such that $e\not\in\E_t$. On the other hand, if $e\in\E_t$ then
$|c'_e-c(t)_e|=\frac{1}{2}$. Let $\mathfrak o'$ be the orientation of
$E_t$ such that $e\in E_t^{\mathfrak o'}$ if and only if
$c(t)_e-c'_e=\frac{1}{2}$. Then 
$c-c'=(\mathfrak o',\mathfrak o^c_t)$, and it follows from
  Proposition~\ref{cappc} that $t\in \IP^{c'}_{\IJ^{\bfb}}$.

Finally, let $s\in\IJ^{\bfb}$ be the point lying under $t$. Let $c\in
C^1(G,\Z)$ such that $t\in \IP^{c}_{\IJ^{\bfb}}$. Then $\mathcal I^c(t)$ generates the subsheaf
$$
\mathcal L^c_v(s)\Big(-\sum_{\substack{e\in E^c_t\\ \he_e=v}}N_e\Big)
$$
of $\mathcal L^c_v(s)$ for each $v\in V$. But this subsheaf is 
$\mathcal L^{c(t)}_v(s)$.
\end{proof}

\subsection{Group action.} There is a natural action on
$\IR_{\IJ^{\bfb}}$ by $H^1(G,\Z)$. 
Indeed, given $\gamma \in H^1(G,\Z)$ and a point $s$ of $\IJ^{\bfb}$, 
corresponding 
to a collection of torsion-free, rank-one sheaves $(L_v\,;\,v\in V)$,
we associate $s'\in\IJ^{\bfb}$,
corresponding to the tuple of sheaves $(L_v^{-\gamma}\,;\,v\in V)$, where
$$
L_v^{-\gamma}:=L_v(-\sum_{\substack{e\in\mathbb E\\ \he_e=v}}\gamma_eN_e).
$$
This association gives a map $\tau^\gamma\:\IJ^{\bfb}\to\IJ^{\bfb}$
which sends $s$ to $s'$. Clearly, $\gamma\mapsto\tau^{\gamma}$ gives a
group homomorphism $H^1(G,\Z)\to\text{Aut}(\IJ^{\bfb})$. Also, 
$$
\L^{c}_v(\tau^\gamma(s))=\L^{c+\gamma}_v(s)\quad\text{for each
}s\in\IJ^{\bfb},\, c\in C^1(G,\Z) \text{ and }v\in V.
$$
By
construction, for each point $s\in\IJ^{\bfd}$, the restrictions of $\I^{c+\gamma}$ over points on the 
fiber of $\IP^{c+\gamma}_{\IJ^{\bfb}}$ 
over $s$ 
are exactly the
same restrictions of $\I^{c}$ over points on the fiber of
$\IP^{c}_{\IJ^{\bfb}}$ over $\tau^\gamma(s)$. 
Thus, the translation $\tau^\gamma\:\IJ^{\bfb}\to\IJ^{\bfb}$ 
lifts to an 
automorphism $\wt\tau^\gamma\:\IR_{\IJ^{\bfb}}\to\IR_{\IJ^{\bfb}}$ sending 
$\IP^{c+\gamma}_{\IJ^{\bfb}}$ to $\IP^{c}_{\IJ^{\bfb}}$ for each $c\in
C^1(G,\Z)$. More precisely, given $t\in \IP^{c+\gamma}_{\IJ^{\bfb}}$ over
$s\in\IJ^{\bfb}$, corresponding to the subsheaf
$$
\I^{c+\gamma}(t)\subseteq\bigoplus_v\L^{c+\gamma}_v(s),
$$
$\tau^\gamma$ sends $s$ to the point $s'$ corresponding to the tuple
$(\L_v(s)^{-\gamma}\,;\,v\in V)$ and $\wt\tau^\gamma$
sends $t$ to the point on $\IP^{c}_{\IJ^{\bfb}}$ over $s'$ which 
corresponds to the 
same subsheaf $\I^{c+\gamma}(t)$. In other words,
$$
\I^{c}(\wt\tau^{\gamma}(t))=\I^{c+\gamma}(t)
\quad\text{as subsheaves of}\quad
\bigoplus_v\L^{c}_v(\tau^{\gamma}(s))=\bigoplus_v\L^{c+\gamma}_v(s).
$$

\begin{prop} Notations as above, the 
assignment $\gamma \mapsto\wt\tau^\gamma$ defines an injective group 
homomorphism $H^1(G,\Z)\to\text{\rm Aut}(\IR_{\IJ^{\bfb}})$. Moreover, 
$\wt\tau^\gamma$ has fixed points only if $\gamma=0$.
\end{prop}

\begin{proof}  By construction, the assignment $\wt\tau^{\,\sbullet} :
  H^1(G, \mathbb Z) \to \text{Aut}(\IR_{\IJ^{\bfb}})$ is clearly a
  group homomorphism. We prove the second statement, from which the
  injectivity of $\wt\tau^{\,\sbullet}$ follows.

Let $\gamma \in H^1(G, \mathbb Z)$, and suppose $\wt\tau^\gamma$ has a fixed point 
$t\in\IP^{c}_{\IJ^{\bfb}}$ for a certain $c\in C^1(G,\Z)$. Since $\wt\tau^\gamma$ sends 
$\IP^{c'+\gamma}_{\IJ^{\bfb}}$ to $\IP^{c'}_{\IJ^{\bfb}}$ for each
$c'\in C^1(G,\Z)$, it follows that $t$ lies on $\IP^{c-n\gamma}_{\IJ^{\bfb}}$ for
every integer $n\geq 0$. This is possible, by Proposition~\ref{cappc},
only if $\gamma=0$.
\end{proof}

The action of $H^1(G,\Z)$ on $\IR_{\IJ^{\bfb}}$ is thus free. We may
consider the quotient stack:
$$
\IS^d:=\Bigg[\frac{\IR_{\IJ^{\bfb}}}{H^1(G,\Z)}\Bigg].
$$

Since $\I(t)=\I(\wt\tau^\gamma(t))$ for each $t\in\IR^E_{\IJ^{\bfb}}$ and $\gamma\in H^1(G,\Z)$, the sheaf $\I$
descends to a torsion-free, rank-one, degree-$d$ sheaf on
$X\times\IS^d/\IS^d$, which we shall also denote by $\I$.

\smallskip

We call $\IS^d$ the \emph{moduli of embedded sheaves}. 
We have the following moduli description.

\begin{prop}\label{Sd}
The stack $\IS^d$ parameterizes the data of torsion-free, 
rank-one sheaves $M_v$ on $X_v$ for each $v\in V$, and subsheaves 
$I\subseteq\oplus M_v$ of degree $d$ such that all the induced maps 
$h_v\:I\to M_v$ are surjective. 
\end{prop}

\begin{proof}
For each such data, let $E_I$ be the collection of edges $e\in E$ for
which $I$ fails to be invertible at $N_e$. Let $\mathfrak
u\:E_I\to\E$ be an orientation. For each $v\in V$, let
$$ 
\widetilde M_v:=M_v\Big(\sum_{\substack{e\in E_I^{\mathfrak u}\\\he_e=v}} N_e\Big).
$$ 
Then 
$\sum_{v\in V}\deg(\widetilde M_v)=d$, whence there are $c\in C^1(G,\Z)$ and 
$s\in\IJ^{\bfd}$ such that $\widetilde M_v\cong\L_v^{c}(s)$ for each $v\in
V$. (Given $c$, the point $s$ is unique.) Now, since
$$
I\subseteq\bigoplus_{v\in V} M_v\subseteq\bigoplus_{v\in V}\L_v^{c}(s),
$$
it follows that the data corresponds to a unique point $t$ on 
$\IP^{c}_{\IJ^{\bfd}}$ over $s$. Since the $h_v$ are surjective,
$\mathfrak o^c_t=\mathfrak u$.

A choice $\mathfrak u'\:E_I\to\E$ different from $\mathfrak u$ would correspond to a
point $t'$ on $\IP^{c'}_{\IJ^{\bfd}}$ over the same $s$, for $c':=c-(\mathfrak u',\mathfrak u)$. 
Note however that $t'=t$ on $\IR_{\IJ^{\bfb}}$ by Proposition~\ref{cappc}.

Furthermore, for a fixed $\mathfrak u$, any other
choice of $c$ differs from the one above by an element of
$H^1(G,\Z)$. More precisely, if $c'\in C^1(G,\Z)$ and
$s'\in\IJ^{\bfb}$ are such that $\widetilde M_v\cong\L^{c'}_v(s')$ for
every $v\in V$, then $\gamma:=c-c'\in H^1(G,\Z)$ and 
$s'=\tau_{\gamma}(s)$. As the construction yields the same subsheaf $I$,
we have that the corresponding point $t'$ on $\IP^{c'}_{\IJ^{\bfd}}$
over $s'$ satisfies $t'=\wt\tau_{\gamma}(t)$. Then the images
of $t$ and $t'$ on $\IS^d$ are the same.

Conversely, given a point on $\IS^d$, let 
$t\in\IR_{\IJ^{\bfb}}$ be a lifting, and 
$s\in\IJ^{\bfb}$ its image. Let $c$ be an element of $C^1(G,\Z)$ such that
$t\in\mathbf P^c_{\IJ^{\bfb}}$. Then $s$ corresponds to 
the sheaves $\L_v(s)$ for $v\in V$ and $t$ to the subsheaf
\begin{equation}\label{ict}
\I^c(t)\subseteq\bigoplus_{v\in V}\L_v^{c}(s).
\end{equation}
For each $v\in V$, the sheaf $\I^c(t)$ generates the subsheaf 
$$
M_v:=\L_v^c(s)(-\sum_{\substack{e\in E_t^c\\\he_e=v}} N_e),
$$
Then
$$
I:=\I^c(t)\subseteq \bigoplus_{v\in V}M_v
$$
with surjective induced maps $I\to M_v$.

Picking a different $c$ will not change the data of the $M_v$ and
$I\subseteq\oplus M_v$, since replacing $\oplus\L_v^{c}(s)$ by the
larger sheaf $\oplus_v\L_v^{(n)}(s)$ for large $n$ does not change the sheaves
$M_v$ obtained.

If $t'\in\IR_{\IJ^{\bfb}}$ is another lifting, lying over
$s'\in\IJ^{\bfb}$, then $s'=\tau_{\gamma}(s)$ and
$t'=\wt\tau_{\gamma}(t)$ for a certain $\gamma\in H^1(G,\Z)$.
Since $t\in\P^c_{\IJ^{\bfb}}$ we have $t'\in\P^{c-\gamma}_{\IJ^{\bfb}}$.
Also, the
inclusion in \eqref{ict} is the same as
$$
\I^{c-\gamma}(t')\subseteq\bigoplus_{v\in V}\L_v^{c-\gamma}(s').
$$
Thus we obtain the same data of sheaves $M_v$ on $X_v$ and subsheaf of
$\oplus M_v$.

We omit the simple verification that the two procedures described
above are inverse to each other.
\end{proof}

Observe that a surjection $h_v\:I\to M_v$ from a torsion-free, rank-one sheaf on $X$ to one on $X_v$ induces an isomorphism 
$\wt h_v\: I_v\to M_v$, where $I_v$ is the restriction of $I$ to $X_v$
modulo torsion. So we may view the data
defining a point on $\IS^d$ as 
the equivalence class of the 
data of an abstract torsion-free, rank-one degree-$d$ sheaf $I$ on $X$ together with 
elements $z_v\in\Gm(\k)$ for each $v\in V$; these data give rise to a map
$$
\begin{CD}
I\hooklongrightarrow\oplus I_v @>(z_v\,;\,v\in V)>> 
\oplus I_v,
\end{CD}
$$
and the equivalence identifies data with the same image.

\medskip

Another way to put this is as follows. 
Observe that there is a natural action of the character group of
$C^0(G,\Z)$, henceforth denoted $\mathbf G^V_{\mathbf m}$, on
$\IR_{\IJ^{\bfd}}$. In fact, the action of each $z\in\mathbf
G^{V}_{\mathbf m}(\k)$ 
is the restriction of the one on 
$\Quot_{X\times\IJ^{\bfd}/\IJ^{\bfd}}\Bigl(\,\bigoplus_{v\in V}(\iota_v\times
1_{\IJ^{\bfb}})_*\L^{(n)}_v\,\Bigr)$ induced by the
isomorphisms
$$
\L^{(n)}_v  \stackrel{z_v}{\longrightarrow} \L^{(n)}_v\text{ for $v\in V$}.
$$

If all the $z_v$ are equal, the action is the identity. We have thus
defined an action of $\mathbf G_{\mathbf m}^{V}/\Gm$, where $\Gm$ is
embedded diagonally in $\mathbf G_{\mathbf m}^{V}$, on 
$\IR_{\IJ^{\bfd}}$ that sends a fiber over $\IJ^{\bfb}$ to
itself. Moreover, the isomorphism between a fiber of
$\IR^E_{\IJ^{\bfb}}$ over $\IJ^{\bfb}$ and $\IR$ is
equivariant with respect to the action of $\mathbf G_{\mathbf
  m}^{V}/\Gm$ on $\IR$ defined in \cite{AE2}, Subsection 5.1, and recalled
below in Section~\ref{stackydeg}. 
Also, the action on $\IR_{\IJ^{\bfd}}$ restricts to an action on $\mathbf
P^c_{\IJ^{\bfb}}$ for each $c\in C^1(G,\Z)$. Finally, this action commutes with the action of $H^1(G,\Z)$ and thus induces
an action on $\IS^d$.

The (abstract) torsion-free, rank-one sheaves $I$ of degree $d$ on $X$ 
correspond to orbits of $\IS^d$ under $\mathbf G_\mathbf m^V/\Gm$. The
specific point on the orbit tells us how $I$ is viewed 
inside $\oplus I_v$. So, we may view the space of 
orbits, the quotient stack 
$$
\mathbf J^d:=\Bigg[\frac{\IS^d}{\mathbf G_\mathbf m^V/\Gm}\Bigg],
$$
as  a parameter space for 
torsion-free, rank-one sheaves of degree $d$ on $X$. The sheaf $\I$ on $X\times\IS^d$ is 
$\mathbf{G}_{\m}^V/\mathbf{G}_\m$-invariant, and thus descends to a
relative torsion-free, rank-one,
degree-$d$ sheaf on $X\times \mathbf J^d/\mathbf J^d$, a universal
torsion-free, rank-one, degree-$d$ sheaf.

However,  the quotient $\mathbf J^d$ is not 
well-behaved because the orbits of the action have variable dimension. 

\begin{prop} The orbits of maximum dimension of $\IS^d$ under the
  action of $\mathbf G_\mathbf m^V/\Gm$ have 
dimension $|V|-1$ and correspond to simple torsion-free, rank-one
sheaves of degree $d$.
\end{prop} 

\begin{proof} Simple, torsion-free, rank-one sheaves on $X$ are
  direct images of invertible sheaves on \emph{connected} partial
  desingularizations of $X$. So
  the proof reduces to considering invertible sheaves on $X$, and to
  showing that the subscheme of the Picard scheme of $X$ parameterizing 
invertible sheaves with given restrictions to the components $X_v$ is
isomorphic to $\mathbf G_{\mathbf m}^{|E|-|V|+1}$. This has been shown by a number
of authors --- see e.g.~\cite[10.2]{OS79} --- and follows
formally from the exact sequence 
\[
1 \rightarrow H^1(G, \mathbf G_{\mathbf m}) \rightarrow
\mathrm{Pic}^0(X) \rightarrow \prod_{v \in V} \mathrm{Pic}^0(X_v)
\rightarrow 1.
\]
\end{proof}

Thus, restricting to
orbits of dimension $|V|-1$, and considering the orbit space, 
we get the moduli space of
simple, torsion-free, rank-one sheaves of degree $d$ on $X$. This space
has been constructed by Altman and Kleiman as an algebraic
space~\cite{AK80,AK79}. (Actually, they prove the representability by an algebraic
space of a functor parameterizing far
more general objects.) The second author of the present article has shown later that
this space is actually a scheme locally of finite type and universally closed over the
field, though not separated; see~\cite{Esteves01}.

The second author has also considered certain open subspaces of the moduli
space of simple sheaves which are actually proper over the field,
actually projective \cite{Espine}, thus
producing various compactifications of the Jacobian, even for general families
of curves satisfying certain mild conditions \cite{Esteves01}. Instead of doing this,
however, we will consider clusters of orbits. The following
proposition, which collects statements proved along the section and
states new ones, 
will serve us in describing these clusters.
 
\begin{prop}\label{embshfc}  Let $I$ be a torsion-free, 
rank-one, degree-$d$ sheaf on
  $X$. Let $E_I$ be the collection of edges $e\in E$ for which $I$
  fails to be invertible at $N_e$. For each $v\in V$, let $I_v$ denote the
  maximum torsion-free quotient of $I|_{X_v}$. Then the following
  statements hold:
\begin{enumerate}
\item For each orientation $\mathfrak u\:E_I\to\E$, there are 
$c\in C^1(G,\Z)$ and $s\in\IJ^{\bfd}$ such that 
\begin{equation}\label{ocs}
\L^{c}_v(s)\cong I_v\Big(\sum_{\substack{e\in E_I^{\mathfrak u}\\ \he_e=v}}
N_e\Big)\quad\text{for each }v\in
V.
\end{equation}
\item For each $\mathfrak u$, $c$ and $s$ satisfying Equations~\eqref{ocs}, there is a
point $t\in \IP^{c}_{\IJ^{\bfd}}$ over
  $s$ representing $I$ with $\mathfrak o_t^c=\mathfrak u$.
  \smallskip
  
\item For each point $s\in\IJ^{\bfd}$ and each element $\mathfrak c\in
  C^1(G,\frac{1}{2}\Z)$ 
satisfying $e\in E_I$ if and only if $\mathfrak c_e\not\in\Z$, and such that 
\begin{equation}\label{ofrakcs}
\L^{\mathfrak c}_v(s)\cong I_v \quad\text{for each }v\in V,
\end{equation}
there is a point $t\in \mathbf P^{\mathfrak c}_{\IJ^{\bfd}}$ over $s$ representing $I$ with $c(t)=\mathfrak c$.
\smallskip

\item Conversely, given $c\in C^1(G,\Z)$ and $s\in\IJ^{\bfb}$, and
  a point $t$ on the fiber of $\IP^{c}_{\IJ^{\bfd}}$ over
  $s$ representing $I$, Equations~\eqref{ocs} hold for
  $\mathfrak u:=\mathfrak o_t^c$ and Equations~\eqref{ofrakcs} hold
  for $\mathfrak c:=c(t)$.
\smallskip

\item For each two sets of data $(\mathfrak u_1,c_1,s_1)$ and
  $(\mathfrak u_2,c_2,s_2)$ satisfying \eqref{ocs}, and each two points
  $t_1$ and $t_2$ representing $I$ on the fibers of $\IP^{c_1}_{\IJ^{\bfd}}$ and
  $\IP^{c_2}_{\IJ^{\bfd}}$ over $s_1$ and $s_2$, respectively,
  with $\mathfrak o^{c_1}_{t_1}=\mathfrak u_1$ and
  $\mathfrak o^{c_2}_{t_2}=\mathfrak u_2$, we have that
$$
\gamma:=c_1-c_2+(\mathfrak u_1,\mathfrak u_2)=c(t_1)-c(t_2)
$$
is in $H^1(G,\Z)$, satisfies $\tau^\gamma(s_1)=s_2$ and is such that
$t_2$ differs from $\wt\tau^{\gamma}(t_1)$ by the action of an element
of $\mathbf G_{\mathbf m}^V/\mathbf G_{\mathbf m}(\k)$.
\end{enumerate}
\end{prop} 

\begin{proof} For each orientation $\mathfrak u\:E_I\to\E$, 
let $f^{\mathfrak u}\in C^0(G,\Z)$ be defined by
$$
f^{\mathfrak u}(v):=\deg(I_v)+|E^{\mathfrak u}_I(v)| \text{ for each
 $v\in V$, where }E^{\mathfrak u}_I(v):=\{e\in E_I^{\mathfrak u}\,|\,
\he_e=v\}.
$$
Then \eqref{ocs} is satisfied for $\mathfrak u$, $c$ and $s$ only if 
\begin{equation}\label{cslb}
f^{\mathfrak u}= d^*(\bfc-c)+\bfb.
\end{equation}
 Furthermore, if the latter
  holds for given $\mathfrak u$ and $c\in C^1(G,\Z)$, then there is a unique
  $s\in\IJ^{\bfb}$ such that \eqref{ocs} holds for $\mathfrak
  u$, $c$ and $s$. Indeed, as the points
  $s\in\IJ^{\bfb}$ parameterize all tuples of sheaves with degrees given
  by $\bfb$, the tuples $(\L^{c}_v(s)\,|\,v\in V)$ run through all
  tuples of sheaves with degrees given by
  $d^*(\bfc-c)+\bfb$.

Now, given $\mathfrak u$, since the union of all $E^{\mathfrak u}_I(v)$ for $v\in V$ is
$E_I^{\mathfrak u}$, it follows that $f^{\mathfrak u}$ has degree $d$,
the same as $\bfb$. So there is $c\in C^1(G,\Z)$ such that 
$f^{\mathfrak u}=d^*(\bfc-c)+\bfb$. This finishes the proof of Statement~(1). 

Statement~(2)~and the first part of Statement~(4) have already been addressed before, in the proof of
Proposition~\ref{Sd}. The remaining of Statement (4) follows from the
definition of $c(t)$ for $t\in \IP^{c}_{\IJ^{\bfd}}$. 

As for Statement~(3), let $\mathfrak u\:E_I\to\E$ be any
orientation. And let $c\in C^1(G,\Z)$ such that $c_e=\mathfrak c_e$
for each $e$ such that $\mathfrak c_e\in\Z$ and $c_e=\lfloor\mathfrak c_e\rfloor$
for each $e\in E^{\mathfrak u}_I$. Then 
$$
\L^{c}_v(s)\cong\L^{\mathfrak c}_v(s)\Big(\sum_{\substack{e\in
    E_I^{\mathfrak u}\\\he_e=v}}N_e\Big)\quad\text{for each $v\in V$,}
$$
from which it follows that $\mathfrak u$, $c$ and $s$ satisfy
Equations~\eqref{ocs}. By Statement~(2), there is a
point $t\in \IP^{c}_{\IJ^{\bfd}}$ over
  $s$ representing $I$ with $\mathfrak o_t^c=\mathfrak u$. From the
  definition of $c(t)$, it follows that $c(t)=\mathfrak c$.

Finally, we prove Statement~(5). First observe that
$$
f^{\mathfrak u_1}(v)-f^{\mathfrak u_2}(v)=
|E^{\mathfrak u_1}_I(v)|-|E^{\mathfrak u_2}_I(v)|=
d^*((\mathfrak u_1,\mathfrak u_2))(v) \text{ for each }v\in V.
$$
Since \eqref{ocs} holds for the 
two sets of data, $(\mathfrak u_1,c_1,s_1)$ and
$(\mathfrak u_2,c_2,s_2)$, instead of $(\mathfrak u,c,s)$,
Equations~\eqref{cslb} hold, that is,
$$
f^{\mathfrak u_1}= d^*(\bfc-c_1)+\bfb
\quad\text{and}\quad 
f^{\mathfrak u_2}= d^*(\bfc-c_2)+\bfb.
$$
Thus
$$
d^*(c_2-c_1)=f^{\mathfrak u_1}-f^{\mathfrak
  u_2}=d^*((\mathfrak u_1,\mathfrak u_2)),
$$
and hence $\gamma$, as defined, is in $H^1(G,\Z)$. The equality 
$$
c_1-c_2+(\mathfrak u_1,\mathfrak u_2)=c(t_1)-c(t_2)
$$
follows from the definition of $c(t_1)$ and $c(t_2)$, 
using that $\mathfrak o_{t_1}^{c_1}=\mathfrak u_1$ and 
$\mathfrak o_{t_2}^{c_2}=\mathfrak u_2$.

Now, since \eqref{ocs} holds for $(\mathfrak u_1,c_1,s_1)$
and for $(\mathfrak u_2,c_2,s_2)$, instead of $(\mathfrak u,c,s)$, we have
\begin{align*}
\L^{c_2}_v(\tau^{\gamma}(s_1))\cong\L^{c_2+\gamma}_v(s_1)&=
\L^{c_1-(\mathfrak u_2,\mathfrak u_1)}_v(s_1)\\
&\cong
\L_v^{c_1}(s_1)\Big(\sum_{\substack{e\in\E\\\he_e=v}}(\mathfrak
u_2,\mathfrak u_1)_eN_e\Big)\\
&\cong I_v\Big(\sum_{e\in E_I^{\mathfrak
  u_1}(v)}N_e+
\sum_{e\in E_I^{\mathfrak
  u_2}(v)}N_e-\sum_{e\in E_I^{\mathfrak
  u_1}(v)}N_e\Big)\\
&\cong
I_v\Big(\sum_{e\in E_I^{\mathfrak u_2}(v)}N_e\Big)\cong\L_v^{c_2}(s_2)
\end{align*}
for each $v\in V$, 
whence $\tau^\gamma(s_1)=s_2$.

Also, $\wt\tau^{\gamma}(t_1)$ represents the same sheaf as
$t_1$, and $\wt\tau^{\gamma}(t_1)$ lies on $\mathbf
P^{c_1-\gamma}_{\IJ^{\bfb}}$, on the same fiber of $\mathbf R_{\IJ^{\bfb}}$ over
$\IJ^{\bfb}$ as $t_2$, that over $s_2$. We may thus assume $s_1=s_2$
and $\gamma=0$. 

Now, $\mathfrak o_{t_1}^{c_1}=\mathfrak u_1$ and $\mathfrak o_{t_2}^{c_2}=\mathfrak
u_2$. Since $c_2=c_1+(\mathfrak
u_1,\mathfrak u_2)$, it follows from Proposition~\ref{cappc} that $t_1\in \mathbf
P^{c_2}_{\IJ^{\bfb}}$ as well. We may thus assume $c_1=c_2$, whence 
$\mathfrak u_1=\mathfrak u_2$. Set $c:=c_1$.

Since $t_1$ and $t_2$ represent $I$, there is an isomorphism
$\I^c(t_1)\to\I^c(t_2)$. This isomorphism induces isomorphisms 
$\I^c(t_1)_v\to\I^c(t_2)_v$ for each $v\in V$ making the diagram
commute:
$$
\begin{CD}
\I^c(t_1) @>>> \bigoplus \I^c(t_1)_v\\
@VVV @VVV\\
\I^c(t_2) @>>> \bigoplus \I^c(t_2)_v.
\end{CD}
$$
Since $\mathfrak o_{t_1}^c=\mathfrak u_1=\mathfrak u_2=\mathfrak o_{t_2}^c$,
it follows from \eqref{ocs} that 
the isomorphisms $\I^c(t_1)_v\to\I^c(t_2)_v$ extend to isomorphisms 
$\L^c_v(s)\to\L^c_v(s)$ making the extended diagram commute:
$$
\begin{CD}
\I^c(t_1) @>>> \bigoplus_{v\in V} \I^c(t_1)_v @>>> \bigoplus_{v\in V} \L^c_v(s)\\
@VVV @VVV @VVV\\
\I^c(t_2) @>>> \bigoplus_{v\in V} \I^c(t_2)_v @>>> \bigoplus_{v\in V} \L^c_v(s).
\end{CD}
$$
But an automorphism of $\L^c_v(s)$ is
multiplication by a certain $z_v\in\mathbf
G_{\mathbf m}(\k)$. It follows that $\I^c(t_2)=(z_v\,|\,v\in V)\, \I^c(t_1)$
as subsheaves of $\oplus \L^c_v(s)$, and hence that $t_1$ and $t_2$
differ by the action of a $\k$-point on 
$\mathbf G_{\mathbf m}^V/\mathbf G_{\mathbf m}$.
\end{proof}


\section{Degenerations of line bundles I}\label{degbdles}

\subsection{Admissible extensions.}
Let $B$ be the spectrum of the power series 
ring $\k[[t]]$. Let $\pi\colon\X\to B$ be a flat, projective map with
smooth generic fiber and special fiber isomorphic to $X$. We say 
$\pi$ is a \emph{smoothing} of $X$. We will identify the special fiber
with $X$, through a fixed isomorphism. For each $e\in\mathbb E$, 
let $\ell_e$ be the \emph{singularity degree} of $\X$ at $N_e$. More
precisely, $\ell_e$ is the unique positive integer such that the
completion of the local ring of $\X$ at $N_e$ is $\k[[t]]$-isomorphic
to $\k[[x,y,t]]/(xy-t^{\ell_e})$. Let 
$\ell\colon E\to\mathbb N$ be the function assigning $\ell_e$ to $e$ for each
$e\in E$. For each $v\in V$, the irreducible component $X_v$ of $X$
fails to be a Cartier divisor of $\X$ at the nodes $N_e$ for which
$\ell_e>1$. 

By a sequence of blowups supported over the nodes $N_e$, we obtain a map 
$\sigma\:\wt\X\to\X$ for which the strict transform of $X_v$ is
Cartier in $\wt\X$ for each $v\in V$. Then all irreducible components
of $\sigma^{-1}(X)$ are Cartier in $\wt\X$. We choose $\wt\X$ minimal with
this property. (Observe that $\sigma$ is an isomorphism over any node
of $X$ other than the $N_e$, so $\wt\X$ is not regular if the
singularity degree of any such node is greater than 1.) We call
$\sigma$ the \emph{Cartier reduction} of $\pi$. 
It is the semistable reduction if $\wt\X$ fails to be regular
only at the $N_e$. 

Put $\wt\pi:=\pi\sigma$ and let $X^\ell$ be the special
fiber of $\wt\pi$. It is obtained by 
replacing each $N_e$ by a chain $Z^e:=\sigma^{-1}(N_e)$ of rational smooth
curves of length $\ell_e-1$. For each $e=uv$ in $\mathbb E$ 
we will order the components of $Z^e$
as $Z^{e}_{1},\dots,Z^{e}_{\ell_e-1}$, where $Z^{e}_{1}$
intersects $X_u$ and $Z^{e}_{\ell_e-1}$ intersects
$X_v$. Sometimes it will be convenient to set 
$Z^{e}_0:=X_u$ and $Z^{e}_{\ell_e}:=X_v$. 
We will identify the remaining components of 
$X^\ell$ with their corresponding images on $X$. The curve 
$X^\ell$ and the restriction 
$\sigma^\ell:=\sigma|_{X^\ell}\:X^\ell\to X$ depend only on $\ell$.

Let $H$ be the graph obtained as a subdivision of $G$ by inserting
$\ell_e-1$ new vertices in the middle of each edge, that
is, by replacing each $e=uv$ in $\mathbb E$ by a path 
$P_e=uz_1^ez_2^e\dots z^e_{\ell_e-1}v$, for new vertices
$z_1^e,\dots,z^e_{\ell_e-1}$. Sometimes it will be convenient to set 
$z^{e}_0:=u$ and $z^{e}_{\ell_e}:=v$. 
The graph $H$ coincides precisely with the graph obtained from
the dual graph of $X^\ell$ by removing all loop edges: the new
components $Z^{e}_i$ correspond to the vertices $z^e_i$ and the nodes 
$Z^{e}_i\cap Z^{e}_{i+1}$ to the edges $\{z^e_i,z^e_{i+1}\}$ on the path $P_e$. 
Abusing the notation, given a
vertex $v$ (resp.~edge $e$) of $H$, we will also denote by $X_v$
(resp.~$N_e$) the corresponding component (resp.~node) of $X^\ell$. To
avoid confusion, we will sometimes denote by $V(G)$ and $E(G)$ the
vertex and edge sets of $G$, likewise for $H$. We will also view
$V(G)$ as a subset of $V(H)$ in the natural way.

To an \emph{almost invertible} sheaf $L$ on $X$, that is, a
torsion-free, rank-one sheaf on $X$ that is invertible at all the
nodes $N_e$ for $e\in E(G)$, we may associate a divisor
$D\in\text{Div}(G)$, 
namely, $D:=\sum\deg(L|_{X_v})v$. Likewise for sheaves on
$X^\ell$ and divisors on $H$. Recall from~\cite[Subsection 2.3]{AE1} that a divisor $D$ on $H$ is
called \emph{$G$-admissible} if $D$ contains in its support at most one vertex among the new
vertices $z_1^e,\dots,z_{\ell_e-1}^e$ for each edge $e\in\E(G)$, and if so,
with coefficient equal to 1. We set
$$
t^D_e:=\sum_{j=1}^{\ell_e}(\ell_e-j)D(z^e_j)\quad\text{for each
}D\in\text{Div}(H)\text{ and }e\in\E(G)
$$
to identify the vertex. Similarly, we say an almost invertible sheaf $L$ on $X^\ell$ is
\emph{$\sigma^\ell$-admissible} if its associated divisor is
$G$-admissible. If $L$ is a $\sigma^\ell$-admissible almost 
invertible sheaf of degree $d$, then the push-forward $\sigma^\ell_*L$ is a
torsion-free, rank-one sheaf on $X$ of degree $d$; see
e.g.~\cite{EP16}, Thm.~3.1, p.~63 
(where a more general notion of admissibility, for which the statement 
still holds, is considered).

Let now $L_\eta$ be an invertible sheaf on the generic fiber of $\pi$,
which is the same as that of $\wt\pi$. Let $d:=\deg L_\eta$. 
Since $\sigma$ is the Cartier reduction of $\pi$, the sheaf $L_\eta$
extends 
to an \emph{almost invertible} sheaf on 
$\wt\X$, that is, a relative torsion-free, rank-one sheaf on $\wt\X/B$
whose restriction to the special fiber $X^{\ell}$ is almost invertible.
The extension is not unique though: Given an almost invertible sheaf
$\L$ on $\wt\X$, all the sheaves of the form
$$
\L(f):=\L\otimes\O_{\wt\X}\Big(\sum_{v\in V(H)}f(v)X_v\Big)
$$
for $f\in C^0(H,\Z)$ have the same restriction to $X_\eta$ as $\L$, and these
are the only almost invertible sheaves with this property. The
invertible sheaves
$$
\T_f:=\O_{\wt\X}\Big(\sum_{v\in V(H)}f(v)X_v\Big)
$$
are called \emph{twisters}. Given $f\in C^0(H,\Z)$, the divisor on $H$
associated to $\T_f|_{X^\ell}$ 
does not depend on the choice of the smoothing $\pi$. In addition, two
twisters $\T_f$ and $\T_h$ yield the same divisor if and
only if $f-h$ is constant.

The divisor on $H$ associated to $\T_f|_{X^{\ell}}$ is the principal divisor defined by
$$
\text{\rm div}(f):=\sum_{e\in\mathbb E(H)}(f(\te_e)-f(\he_e))\he_e
$$
for each $f\in C^0(H,\Z)$. The assignment defines a group homomorphism 
$$
\text{div}\: C^0(H,\Z)\to\text{Div}(H)
$$
whose kernel is the subgroup
of constant functions and whose image defines an equivalence relation
on $\text{Div}(H)$, called \emph{linear
  equivalence}~\cite{BHN97, BN06}. It follows that the divisors in
$\text{Div}(H)$ associated to the $\L|_{X^\ell}$ for all almost invertible
extensions $\mathcal L$ of $L_\eta$ to $\wt\X$ are those in a 
certain linear equivalence class.

\medskip

We are interested in special extensions of $L_\eta$ to $\wt\X$, those
that yield meaningful extensions of $L_\eta$ to $\X$ as well, 
so we make the following definition.

\begin{defi} \rm An almost invertible sheaf $\L$ on $\wt\X$ 
is called $\sigma$-\emph{admissible}
if $\L|_{X^\ell}$ is $\sigma^\ell$-admissible.  
\end{defi}

The admissible sheaves we consider are special in the sense that the
last statement of the following proposition holds.

\begin{prop}\label{ILv} Let $\L$ be a $\sigma$-admissible almost invertible
  sheaf on $\wt\X$. Then 
$R^1\sigma_*\L=0$ and $\I:=\sigma_*\L$ is a (relative) 
torsion-free, rank-one sheaf 
on $\X/B$, with
formation commuting with base change. In particular, the restriction of
$\I$ to $X_\eta$ is the same as that of $\L$. Furthermore, $I_v\cong\L|_{X_v}$ for
  each $v\in V(G)$, where  $I_v$ is the restriction of $\I$ to
$X_v$ modulo torsion.
\end{prop}

\begin{proof} The first statement follows from \cite{EP16}, Thm.~3.1,
  p.~63. As for the last statement, set $L:=\L|_{X^\ell}$ and
  $I:=\I|_X$. By the first statement, $I=\sigma^\ell_*L$. Let
$Z\subseteq X^{\ell}$ be the union of the chains $Z^e$ of exceptional components
over which $L$ has total degree 1, and 
$W\subseteq X^{\ell}$ the union of the remaining components. Let
$\mu\:X'\to X$ be the partial normalization of $X$ along the nodes $N_e$
at which $I$ fails to be invertible.  Since $I=\sigma^\ell_*L$, it
follows from Thm.~3.1 in loc.~cit.~that $\sigma^{\ell}$
induces upon restriction a map $\nu\:W\to X'$. It follows
from Lemma~2.1 in loc.~cit.~that applying $\sigma^{\ell}_*$ 
to the short exact sequence
$$
0\longrightarrow L|_Z\Big(-\sum_{P\in Z\cap W}P\Big)
\longrightarrow L \longrightarrow L|_{W} \longrightarrow 0
$$
yields an isomorphism $I\cong\mu_*\nu_*(L|_W)$. 
As $W$ is a semistable model of $X'$, and $L|_W$ has degree 0
on every exceptional component of $W$, it follows
again from Thm.~3.1 in loc.~cit.~that $\nu_*(L|_W)$ is almost 
invertible and $\nu^*\nu_*(L|_W)=L|_W$. In particular, 
$\nu_*(L|_W)|_{X_v}\cong L|_{X_v}$ for each
$v\in V$. Finally, since $\mu$ is finite, the natural map
$\mu^*\mu_*\nu_*(L|_W)\to \nu_*(L|_W)$ is surjective. Thus,
restricting it to $X_v$ 
we get a
surjection $\mu_*\nu_*(L|_W)|_{X_v}\to L|_{X_v}$. Since $L$ is almost invertible,
whence $L|_{X_v}$ is torsion-free, and $I\cong\mu_*\nu_*(L|_W)$, it follows
that $I_v\cong L|_{X_v}$, finishing the proof.
\end{proof}

Admissible extensions exist, according to the following proposition.

\begin{prop}\label{ext} There is a $\sigma$-admissible almost
  invertible extension $\L$ of $L_\eta$ 
to $\wt\X$. Furthermore, if $L_\eta$ admits a nonzero section, there is such an
extension whose restriction to $X^\ell$ has an effective associated
divisor in $\text{\rm Div}(H)$.
\end{prop}

\begin{proof} The first statement is essentially a restatement of 
the results proved in~\cite[Section 2]{AE1}. Indeed, it is enough to see
that any divisor on $H$ is linearly equivalent to a $G$-admissible
divisor. This is proved in~\cite[Thm.~2.10]{AE1}. 
As for the second statement,
first observe that for any almost invertible extension $\L$ of $L_\eta$, any nonzero section of $L_\eta$ extends to a 
section of $\L$ with a nonzero restriction 
to $\L|_{X^\ell}$. For each component $X_v$ of $X^\ell$, this section 
vanishes to a certain finite 
order, say $-f(v)$, on $X_v$. So it 
induces a section of $\M:=\L\otimes\T_f$ whose zero scheme is 
finite over $B$. In particular, the divisor $D\in\text{Div}(H)$ associated to 
$\M|_{X^\ell}$ is effective. Hence $t_e^D\geq 0$ and thus 
$\delta_e(0;t^D):=\lfloor\frac{t_e^D}{\ell_e}\rfloor\geq 0$ 
for each $e\in\E(G)$. The statement follows now 
from \cite[Prop.~2.9]{AE1}, which yields that 
the $G$-admissible divisor $D'= D+ \div(h)$ is effective, where $h\in
C^0(H,\Z)$ is the canonical extension of the zero function in
$C^0(G,\Z)$ with respect to $D$.
\end{proof}

Special twisters allow us to pass from one $\sigma$-admissible almost
invertible sheaf on $\wt\X$ to any other with the same restriction to $X_\eta$.

First, we introduce notation. By \cite{AE1}, Prop.~2.7, for each divisor
$D\in\text{Div}(H)$ and each $f\in C^0(G,\Z)$, there is a unique  extension
$\wt f\in C^0(H,\Z)$ of $f$ such that $D+\text{div}(\wt f)$ is
$G$-admissible. The function $\wt f$ is called the \emph{canonical
  extension} of $f$ with respect to $D$. Thus, for each almost
invertible sheaf $\L$ on $\wt\X$ with $D$ as the associated divisor to
$\L|_{X^{\ell}}$, the sheaf $\L(\wt f)$ is
$\sigma$-admissible. As $\wt f$ depends in this case on $\L$ and $f$,
we abuse the notation by setting $\L(f):=\L(\wt f)$.

We will often deal with the special case where $D\in\text{Div}(H)$ is
$G$-admissible and $f=\chi_{\indm v}\in C^0(G,\Z)$ for a vertex $v\in V$. In
this case we denote $\wt f$ by $f_{D,v}$. We have that 
$f_{D,v}(z^e_i)=1$ for each $e\in\mathbb E$ with 
$\he_e=v$ and each $i=\ell_e-t^D_e,\dots,\ell_e$, whereas $f_{D,v}(w)=0$ for all other
vertices $w\in V(H)$. Also, for a $\sigma$-admissible almost
invertible sheaf $\L$ on $\wt\X$ with $D$ as the associated divisor to
$\L|_{X^{\ell}}$, we denote $\L(f)$ by $\mathcal M_v(\L)$.
The divisor on $H$ associated to $\mathcal M_v(\L)|_{X^\ell}$ is the $G$-admissible chip
firing move of $D$ at $v$, denoted $M_v(D)$ in~\cite[Subsection~2.6]{AE1}. 
So the notations are coherent.

\begin{prop}\label{admtw} 
Let $\L$ be a $\sigma$-admissible almost invertible sheaf on $\wt\X$. Then: 
\begin{enumerate}
\item For each pair of vertices $v,w\in V(G)$, we have $\M_v(\M_w(\L))=\M_w(\M_v(\L))$.
\item Enumerating the vertices of $G$ as $v_1, \dots, v_n$, we have 
$$
\M_{v_1}(\M_{v_2}(\cdots\M_{v_n}(\L)\cdots))\cong\L.
$$
\item For each $\sigma$-admissible almost invertible $\M$ on $\wt\X$ with the same
  restriction to $X_\eta$ as $\L$, there 
exists a sequence $v_1,\dots,v_m$ of vertices of $G$ such that 
$$
\M\cong \M_{v_1}(\M_{v_2}(\cdots\M_{v_m}(\L)\cdots)).
$$
Furthermore, 
the sequence is unique up to reordering the vertices and adding or
subtracting all
the vertices of $G$.
\end{enumerate}
\end{prop}

\begin{proof} Statements (1) and (2) are easy to check, and follow from
  \cite{AE1}, Prop.~2.11. The existence part
  in Statement (3) follows from Prop.~2.14 in loc.~cit.

  As for the 
  uniqueness part in Statement (3), the first two statements yield that the mentioned
  operations to a sequence $v_1,\dots,v_m$ do not change the resulting
  sheaf. In addition, the proof of Prop.~2.14 in loc.~cit.~shows that
  $D'=D+\div(h)$, where $D$ (resp.~$D'$) is the divisor in $H$ associated to
  $\L|_{X^\ell}$ (resp.~$\M|_{X^\ell}$),   and $h\in C^0(H,\Z)$ is the
  canonical extension of $f:=\chi_{\indmbi{v}{1}}+\cdots+\chi_{\indmbi{v}{m}}$ with
  respect to $D$. Thus, for any other sequence
  $v'_1,\dots,v'_p$ of vertices of $G$ with the same property as
  $v_1,\dots,v_m$, we have $\div(h)=\div(h')$, where $h'$ is the
  canonical extension of 
$f':=\chi_{v'_1}+\cdots+\chi_{v'_p}$ with respect to $D$. Then 
$h'-h$ is constant. Since $h'-h$ extends $f'-f$, it
follows that
$f'-f$ is constant, finishing the proof.
\end{proof}

Let $\L$ be a $\sigma$-admissible almost invertible sheaf on $\wt\X$. 
Set $\I(f):=\sigma_*\L(f)$ for each
$f\in C^0(G,\Z)$.  It will be useful to have a direct interpretation of
the $\I(f)$. It is given by Proposition~\ref{kerI} below, from which
we see how to obtain the $\I(f)$ recursively from 
$\I:=\sigma_*\L$. 

\begin{prop}\label{kerI} Notations as above, for each $v\in V(G)$ and $f\in
  C^0(G,\Z)$, the sheaf $\I(f)$ is the
kernel of the surjection
\begin{equation}\label{Iv}
\mathcal I(f+\chi_{\indm v}) \longrightarrow \mathcal I(f+\chi_{\indm v})_v,
\end{equation}
where the sheaf to the right is the torsion-free sheaf generated on
$X_v$ by $\mathcal I(f+\chi_{\indm v})$.
\end{prop} 

\begin{proof} The kernel of
the surjection~\eqref{Iv} is a relatively torsion-free, rank-1 sheaf on $\mathcal
X/B$ by an argument analogous to the one found in \cite{Langton},
Prop.~6, p.~100; see \cite{Esteves01}, Section~3. 
On the other hand, there is a natural exact sequence
on $\wt{\mathcal X}$,
\begin{equation}\label{IL}
\begin{CD}
0 @>>> \L(f) @>>> \L(f+\chi_{\indm v}) @>>> \L(f+\chi_{\indm v})|_Y @>>> 0,
\end{CD}
\end{equation}
where $Y$ is the subcurve of $X^{\ell}$ which is the union of $X_v$ and the
components $Z^e_j$ for all $e\in\E$
with $\he_e=v$ and $j=\ell_e-t^D_e+1,\dots,\ell_e-1$, 
where $D\in\text{Div}(H)$ is the divisor associated to 
$\L(f+\chi_{\indm v})|_{X^\ell}$. Since $\L(f+\chi_{\indm v})$ has degree 0 on these
$Z^e_j$, we have 
$$
\sigma_*(\L(f+\chi_{\indm v})|_Y)=\L(f+\chi_{\indm v})|_{X_v}=\mathcal
I(f+\chi_{\indm v})_v,
$$
and the surjection in \eqref{IL} restricts to the
surjection \eqref{Iv}. It follows that 
$\sigma_*\L(f)$ is the kernel of  \eqref{Iv}, as claimed.
\end{proof}

\subsection{Generalized enriched structures}\label{ges} We apply the
above recursive interpretation to describe the $\mathcal O_{\X}(f)$. (Here 
$\mathcal O_{\X}(f):=\sigma_*\mathcal O_{\wt\X}(f)$ for each $f\in
C^0(G,\Z)$.) That interpretation allows us to view
$\mathcal O_{\mathcal X}(f)$ as a sheaf of
fractional ideals of $\mathcal X$, which we will describe below.

\smallskip

In her thesis~\cite{Maino}, Main\`o defined the notion of enriched structures
over $X$, and constructed the moduli space of enriched
curves, that is, stable curves with enriched structures. 
Under the interpretation given by the second author and Medeiros~\cite{EM}, an enriched
structure on $X$ arises from a smoothing $\mathcal X\to B$ with
regular total space $\mathcal X$, thus $\ell_e=1$ for each $e\in
E(G)$, as the group homomorphism $L\:C^0(G,\Z)\to\text{Pic}^0(X)$ given
by $L(f):=\mathcal O_{\mathcal X}(f)|_X$ for each $f\in C^0(G,\Z)$. 

In the general case of higher singularity degrees, we may thus see the
$\mathcal O_{\X}(f)|_X$ as part of a generalized enriched structure over
$X$; the precise definition is given in Definition~\ref{stlim}. Main\`o constructed in
loc.~cit.~a quasiprojective variety parameterizing enriched curves 
over the moduli space of stable curves. One of our goals, to be pursued in
a subsequent work with the contributions given here, is to compactify
this variety in a meaningful way. 
(We refer to recent work by Biesel and Holmes~\cite{BH16} 
for a compactification of Main\`o's moduli space following a different
approach.)

Let $V/M$ be the versal deformation of $X$. As explained in
\cite{DM69}, pp.~79--81 and reviewed in \cite{EM}, 
pp.~288--9, we have $M=\text{Spec}(R)$, where $R$ is the power series ring over $\k$ 
in the variables $t_e$, for $e\in E^{\mathfrak o}$, and $s_1,\dots,s_p$, 
for a certain integer $p$. Furthermore, the variables can be chosen so
that for each $e\in E^{\mathfrak o}$ 
we have an isomorphism of $R$-algebras,
\[
\psi_e\:\widehat{\mathcal O}_{V,N_e}\to R[[z_e,w_e]]/(z_e w_e-t_e).
\]
The versal deformation comes with an identification of the special
fiber of $V/M$ with $X$, thus we may view $X\subseteq V$. We may assume that
$z_e=0$ corresponds to the component $X_v$ and $w_e=0$ to
$X_u$, where $v:=\he_e$ and $u:=\te_e$. 
Letting $\hat z_e$ and $\hat w_e$ denote the elements of
$\widehat{\mathcal O}_{V,N_e}$ corresponding to $z_e$ and $w_e$, we have that
$\hat z_e$ restricts to a local parameter $z_e$ of $\widehat{\mathcal O}_{X_u,N_e}$, whereas $\hat w_e$ restricts to a local parameter $w_e$ of 
$\widehat{\mathcal O}_{X_v,N_e}$. 

As $\mathcal X/B$ is a deformation of $X$, there is a natural
Cartesian diagram factoring the inclusion of $X$ in $V$:
$$
\begin{CD}
X @>>> \mathcal X @>>> V\\
@VVV @V\pi VV @VVV\\
\text{Spec}(\k) @>>> B @>>> M
\end{CD}
$$
The map $B\to M$ sends via pullback $t_e$ to $a_et^{\ell_e}\xi_e$,
for certain $a_e\in \k$ and $\xi_e\in \k[[t]]$ with $\xi_e(0)=1$. If we
denote by $\tilde z_e$ and $\tilde w_e$ the pullbacks under the map $\mathcal
X\to V$ of $\hat z_e$ and $\hat w_e$, respectively, we have $\tilde
z_e\tilde w_e=a_et^{\ell_e}\xi_e$ in $\widehat
\O_{\mathcal X,N_e}$ for each $e\in E^{\mathfrak o}$.

\smallskip

Now, for each $f\in C^0(G,\Z)$ define
$$
\delta_{e}(f):=\Big\lfloor\frac{f(v)-f(u)}{\ell_{e}}
\Big\rfloor\quad\text{for each }e=uv\in\E.
$$
We claim that, for each $e=uv\in E^{\mathfrak o}$, the sheaf of
fractional ideals $\mathcal O_{\mathcal X}(f)$ is 
generated locally analytically at $N_e$ by 
\begin{equation}\label{fracid}
(t^{-f(u)}\tilde z_e^{-\delta_e(f)},t^{-f(v)}\tilde w_e^{-\delta_{\ol
    e}(f)}).
\end{equation}
Indeed, suppose $f(v)\geq f(u)$. As $t=0$ gives $X$, we
have that 
$$
\mathcal O_{\mathcal X}(f)=t^{-f(v)}\mathcal O_{\mathcal
  X}\big(f-f(v)\chi_V\big), \text{ where }\chi_V:=\sum\chi_{\indm v}.
$$
And, as shown in \cite{CEG}, p.~14, locally
analytically at $N_e$ the sheaf $\mathcal O_{\mathcal
  X}(f-f(v)\chi_V)$ is the ideal generated by $(\tilde
w_e^{q+1},\tilde w_e^qt^r)$,
where $q$ and $r$ are the quotient and the remainder of the Euclidean
division of $f(v)-f(u)$ by $\ell_e$. In other words, $q=\delta_{e}(f)$
and $r=f(v)-f(u)-\ell_eq$. It is now easy to check, using $\tilde
z_e\tilde w_e=a_et^{\ell_e}\xi_e$, that
$t^{-f(v)}(\tilde w_e^{q+1},\tilde w_e^qt^r)$ is the fractional ideal
\eqref{fracid}. An analogous argument works when $f(v)\leq
f(u)$. Notice that 
$$
t^{-f(v)}\tilde w_e^{-\delta_{\ol e}(f)}
=(a_e\xi_e)^{\delta_e(f)}t^{-f(u)}\tilde z_e^{-\delta_e(f)}
$$ 
if $\ell_e$ divides $f(v)-f(u)$; in particular, $\mathcal
O_{\mathcal X}(f)$ is principal at $N_e$ in this case.

For each $f\in C^0(G,\Z)$, the sheaf of fractional ideals
$t^{f(v)}\mathcal O_{\mathcal X}(f)$ of $\mathcal X$ generates a
sheaf of fractional ideals $\mathcal J_v(f)$ of $X_v$ for each $v\in V(G)$. Given
the above description, for each $e=uv\in E^{\mathfrak o}$ we have
that, locally analytically at $N_e$, the element
$\tilde z_e^{-\delta_e(f)}$ is mapped to $z_e^{-\delta_e(f)}$ in
$\mathcal J_u(f)$, whereas $\tilde w_e^{-\delta_{\ol e}(f)}$ is mapped
to $w_e^{-\delta_{\ol e}(f)}$ in 
$\mathcal J_v(f)$. On the other hand, if $\ell_e$ does not divide
$f(v)-f(u)$, then $t^{f(u)-f(v)}\tilde
w_e^{-\delta_{\ol e}(f)}$ and $t^{f(v)-f(u)}\tilde z_e^{-\delta_e(f)}$
are mapped to 0 in $\mathcal J_u(f)$ and $\mathcal J_v(f)$,
respectively. It follows that
$$
\mathcal J_v(f)=\mathcal O_{X_v}\Big(\sum_{\substack{e\in\E \\ \he_e=v}}\delta_{\ol e}(f)N_e\Big)
$$
for each $v\in V(G)$.

\smallskip

As in \cite{EM}, p.~292, we may
consider the image $\mathcal J(f)$ of the composition
$$
\begin{CD}
\mathcal O_{\mathcal X}(f) @>>> \bigoplus_{v\in V}t^{f(v)}\mathcal
O_{\mathcal X}(f) @>>> \bigoplus_{v\in V}\mathcal J_v(f).
\end{CD}
$$
Then $\mathcal J(f)$ is a sheaf of fractional ideals of $X$ isomorphic to
$\mathcal O_{\mathcal X}(f)|_X$. Clearly, $\mathcal J(f)$ generates
$\mathcal J_v(f)$ for each $v\in V(G)$. As the $\mathcal J_v(f)$ are
described above, to describe the subsheaf $\mathcal J(f)$ we need only
describe it locally analytically at the $N_e$ where $\mathcal J(f)$ is
invertible, that is, for $e=uv\in E^{\mathfrak
  o}$ such that $\ell_e$ divides $f(v)-f(u)$.

Let $K:=\prod_{v\in V}\k(X_v)$, the product of the fields of rational functions of the irreducible components of $X$. And put
$\mathcal K:=\prod_{v\in V}\mathcal K_v$, the product of the constant
sheaves of rational functions of the irreducible components of
$X$. For each $e\in E$ and $v\in e$, let $\widehat K_{v,e}$ denote the
field of fractions of $\widehat{\mathcal
  O}_{X_v,N_e}$; it contains the field of fractions of $\mathcal
  O_{X_v,N_e}$, which is $\k(X_v)$. We may thus view any local description of
  a sheaf of fractional ideals of $X$ at $N_e$ in $\widehat
  K_{u,e}\times\widehat K_{v,e}$, for each $e=uv\in E^{\mathfrak
    o}$. 

In particular, for each $f\in C^0(G,\Z)$ and each $e=uv\in
E^{\mathfrak o}$, the sheaf of fractional
ideals $\mathcal J(f)$ is generated at $N_e$ locally analytically in
$\widehat K_{u,e}\times\widehat K_{v,e}$ by
$(a_e^{\delta_e(f)}z_e^{-\delta_e(f)},w_e^{-\delta_{\ol e}(f)})$ if $\ell_e$ divides
$f(v)-f(u)$ and by $((z_e^{-\delta_e(f)},0),(0,w_e^{-\delta_{\ol
    e}(f)}))$ otherwise.

For each $e=uv\in E^{\mathfrak o}$, we may use $z_e$ to establish 
isomorphisms $\alpha_{e,m}\:\mathcal
O_{X_u}(mN_e)|_{N_e}\cong \k$ for each $m\in\Z$, by analytically identifying
$\mathcal O_{X_u}(mN_e)$ at $N_e$ with the fractional ideal generated
by $z_e^{-m}$, and taking $z_e^{-m}$ to 1. Doing the same for $w_e$, 
we get isomorphisms $\beta_{e,m}\:\mathcal
O_{X_v}(mN_e)|_{N_e}\cong \k$. Under these isomorphisms the sheaf
$\mathcal J(f)$ is the subsheaf of $\bigoplus_{v\in V}\mathcal J(f)_v$ given
locally at $N_e$, for each 
$e=uv\in E^{\mathfrak o}$ such
that $\ell_e$ divides $f(v)-f(u)$, as the kernel of the surjection
$$
\begin{CD}
\mathcal O_{X_u}(\delta_{e}(f)N_e)|_{N_e}\oplus \mathcal
O_{X_v}(\delta_{\ol e}(f)N_e)|_{N_e}
@>(\alpha_{e,\delta_e(f)},\beta_{e,\delta_{\ol e}(f)})>\cong >
\k\oplus \k @>(1,-a_e^{\delta_e(f)})>> \k.
\end{CD}
$$

Notice that as $\mathcal X/B$ varies among smoothings of $X$ with
singularity degree function $\ell$, the $a_e$ on which the $\mathcal
O_{\mathcal X}(f)$ ultimately depends vary freely. The $\mathcal
O_{\mathcal X}(f)$ depends thus on the free choice
of a homomorphism $a\: C^1(G,\Z)\to {\mathbf G}_{\mathbf m}$. 
Also, by Nakayama Lemma, giving an isomorphism $\alpha_e\:\mathcal
O_{X_v}(N_e)|_{N_e}\cong \k$ for each $e\in E$ and $v\in e$ is the same
as choosing an analytic local parameter
for $X_v$ at $N_e$. We have chosen above, as a result of considering
a versal deformation of $X$, certain analytic local parameters
$z_e$ and $w_e$ for $X_u$ and $X_v$ at $N_e$  for each $e=uv\in
E^{\mathfrak o}$. Different choices $z'_e$ and $w'_e$ can be expressed as
power series $z'_e=\tau_ez_e+\cdots$ and $w'_e=\sigma_ew_e+\cdots$ for
$\tau_e,\sigma_e\in \k$, and we would obtain the same subsheaf
$\mathcal J(f)$ of 
$\bigoplus_{v\in V}\mathcal J(f)_v$ by replacing the $a_e$ by
$\tau_ea_e\sigma_e$ for each $e\in E^{\mathfrak o}$. 

\section{Degenerations of line bundles II}\label{stackydeg}

Let $B$ be the spectrum of the power series ring $\k[[t]]$. 
Let $\pi\colon\X\to B$ be a smoothing of $X$. 
Let $\ell\colon E\to\mathbb N$ be the function assigning to $e$ the
singularity degree $\ell_e$ of $\X$ at $N_e$. Let $\sigma\:\wt\X\to\X$ be the Cartier reduction of 
$\X$. Put $\wt\pi:=\pi\sigma$. Let $X^\ell$ be the special fiber of
$\wt\pi$. Keep the remaining notation of Section~\ref{degbdles}. 
As $f$ runs in
$C^0(G,\Z)$, the sheaf $\mathcal I(f)|_X$ runs through what we call \emph{limits} of
$L_\eta$. These are not all of what we call \emph{stable limits}
though. To obtain all of them, we consider base changes $t\mapsto t^n$
for positive integers $n$, as explained below.

For each $n\in\N$, let $\mu_n\:B\to B$ be the
base change map given by $t\mapsto t^n$, and let $\pi^n\:\mathcal
X^n\to B$ be the base
extension of $\pi$. Let $\sigma^n\:\wt{\mathcal X}^n\to\mathcal X^n$ be
the Cartier reduction. If $m\in\N$ divides $n$, then we have the
following commutative diagram of maps:
$$
\begin{CD}
\wt{\mathcal X}^n @>\sigma^n>> \mathcal X^n @>\pi^n>> B\\
@VVV @VVV @V\mu_{n/m}VV\\
\wt{\mathcal X}^m @>\sigma^m >> \mathcal X^m @>\pi^m >> B.
\end{CD}
$$
The square to the right is Cartesian, so the geometric fibers of $\mathcal
X^n/B$ are the same as those of $\mathcal X/B$ but the one to the
left is not. The singularity degree of $\mathcal X^n$ at $N_e$ is now
$n\ell_e$ for each $e\in E(G)$, so the special fiber of $\wt{\mathcal
  X}^n/B$ is $X^{n\ell}$. We have the same configuration as before. 
The only difference is that $\ell$ is
replaced by $n\ell$.

Given an invertible sheaf $L_\eta$ on the generic fiber of $\pi^m$,
for given $m\in\N$, and given $n\in\N$ divisible by $m$, we
may pull $L_\eta$ back to an invertible sheaf $L_\eta^n$ on the generic fiber of
$\pi^n$, and consider its extensions to 
$\wt{\mathcal X}^n$, as we did before for the case $n=1$. Given an
extension $\mathcal L$ of $L_\eta$ to $\wt{\mathcal X}^m$ we may pull it
back to an extension $\mathcal L^n$ of $L_\eta^n$ to $\wt{\mathcal
  X}^n$. If $\mathcal L$ is
$\sigma^m$-admissible, 
then  $\mathcal L^n$
is $\sigma^n$-admissible. We will also denote by $\mathcal I^n$ the
pullback of a relative torsion-free rank-one sheaf $\mathcal I$ on $\mathcal X^m/B$; it
is one on $\mathcal X^n/B$. If $\I=\sigma^m_*\L$ then
$\I^n=\sigma^n_*\L^n$. 

As before, to each $f\in C^0(G,\Z)$ and each $\sigma^m$-admissible almost invertible 
extension $\mathcal L$ of $L_\eta$ to $\wt{\mathcal X}^m$, we associate a
$\sigma^n$-admissible extension $\mathcal L^n(f)$ of $L_{\eta}^n$. 
Again by \cite{EP16}, Thm.~3.1, p.~63, the pushforward 
$$
\mathcal I^n(f):=\sigma^n_*\mathcal L^n(f)
$$
is a relative torsion-free rank-one sheaf on $\mathcal X^n/B$. The
notation is consistent, as $\I^n(0)=\I^n$.

\begin{defi}\label{stlim} \rm Let $m\in\N$ and $L_\eta$ be an invertible sheaf on
  the generic fiber of $\mathcal X^m/B$. Let $\L$ be a $\sigma^m$-admissible
  almost invertible extension
  of $L_\eta$ to $\wt{\mathcal X}^m$. For each $n\in\N$
  divisible by $m$ and $f\in
  C^0(G,\Z)$, let 
$$
L^n(f):=\mathcal L^n(f)|_{X^{n\ell}}\quad\text{and}\quad
I^n(f):=\mathcal I^n(f)|_X.
$$
We call 
$$
\mathfrak I:=\{I^{n}(f)\,|\, n\in m\mathbb N,\, f\in C^0(G,\Z)\}
$$
the collection of
\emph{stable limits} of $L_\eta$. In case $L_\eta$ is the structure sheaf of the generic fiber of
$\X/B$, we call $\mathfrak I$ a \emph{generalized enriched structure}.
\end{defi}

If $\ell$ is the constant function 1, and $L_\eta$ is the structure sheaf of the generic fiber of
$\X/B$, Main\`o called an enriched
structure the subset $\{I^1(\chi_{\indm v})\,|\,v\in V(G)\}$. In this
case, the full set $\mathfrak I$ is obtained from this subset by
tensor products and degeneration.

\smallskip

In general, $\mathfrak I$ is a set of torsion-free, rank-one sheaves
of degree equal to that of $L_\eta$. If $\k$ has characteristic 0, then
the field of Puiseux series is the algebraic closure of the field of
Laurent series, the field of fractions of $\k[[t]]$, and we may thus
think of $\mathfrak I$ as the
collection of all the limits of the pullback of $L_\eta$ to the
geometric generic fiber of $\mathcal X/B$.

For each $d\in\Z$, denote by $\mathbf J^d$ the stack parameterizing 
torsion-free, rank-one sheaves of degree $d$ on
$X$; see Section~\ref{embsheaves}. If $L_\eta$ has degree $d$, 
we may view $\mathfrak I$ as a subset of $\mathbf J^d$. The main 
theorem of this section, Theorem~\ref{thm:main4}, 
asserts that the collection
$\mathfrak I$ is the support of a closed substack of $\mathbf J^d$,
and describes thoroughly the structure of this substack.

\smallskip

We need a few preliminary results though. Recall the notation
introduced in \cite{AE1} and \cite{AE2}: For each $n\in\N$, each $\mathfrak
m\in C^{1}(G,\Q)$ and $f\in C^0(G,\Z)$, we set
$$
\delta^{\mathfrak m,n}_{\ell,e}(f):=\Big\lfloor\frac{f(\he_e)-f(\te_e)+n\mathfrak
  m_e}{n\ell_e}\Big\rfloor\quad\text{for each $e\in\E$,}
$$
and let $\mathfrak d_{\ell,f}^{\mathfrak m,n}\in C^1(G,\frac{1}{2}\Z)$ be
defined by putting
$$
\mathfrak d_{\ell,f}^{\mathfrak m,n}(e)=\frac{1}{2}
\Big(\delta^{\mathfrak m,n}_{\ell,e}(f)-\delta^{\mathfrak m,n}_{\ell,\ol e}(f)\Big)
$$
for each $e\in\mathbb E$. If $\ell$ is fixed we drop the subscript $\ell$.

In addition, let $H^n$ denote the graph obtained
from the dual graph of $X^{n\ell}$ by removing the self
loops. For each $e\in E(G)$, we let $Z^{e,n}\subseteq X^{n\ell}$ denote the chain of
rational smooth curves of length $n\ell_e-1$ lying over the node $N_e$
of $X$. For each $e=uv\in\mathbb E$, order the components 
$Z^{e,n}_1,\dots,Z^{e,n}_{n\ell_e-1}$ of $Z^{e,n}$ by assuming that
$Z^{e,n}_i$ intersects $Z^{e,n}_{i+1}$ for $i=1,\dots,n\ell_{e}-2$ and
$Z^{e,n}_1$ intersects $X_u$ at $N_e$ (whence $Z^{e,n}_{n\ell_e-1}$
intersects $X_v$ at $N_e$). For convenience, we set $Z^{e,n}_0:=X_u$
and $Z^{e,n}_{n\ell_e}:=X_v$. We denote by $z_j^{e,n}$ the vertex of
$H^n$ associated to $Z_j^{e,n}$ for each $j=0,\dots,n\ell_e$. As with
$H$ we may view $V(G)$ as a subset of $V(H^n)$, and more generally
$V(H^m)$ as a subset of $V(H^n)$ for each $m\in\N$ with $m|n$.

If $n=1$, the overscript $n$ is dropped throughout.

\begin{prop}\label{IfI0}
Let $n\in\N$ and $D$ be a $G$-admissible divisor of $H^n$. For each 
$e\in E^{\mathfrak o}$, put 
$$
\mathfrak m_e:=\frac{1}{n}\sum_{i=1}^{n\ell_e-1}iD(z^{\ol e,n}_i),
$$
and let $\m$ be the element of $C^1(G,\Q)$ satisfying this. 
Then, for each $f\in C^0(G,\Z)$ and $v\in V$, we have
$$
\O_{\wt{\mathcal X}^n}\big(\,\wt f\,\big)|_{X_v}
\cong 
\O_{X_v}\Big(\sum_{\substack{e\in\E\\ \te_e=v}}
\lfloor \mathfrak d_f^{\mathfrak m,n}(e)\rfloor N_e
+\sum_{\substack{e\in\E\\ \te_{e}=v\\ \mathfrak m_e<0}} N_e\Big),
$$
where $\wt f\in C^0(H^n,\Z)$ is the canonical extension of $f$ with
respect to $D$.
\end{prop}

(Recall that $\wt f$ is the unique extension of $f$ to $V(H^n)$ such
that $D+\div(\wt f)$ is $G$-admissible; see~\cite[Subsection~2.5]{AE1}.)

\begin{proof} Since $D$ is $G$-admissible, $0$ is the
canonical extension of 0 with respect to $D$. Clearly, 
$\O_{\wt{\mathcal X}^n}(0)$ is trivial. On the other hand, if $e\in\E$ is
  such that $\te_e=v$, then $\mathfrak d_0^{\mathfrak m,n}(e)=0$
  unless $\mathfrak m_e\neq 0$. In this case, since $|\m_e|<\ell_e$, 
if $\mathfrak m_e>0$ then  $\mathfrak d_0^{\mathfrak m,n}(e)=1/2$. 
And if $\mathfrak m_e<0$ then 
$\mathfrak d_0^{\mathfrak m,n}(e)=-1/2$. In any case, the
proposition holds for $f=0$.

Observe that $f$ can be replaced by $f+b\chi_{\ind V}$ for any
$b\in\Z$, where $\chi_V:=\sum\chi_v$. 
We may thus assume $f\geq 0$. We argue by induction on
$\deg(f)$, the initial case, $\deg(f)=0$, having just been
considered. 

Suppose first that $f=\chi_w$ for a certain $w\in V(G)$. Then 
$$
\O_{\wt{\mathcal X}^n}\big(\,\wt f\,\big)=\O_{\wt{\mathcal
    X}^n}\Big(\,X_w+\sum_{\substack{e\in\E\\
    \te_e=w}}\sum_{i=1}^{r_e}Z^{e,n}_i\,\Big),
$$
where $r_e=0$ if $\m_e=0$; otherwise $r_e$ is the only integer $i$
such that $D(z^{e,n}_i)=1$. 

If $v\neq w$ then 
$$
\O_{\wt{\mathcal X}^n}\big(\,\wt f\,\big)|_{X_v}\cong\O_{X_v}\Big(
\sum N_e\Big),
$$
where the sum is over the $e\in\E$ such that $e=wv$ and either
$n\ell_e=1$ or $D(z_1^{\bar e,n})=1$. For such $e$ we have that 
$\mathfrak d_f^{\mathfrak m,n}(\bar e)=1$ and $\m_{\bar e}=0$ if
$n\ell_e=1$; and if $n\ell_e>1$ we have that  $\m_{\bar e}=-1/n$ and 
$\mathfrak d_f^{\mathfrak m,n}(\bar e)=0$ if
$e\in E^{\mathfrak o}$, whereas $\m_{\bar e}=(n\ell_e-1)/n$ and 
$\mathfrak d_f^{\mathfrak m,n}(\bar e)=1$ if
$\bar e\in E^{\mathfrak o}$. In all cases, $N_e$ is in the support of 
$$
U:=\sum_{\substack{e\in\E\\ \te_e=v}}
\lfloor \mathfrak d_f^{\mathfrak m,n}(e)\rfloor N_e
+\sum_{\substack{e\in\E\\ \te_{e}=v\\ \mathfrak m_e<0}} N_e
$$
with multiplicity 1. Furthermore, suppose $N_e$ is in the support of $U$
for a certain $e\in\E$. Then $v\in e$ and we may suppose $\he_e=v$. If
$\te_e\neq w$, then, as before, $\lfloor\mathfrak d_f^{\mathfrak m,n}(\bar e)\rfloor=0$ if
$\m_{\bar e}\geq 0$ and $\lfloor\mathfrak d_f^{\mathfrak m,n}(\bar
e)\rfloor=-1$ otherwise. Since $N_e$ is in the support of $U$ we must
then have $e=wv$. Suppose $n\ell_e>1$. If $\m_{\bar e}\geq 0$, since
$|\m_{\bar e}|<\ell_e$, we have  $\lfloor\mathfrak d_f^{\mathfrak
  m,n}(\bar e)\rfloor=0$ unless $n\m_{\bar e}=n\ell_e-1$, in which
case $\lfloor\mathfrak d_f^{\mathfrak
  m,n}(\bar e)\rfloor=1$. But this can only happen if $\bar e\in E^{\mathfrak
  o}$ and $D(z_1^{\ol e,n})=1$. On the other hand, if $\m_{\bar e}< 0$
then $\lfloor\mathfrak d_f^{\mathfrak
  m,n}(\bar e)\rfloor=-1$ unless $n\m_{\bar e}=-1$, in which
case $\lfloor\mathfrak d_f^{\mathfrak
  m,n}(\bar e)\rfloor=0$. But this can only happen if $e\in E^{\mathfrak
  o}$ and $D(z_1^{\ol e,n})=1$. At any rate, the statement of the
proposition holds.

If $v=w$ then 
$$
\O_{\wt{\mathcal X}^n}\big(\,\wt f\,\big)|_{X_v}\cong\O_{X_v}\Big(
-\sum N_e\Big),
$$
where the sum is over the $e\in\E$ such that $\te_e=v$ and $\m_e=0$. 
For such $e$ we have that $\mathfrak d_f^{\mathfrak m,n}(e)=-1$, thus 
$N_e$ is in the support of $U$ with multiplicity -1. 
Furthermore, suppose $N_e$ is in the support of $U$
for a certain $e\in\E$. As before, we may suppose $\te_e=v$. If
$\m_e>0$ then, since $|\m_{\bar e}|<\ell_e$, we have  $\lfloor\mathfrak d_f^{\mathfrak
  m,n}(\bar e)\rfloor=0$. Since $N_e$ is in the support of $U$ we must
then have $\m_e\leq 0$. But if $\m_e<0$ then $\lfloor\mathfrak d_f^{\mathfrak
  m,n}(\bar e)\rfloor=-1$. Thus, since $N_e$ is in the support of $U$
we must have $\m_e=0$. It follows that the statement of the
proposition holds.

In the general case, we may assume that $\deg(f)>0$ and let $w\in
V(G)$ 
such that $f(w)>0$. Let
$g:=f-\chi_{\indm w}$. Then
$$
\O_{\wt{\mathcal X}^n}\big(\,\wt f\,\big)\cong 
\O_{\wt{\mathcal X}^n}\big(\,\wt g\,\big)\otimes\O_{\wt{\mathcal
    X}^n}\big(\,\wt h\,\big)
$$
where $\wt g\in C^0(H^n,\Z)$ is the canonical extension of $g$ with respect to
$D$, where $D'=D+\text{div}(\wt g)$, and where $\wt h\in C^0(H^n,\Z)$ is the
canonical extension of $\chi_w$ with respect to $D'$. The isomorphism
holds by~\cite[Prop.~2.11]{AE1}. Now, by induction,
$$
\O_{\wt{\mathcal X}^n}\big(\,\wt g\,\big)|_{X_v}\cong \O_{X_v}
\Big(\sum_{\substack{e\in\E\\ \te_e=v}}
\lfloor \mathfrak d_g^{\mathfrak m,n}(e)\rfloor N_e
+\sum_{\substack{e\in\E\\ \te_{e}=v\\ \mathfrak m_e<0}} N_e\Big).
$$
Thus we need only prove that
$$
\O_{\wt{\mathcal
    X}^n}\big(\,\wt h\,\big)|_{X_v}\cong\O_{X_v}
\Big(\sum_{\substack{e\in\E\\ \te_e=v}}
\big(\lfloor \mathfrak d_f^{\mathfrak m,n}(e)\rfloor N_e-
\lfloor \mathfrak d_g^{\mathfrak m,n}(e)\rfloor\big) N_e\Big),
$$
or, using what we have proved above, that for each $e\in\E$ with
$\te_e=v$ we have that 
\begin{equation}\label{TFDU}
\lfloor \mathfrak d_f^{\mathfrak m,n}(e)\rfloor-
\lfloor \mathfrak d_g^{\mathfrak
  m,n}(e)\rfloor=
\begin{cases}1&\text{if $\he_e=w$ and either $n\ell_e=1$ or $D'(z_1^{e,n})=1$;}\\
-1&\text{if $v=w$ and $\m'_e=0$;}\\
0&\text{otherwise,}
\end{cases}
\end{equation}
where $\m'\in C^1(G,\Q)$ satisfies
$$
\mathfrak m'_e:=\frac{1}{n}\sum_{i=1}^{n\ell_e-1}iD'(z^{\ol e,n}_i)
$$
for each $e\in E^{\mathfrak o}$.

Indeed, if $v\neq w$, then both sides in~\eqref{TFDU} are zero unless
$\he_e=w$. Suppose $e=vw$. Then the left-hand side in~\eqref{TFDU} is 1
if $\mathfrak d_f^{\mathfrak m,n}(e)\in\Z$ and 0 otherwise. But 
$\mathfrak d_f^{\mathfrak m,n}(e)\in\Z$ if and only if either $n\ell_e=1$ or
$D'(z_1^{e,n})=1$. 

If $v=w$ then the left-hand side in~\eqref{TFDU} is $-1$
if $\mathfrak d_g^{\mathfrak m,n}(e)\in\Z$ and 0 otherwise. But 
$\mathfrak d_g^{\mathfrak m,n}(e)\in\Z$ if and only if $\m'_e=0$. 

In any case, Equation~\eqref{TFDU} follows.
\end{proof}

Recall the definitions of $\mathbf R$ and the
subschemes $Y_{\ell,\m}^{a,b}\subseteq\mathbf R$ from \cite{AE2},
which we review now.

First, as seen in Section~\ref{embsheaves}, 
$$
\mathbf R:=\prod_{e\in E^{\mathfrak o}}\mathbf R_e,
$$
where $\mathbf R_e$ is the doubly infinite chain of smooth rational
curves. We may order the rational curves in the chain, and denote them by $\mathbf P_{e,i}$ for
$i\in\Z$. We may give them coordinates $(\x_{e,i}:\x_{\ol e,i})$ such
that the point $0_{e,i}$, given by $\x_{e,i}=0$, is the point of
intersection of $\mathbf P_{e,i}$ with $\mathbf P_{e,i+1}$, and 
the point $\infty_{e,i}$, given by $\x_{\ol e,i}=0$, is the point of
intersection of $\mathbf P_{e,i}$ with $\mathbf P_{e,i-1}$. We can
also define $\P_{e,i}$ for each $i\in\R-\Z$ as the point of
intersection of $\P_{e,\lfloor i\rfloor}$ with $\P_{e,\lceil i\rceil}$.
Then 
$$
\mathbf R=\bigcup_{\alpha\in C^1(G,\Z)}\mathbf P_{\alpha},
$$
where, more generally,
$$
\mathbf P_\alpha:=\prod_{e\in E^{\mathfrak
    o}}\P_{e,\alpha_e}\quad\text{for each $\alpha\in
C^1\big(G,\frac{1}{2}\Z\big)$.}
$$
We see that $\mathbf P_{\alpha}\supseteq\mathbf P_{\beta}$ if and only
if $|\beta_e-\alpha_e|\leq\frac{1}{2}$ for each $e\in\E$, with
$\alpha_e\in\Z$ whenever $\beta_e\in\Z$. Removing from each $\mathbf
P_{\alpha}$ all those $\mathbf P_{\beta}$ contained in it, we obtain
the \emph{interior} of $\mathbf
P_{\alpha}$, the open subscheme denoted $\mathbf P_{\alpha}^*$. The
$\mathbf P_{\alpha}^*$ for $\alpha\in C^1(G,\frac{1}{2}\Z)$ form a
stratification of $\mathbf R$.

Let $a\: C^1(G,\Z)\to\Gm(\k)$ and $b\: C^1(G,\Z)\to\Gm(\k)$ be
characters. Let $\ell\: E\to \N$ be an edge length function and
$\mathfrak m\in C^1(G,\Z)$. To each $f \in C^0(G, \mathbb Z)$ we
associate the subvariety $\P_{\ell,\m,f}^{a,b}$ of $\P_{\dl^{\m}_f}$ given by the equations  
\[
\forall \textrm{ oriented cycle } \gamma \textrm { in $G^\m_f$},
\]
\[
\prod_{e\in\bar\gamma\cap E^{\mathfrak o}}b_ea_e^{\dl^\m_f(e)}\prod_{e\in \gamma\cap E^{\mathfrak o}} \x_{e,\dl^\m_f(e)} \prod_{e\in\bar\gamma\cap E^{\mathfrak o}}\x_{\ol e,\dl^\m_f(e)}
=\prod_{e\in \gamma\cap E^{\mathfrak o}} b_ea_e^{\dl^\m_f(e)} \prod_{e\in \bar \gamma\cap E^{\mathfrak o}} \x_{e,\dl^\m_f(e)}\prod_{e\in \gamma\cap E^{\mathfrak o}} \x_{\ol e,\dl^\m_f(e)},
\]
where $G^\m_f$ is the spanning subgraph of $G$ whose edges are those
of $G$ for which $\dl^\m_f$ is an integer. 
The equation corresponding
to $\gamma$ may also be written in the format:
\begin{equation}\label{gammapre}
\prod_{e\in\gamma}\x_{e,\mathfrak d^{\mathfrak m}_f(e^{\mathfrak o})}
=\prod_{e\in \gamma} b_ea_e^{\mathfrak d^{\mathfrak m}_f(e^{\mathfrak o})} 
\prod_{e\in \bar \gamma} \x_{e, \mathfrak d^{\mathfrak m}_f(e^{\mathfrak o})},
\end{equation}
where $e^{\mathfrak o}:=e$ if $e\in E^{\mathfrak o}$ and
$e^{\mathfrak o}:=\ol e$ otherwise.

It is clear from \eqref{gammapre} that $\P_{\ell,\m,f}^{a,b}$ depends
rather on the restriction of $b$ to $H^1(G,\Z)$. As any character of
$H^1(G,\Z)$ extends to one of $C^1(G,\Z)$, we will assume later that
$b$ is rather a character of $H^1(G,\Z)$.

We denote by $Y_{\ell,\m}^{a,b}$ the union of the
$\P_{\ell,\m,f}^{a,b}$ for all $f\in C^0(G,\Z)$.

There is a natural action of the character group of $C^1(G,\Z)$, which
we denote by $\mathbf G^E_{\mathbf m}$, on $\mathbf R$: To a character
$c\:C^1(G,\Z)\to\Gm(\k)$ and a point $p$ on 
$\P_\alpha$ with coordinates $(x_{e,\alpha_e},x_{\ol e,\alpha_e})$, for
each $e\in E^{\mathfrak o}$, we associate the point on the same
$\P_\alpha$ with coordinates $(c_ex_{e,\alpha_e},x_{\ol
  e,\alpha_e})$. The homomorphism $d^*$ induces a homomorphism
from the character group of $C^0(G,\Z)$, which we denote 
$\mathbf G^{V}_{\mathbf m}$, to $\mathbf G^E_{\mathbf m}$. The
  induced action on $\mathbf R$ of $c\in \mathbf G^{V}_{\mathbf m}$
    takes the point $p$ to that on the same
$\P_\alpha$ with coordinates 
$(c_vx_{e,\alpha_e},c_ux_{\ol e,\alpha_e})$ for each $e=uv\in
E^{\mathfrak o}$. Finally, the degree map induces a ``diagonal''
injective homomorphism 
$\mathbf G_{\mathbf m}\hookrightarrow\mathbf G^{V}_{\mathbf m}$, and
the induced action of $\mathbf G_{\mathbf m}$ on $\mathbf R$ is
trivial. We may thus
  speak of the action of the quotient $\mathbf G^{V}_{\mathbf m}/\Gm$ 
on $\mathbf R$. 

It is clear from
Equations~\eqref{gammapre} that the action of $\mathbf G^{V}_{\mathbf m}/\Gm$ 
on $\mathbf R$ leaves each $Y_{\ell,\m}^{a,b}$ invariant. We may thus
describe each $Y_{\ell,\m}^{a,b}$ in terms of its orbits. This was
done in \cite{AE2}, Thm.~5.3, and we reproduce it and its Cor.~5.4 here for later
use.

\begin{thm}\label{unionorbits} 
Let $\ell\: E\to\mathbb N$ be an edge length function,
  $\mathfrak m\in C^1(G,\Z)$, and let
$$
a\: C^1(G,\Z)\to\mathbf G_{\mathbf m}(\k),\quad
b\: C^1(G,\Z)\to \mathbf G_{\mathbf m}(\k)
$$
be characters. For each $n\in\mathbb N$ and $f\in C^0(G,\Z)$,
let $p^n_{f}$ be the point on $\mathbf P_{\mathfrak
  d^{\mathfrak m,n}_f}$ given by the coordinates 
$$
(b_ea_e^{\mathfrak d^{\mathfrak m,n}_f(e)}:1) \quad\text{for each }
e\in E^{\mathfrak o}\text{ with }\mathfrak d^{\mathfrak
  m,n}_f(e)\in\Z.
$$
Then $Y_{\ell,\m}^{a,b}$
is the union of the orbits of the $p^n_f$ under the action of 
$\mathbf G^{V}_{\mathbf m}/\mathbf G_{\mathbf m}$. Furthermore, given 
$\alpha\in C^1(G,\frac 12 \Z)$, we have that
  $\P_\alpha^*\cap Y^{a,b}_{\ell,\m}\neq\emptyset$ if and only if
  $\alpha=\dl^{\m,n}_f$ for certain $n\in\mathbb N$ and $f\in C^0(G,\Z)$, and in this case
  $\P_\alpha^*\cap Y^{a,b}_{\ell,\m}$ is the orbit of $p^n_f$ and is 
dense in $\P_\alpha\cap Y^{a,b}_{\ell,\m}$.
\end{thm} 

Recall from Section~\ref{embsheaves} 
the definition of the schemes $\IJ^{\bfb}$ and $\IR_{\IJ^{\bfd}}$ and the stacks
$\mathbf S^d$ and $\mathbf J^d$. Recall that given a point $s\in\IJ^{\bfb}$ corresponding 
to a collection of torsion-free, rank-one sheaves $(L_v\,;\, v\in V)$
with $L_v$ of degree $\bfb_v$, fixing trivializations
$L_v|_{N_e}\cong \k$ and $\O_{X_v}(N_e)|_{N_e}\cong \k$ for each $e\in
E$ and each $v\in e$, we obtain an isomorphism between the fiber of
$\IR^E_{\IJ^{\bfb}}/\IJ^{\bfb}$ over $s$ and $\IR$. The isomorphism is
uniquely defined once we specify that $\P^{\bfc+c}_{\IJ^{\bfb}}$ is taken to
  $\P_{c}$ for each $c\in C^1(G,\Z)$. For simplification, we assume
from now on that $\bfc=0$. 

The importance of the $Y^{a,b}_{\ell,\m}$ stems from the following proposition.

\begin{prop}\label{corresp} Let $m\in\N$ and $L_\eta$ be an invertible sheaf of
degree $d$ on the 
generic fiber of $\X^m/B$. Let $\L$ be a $\sigma^m$-admissible 
almost invertible extension of $L_\eta$ to $\wt\X^m$. For each
$n\in\N$ divisible by $m$, let $\L^n$ be the pullback of $\L$ 
to $\wt\X^n$, and for each $f\in C^0(G,\Z)$ let 
$$
L^n(f):=\mathcal L^n(f)|_{X^{n\ell}}\quad\text{and}\quad
I^n(f):=\sigma^{n\ell}_*L^n(f).
$$
Let $D^m\in\text{\rm Div}(H^m)$ be the divisor associated to $L^m(0)$. For
each $e\in E^{\mathfrak o}$, let 
$$
\mathfrak m_e:=\frac{1}{m}\sum_{i=1}^{m\ell_e-1}iD^m(z^{\ol e,m}_i),
$$
and let $\m$ be the
unique element of $C^1(G,\Q)$ satisfying this. Let $E_0\subseteq E$ be
the support of $\m$ and put
$$
K_v:=I^m(0)_v\Big(\sum_{\substack{e\in E^{\mathfrak o}_0\\
    \he_e=v}}N_e\Big)\quad\text{for each $v\in V(G)$},
$$
where $I^m(0)_v$ is the torsion-free, rank-one sheaf generated by
$I^m(0)$ on $X_v$. Let $\mathbf b\in C^0(G,\Z)$ satisfying $\mathbf
b_v=\deg (K_v)$ for each $v\in V(G)$. Fix trivializations
$K_v|_{N_e}\cong \k$ and $\O_{X_v}(N_e)|_{N_e}\cong \k$ for
each $e\in E$ and $v\in e$. Then the following three
statements hold:
\begin{enumerate}
\item The degree of $\bfb$ is $d$ and the $K_v$ are represented by a unique
point $s\in\mathbf J^{\mathbf b}$. 
\item For each $n\in\mathbb
N$ with $m|n$ and $f\in C^0(G, \mathbb Z)$, there is a unique 
(modulo action of $\mathbf G_\mathbf m^V/\Gm$) point $t^n_f$ on the fiber of $\IR_{\IJ^{\bfd}}$
over $s$ representing $I^n(f)$ with $c(t^n_f)=\mathfrak d^{\mathfrak
  m,n}_f$.
\item There are characters 
$$
a\: C^1(G,\Z)\to\mathbf G_{\mathbf m}(\k)\quad\text{and}\quad
b\: C^1(G,\Z)\to \mathbf G_{\mathbf m}(\k)
$$
such that, under the chosen trivializations, the equivariant isomorphism of the fiber of
  $\IR_{\IJ^{\bfb}}$ over $s$ with $\mathbf R$, taking $\P^{c}_{\IJ^{\bfb}}$ to
  $\P_c$ for each $c\in C^1(G,\Z)$, takes a point on the
  orbit of $t^n_f$ for each $n\in\N$ divisible by $m$ and $f\in C^0(G,\Z)$ to
  the point on $\mathbf P_{\mathfrak
  d^{\mathfrak m,n}_f}$ given by the coordinates 
$$
(b_ea_e^{\mathfrak d^{\mathfrak m,n}_f(e)}:1) \quad\text{for each }
e\in E^{\mathfrak o}\text{ with }\mathfrak d^{\mathfrak
  m,n}_f(e)\in\Z.
$$
\item In particular, the union of the orbits of all the $t^n_f$ is a closed subset of the
fiber of $\IR_{\IJ^{\bfb}}$ over $s$ isomorphic to $Y^{a,b}_{m\ell,m\m}$.
\end{enumerate}
\end{prop}

\begin{proof} The first statement is immediate from the fact that
  $\deg\mathbf b=\deg I^m(0)=\deg L_\eta=d$.

\smallskip

As for the second statement, let $n\in\N$ divisible by $m$ and $f\in
C^0(G,\Z)$. Let $D^n$ denote the pullback of $D^m$ to $H^n$: we
have $D^n(v)=D^m(v)$ for each $v\in V(G)$ and $D^n(z_j^{\ol e,n})=0$ for
each $e\in\E$ and $j=1,\dots,n\ell_e-1$, unless
$n|jm$ and $D^m(z_{jm/n}^{\ol e,n})=1$, in which case $D^n(z_j^{\ol e,n})=1$. 
Also, $D^n$ is the 
associated divisor to $L^n(0)$. In particular, $D^n$ is
$G$-admissible. Notice that 
$$
\m_e=\frac{1}{n}\sum_{i=1}^{m\ell_e-1}\frac{n}{m}iD^n(z_{ni/m}^{\ol e, n})=
\frac{1}{n}\sum_{i=1}^{n\ell_e-1}iD^n(z_i^{\ol e, n})
$$
for each $e\in E^{\mathfrak o}$. Then, by Propositions~\ref{ILv}~and ~\ref{IfI0},
$$
I^n(f)_v\cong L^n(f)|_{X_v}\cong L^m(0)|_{X_v}\otimes\O_{X_v}
\Big(\sum_{\substack{e\in\E\\ \te_e=v}}
\lfloor \mathfrak d_f^{\mathfrak m,n}(e)\rfloor N_e
+\sum_{\substack{e\in\E\\ \te_{e}=v\\ \mathfrak m_e<0}} N_e\Big)
$$
for each vertex $v\in V(G)$, where $I^n(f)_v$ is the restriction modulo torsion of
$I^n(f)$ to $X_v$. Furthermore, since $\mathfrak m_e>0$ if and only if $e\in E^{\mathfrak
  o}_0$, we have
$$
I^n(f)_v\cong I^m(0)_v\Big(\sum_{\substack{e\in E^{\mathfrak o}_0\\
    \he_e=v}} N_e\Big)\otimes\O_{X_v}
\Big(\sum_{\substack{e\in\E\\ \he_e=v}}
\lfloor \mathfrak d_f^{\mathfrak m,n}(\bar e)\rfloor N_e\Big)\cong
K_v\Big(\sum_{\substack{e\in\E\\ \he_e=v}}
\lfloor -\mathfrak d_f^{\mathfrak m,n}(e)\rfloor N_e\Big)\cong
\L^{\mathfrak d_f^{\mathfrak m,n}}_v(s),
$$
the latter isomorphism following from the definition in
Subsection~\ref{atlas}. In addition, since $I^n(f)=\sigma^{n\ell}_*L^n(f)$, 
it follows that $I^n_f$ fails to be invertible at a node
$N_e$ if and only if $\dl^{\m,n}_f(e)\not\in\Z$. 
The statement follows now from Proposition \ref{embshfc}.

\smallskip

We prove now Statement (3). First, given the trivializations, the $\mathcal O_{\mathcal
  X^n}(f)|_X$ depend on $a\: C^1(G,\Z)\to{\mathbf G}_{\mathbf m}(\k)$ arising
from the deformation $\mathcal X/B$, as explained in
Subsection~\ref{ges}. More explicitly, $\mathcal O_{\mathcal
  X^n}(f)|_X$ is isomorphic to the subsheaf of 
$$
\bigoplus_{v\in V(G)}\mathcal O_{X_v}\Big(
\sum_{\substack{e\in\E\\ \he_e=v}}
\dl^{0,n}_f(\bar e)N_e\Big)
$$
whose cokernel is supported at the nodes $N_e$ for $e=uv\in E^{\mathfrak
  o}$ such that $\dl^{0,n}_f(e)\in\Z$ and is equal in a neighborhood
of such a $N_e$ to the quotient of the vector space 
$$
\mathcal O_{X_u}\big(\dl^{0,n}_f(e)N_e\big)|_{N_e}\oplus
\mathcal O_{X_v}\big(-\dl^{0,n}_f(e)N_e\big)|_{N_e}
$$
by a certain one-dimensional
vector subspace. The space is identified with $\k\oplus \k$,
under our choices of trivializations, and the subspace is that 
generated by $(a^{\dl^{0,n}_f(e)},1)$.

Second, for each $e=uv\in E^{\mathfrak o}$, consider the sheaf
$I^m(m\mathfrak m_e\chi_{\indm u})$. Then 
$\dl^{\m,m}_{m\m_e\chi_u}(e)\in\Z$, whence $I^m(m\mathfrak m_e\chi_{\indm u})$
is invertible at $N_e$. Furthermore, it is naturally in a neighborhood
of $N_e$ a subsheaf of $K_u\oplus K_v$ whose quotient is the same as 
the quotient of a one-dimensional vector subspace 
of $K_u|_{N_e}\oplus K_v|_{N_e}$. The latter is identified with $\k\oplus \k$,
under our choices of trivializations. Define 
$b_e\in \k$ such that the
one-dimensional vector subspace is generated by $(b_e,1)$. Thus, we have an element $b\:C^1(G,\Z)\to{\mathbf G}_{\mathbf m}(\k)$. Its
restriction to $H^1(G,\Z)$ is denoted by $b$ by abusing the notation. 

It follows that, under the above trivializations, 
given $n\in\mathbb N$, $f\in C^0(G,\Z)$ and $c\in C^1(G,\Z)$ such
that $|c_e-\mathfrak d^{\mathfrak m,n}_f(e)|\leq 1/2$ for all $e\in\mathbb
E$, the fiber of $\IP^{c}_{\IJ^{\bfb}}$ over $s$ is identified with
$\mathbf P_{c}$, and the point $t^n_f$ on the fiber 
representing $I^n(f)$ corresponds to the point on
$\mathbf P_{\mathfrak d^{\mathfrak m,n}_f}\subseteq \mathbf P_{c}$
with coordinates
$$
(b_ea_e^{c_e}:1)\quad\text{for each }e\in E^{\mathfrak o}\text{ with }
\mathfrak d^{\mathfrak m,n}_f(e)\in\Z.
$$
Indeed, for each $e=uv\in E^{\mathfrak o}$ with $\mathfrak
d^{\mathfrak m,n}_f(e)\in\Z$ we have that $d^{\mathfrak
  m,n}_f(e)=\dl^{0,n}_g(e)$, where $g:=f-n\mathfrak
m_e\chi_u$. Furthermore, $I^n(f)$ is equal to the tensor product of $I^n(n\mathfrak
m_e\chi_u)$ and $\mathcal O_{\mathcal X^n}(g)|_X$ in a neighborhood
of $N_e$, thus equal in that neighborhood to the subsheaf of 
$$
K_u\Big(\dl^{\mathfrak m,n}_f(e)N_e\Big)\oplus K_v\Big(-\dl^{\mathfrak
  m,n}_f(e)N_e\Big)
$$
whose quotient is, under the identification with $\k\oplus\k$ of the
restriction of the latter sheaf to $N_e$, that of $\k\oplus\k$ by
the subspace generated by $(b_ea_e^{\dl^{\m,n}_f(e)},1)$.

The final statement is now an application of 
Theorem~\ref{unionorbits}, once we observe that 
$$
\dl^{\m,n}_{\ell,f}=\dl^{m\m,n/m}_{m\ell,f}
$$
for each $f\in C^0(G,\Z)$ and each $n\in\N$ divisible by $m$.
\end{proof}

\begin{cor} Notations as in Proposition~\ref{corresp}, the sheaf
  $I^n(f)$ is a flat degeneration of the sheaf $I^q(h)$ if
  $\P_{\dl^{\m,n}_f}\subseteq\P_{\dl^{\m,q}_h}$, or equivalently, for each
  $e\in\E$ the following two conditions are satisfied:
\begin{enumerate}
\item $|\dl^{\m,n}_f(e)-\dl^{\m,q}_h(e)|\leq 1/2$.
\item If $\dl^{\m,n}_f(e)\in\Z$ then $\dl^{\m,q}_h(e)\in\Z$.
\end{enumerate}
\end{cor}

\begin{proof} By Proposition~\ref{corresp}, there is a closed subscheme
  $Y$ of the fiber of $\IR_{\IJ^{\bfb}}$ over $s$, with the induced reduced
  structure, which is the union of all the orbits of the $t^n_f$. 
Consider the quotient stack $\mathfrak I:=[Y/(\mathbf G_{\mathbf
    m}^V/\Gm)]$. Let $s^n_f$ denote the
  image of $t^n_f$ in $\mathfrak I$ for each $n$ and $f$. As observed in
  Section~\ref{embsheaves}, there is a relative
  torsion-free rank-one sheaf $\mathcal I$ on $X\times\mathfrak
  I/\mathfrak I$ such that
  $\mathcal I|_{X\times\{s^n_f\}}\cong I^n(f)$ for each $n$ and $f$. 

By Proposition~\ref{corresp}, there is an equivariant isomorphism from $Y$ to 
$Y^{a,b}_{m\ell,m\m}$ for certain $\m$, $a$ and $b$, taking each $t^n_f$
to the point $p^n_f$ on $\mathbf P_{\dl^{\m,n}_f}^*$ given in 
Theorem~\ref{unionorbits}. The isomorphism is equivariant, whence we may view 
$\mathfrak I$ as a quotient of $Y^{a,b}_{m\ell,m\m}$ where each $s^n_f$
is the image of the point $p^n_f$. Pulling back
$\mathcal I$ we obtain a relative torsion-free rank-one sheaf
$\mathcal L$ on $X\times Y^{a,b}_{m\ell,m\m}/Y^{a,b}_{m\ell,m\m}$ such that 
$\mathcal L|_{X\times\{y\}}\cong I^n(f)$ for each point $y$ on the orbit
of $p^n_f$. 

Let $n,q\in\N$ divisible by $m$ and $f,h\in C^0(G,\Z)$. Suppose 
$\mathbf P_{\dl^{\m,q}_h}\supseteq\mathbf P_{\dl^{\m,n}_f}$. It is now
enough to observe that $\mathbf P_{\dl^{\m,q}_h}^*\cap Y^{a,b}_{m\ell,m\m}$ is
the orbit of $p^q_h$ and is dense in $\mathbf P_{\dl^{\m,q}_h}\cap
Y^{a,b}_{n\ell,m\m}$ by Theorem~\ref{unionorbits}.
\end{proof}

Now we need only one technical statement before stating the main result of the
section.

\begin{prop}\label{f1f2} Let $f_1,f_2\in C^0(G,\Z)$ and $n_1,n_2\in\N$. 
If $\mathfrak
  d^{\mathfrak m,n_2}_{f_2}-\mathfrak
  d^{\mathfrak m,n_1}_{f_1}\in H^1(G,\Z)$, then $\mathfrak
  d^{\mathfrak m,n_2}_{f_2}=\mathfrak
  d^{\mathfrak m,n_1}_{f_1}$.
\end{prop}

\begin{proof} Since $\mathfrak d^{\mathfrak m,n_i}_{f_i}= \mathfrak d^{\mathfrak
    m,pn_i}_{pf_i}$ for each $p\in\N$, we may assume $n_1=n_2$. Then
  the proposition follows from \cite{AE2}, Prop.~5.1.
\end{proof}

\begin{thm}[Degeneration]\label{thm:main4} Let $\mathfrak
  I$ be the collection of stable limits of an invertible sheaf
  $L_\eta$ of degree $d$ on
  the generic fiber of $\X^m/B$ for given integer $m\in\N$ and
  smoothing $\X/B$ with singularity degrees
  $\ell\:E\to\N$ of the nodal curve $X$. Then $\mathfrak
  I$ is
  the support of a reduced closed 
substack of $\IJ^d$, also denoted
  $\mathfrak I$. The inverse image of $\mathfrak I$ in
  $\IR_{\IJ^{\bfd}}$ is an infinite disjoint union of connected closed
  subschemes of certain fibers of $\IR_{\IJ^{\bfd}}$ over
  $\IJ^{\bfd}$. Each connected component is the image of each other,
  under the action of a unique element of $H^1(G,\Z)$, and each maps
  isomorphically to the same closed subscheme $\mathfrak
  Y\subseteq\IS^d$. Furthermore, each connected component, under an
  identification of the fiber containing it with $\IR$, is equal to the subscheme
  $Y_{\ell,\m}^{a,b}$ for certain choices of $\m\in C^1(G,\Z)$ and characters
$a\:C^1(G,\Z)\to\Gm(\k)$ and $b\:H^1(G,\Z)\to\Gm(\k)$. In particular,  
  $\mathfrak Y$ has pure dimension $|V|-1$ and $\mathfrak I$ has pure
  dimension $0$.
\end{thm}

\begin{proof} Using the notation in Proposition~\ref{corresp}, we have
  that $\mathfrak I$ is the collection of the $I^n_f$ for $n\in\N$
  divisible by $m$ and $f\in C^0(G,\Z)$. Let the $K_v$ and $\mathbf b$
  be as in Proposition~\ref{corresp}, as well as $s\in\mathbf
  J^{\bfb}$ and the $t^n_f\in\IR_{\IJ^{\bfb}}$ in Statements (1) and (2)
  therein. 

Let $Y\subseteq \IR_{\IJ^{\bfd}}$ be the union of the orbits of the
$t^n_f$ for all $n\in\N$ and $f\in C^0(G,\Z)$. By
Proposition~\ref{corresp}, it is a closed subset of the fiber of
$\IR_{\IJ^{\bfd}}$ over $s$. We give it
the reduced structure, so that $Y$ becomes isomorphic to
$Y^{a,b}_{\ell,\m}$ for certain choices of $\m\in C^1(G,\Z)$ and characters
$a\:C^1(G,\Z)\to\Gm(\k)$ and $b\:H^1(G,\Z)\to\Gm(\k)$, according to
Proposition~\ref{corresp}. Then $Y$ is a connected closed 
subscheme of $\IR_{\IJ^{\bfd}}$ of pure
dimension $|V|-1$.

We claim the quotient map induces a closed embedding $Y\to\IS^d$. 
Since $\IS^d$ is the quotient of $\IR_{\IJ^{\bfd}}$ by the free action
of the discrete group $H^1(G,\Z)$, 
to prove our claim it is enough to prove that the translate of $Y$ by
any nonzero element 
$\gamma\in H^1(G,\Z)$ does not intersect $Y$. 

So, let $\gamma\in H^1(G,\Z)$. Suppose that there are two points $t_1$
and $t_2$ on $Y$ that differ by the action of $\gamma$. Let $c_1,c_2\in
C^1(G,\Z)$ 
such that $t_1$ lies on the fiber of $\IP^{c_1}_{\IJ^{\bfb}}$ over
$s$ and $t_2$ on the fiber of $\IP^{c_2}_{\IJ^{\bfb}}$ over $s$. 
We may assume that $c_1=c_2+\gamma$, and thus
$c(t_1)=c(t_2)+\gamma$. But, by construction of $Y$, 
there are $n_1,n_2\in\mathbb N$ and 
$f_1,f_2\in C^0(G,\Z)$ such that 
$c(t_i)=\mathfrak d_{f_i}^{\mathfrak m,n_i}$ for $i=1,2$. 
Then Proposition~\ref{f1f2} implies that $\gamma=0$.

Any two points on $\IR_{\IJ^{\bfd}}$ representing $I^n(f)$ for each
$n\in\mathbb N$ and $f\in C^0(G,\Z)$ differ by the action of $H^1(G,\Z)\times \mathbf
G_\mathbf m^V/\Gm$. Since the action of $H^1(G,\Z)$ commutes with that
of $\mathbf G_\mathbf m^V/\Gm$, it follows that all the points on $\IR_{\IJ^{\bfd}}$ representing
$I^n(f)$ are in $\wt\tau^\gamma(Y)$ for
$\gamma\in H^1(G,\Z)$. So the inverse image of $\mathfrak I$ in
$\IR_{\IJ^{\bfd}}$ is the union $\bigcup \wt\tau^\gamma(Y)$ for
$\gamma$ running in $H^1(G,\Z)$. Now, it follows from what we proved
above that the 
intersection of $\wt\tau^\gamma(Y)$ with  $\wt\tau^{\gamma'}(Y)$
for distinct $\gamma,\gamma'\in H^1(G,\Z)$ is empty. Thus the
$\wt\tau^\gamma(Y)$ are the connected components of the inverse image.

It follows that each of the $\wt\tau^\gamma(Y)$ is mapped
isomorphically to the same closed subscheme $\mathfrak
  Y\subseteq\IS^d$, and that $\mathfrak Y$ has pure dimension
  $|V|-1$. Then $\mathfrak I$ has the structure of a closed substack
  of $\mathbf J^d$ of pure dimension 0. All the statements have been proved.
\end{proof}

\begin{remark}\rm Each $\mathfrak I$ should be parameterized by a
  point on $\text{\rm Hilb}_{\IJ^d}$, the Hilbert stack of $\IJ^d$,
  and each $\mathfrak Y$ should be parameterized by a point on 
$\text{\rm Hilb}^{\mathbf G_\mathbf m^V/\Gm}_{\IS^d}$, the
Hilbert stack parameterizing $\mathbf G_\mathbf m^V/\Gm$-invariant closed 
substacks of $\IS^d$, if this Hilbert stack exists! Moreover, the
analysis of a few examples leads us to ask whether these points, as
$\pi$ and $L_\eta$ vary, form a closed substack of
$\text{\rm Hilb}_{\IJ^d}$ or $\text{Hilb}^{\mathbf G_\mathbf
  m^V/\Gm}_{\IS^d}$. 
This substack should be
regarded as a \emph{new compactified Jacobian} for $X$. 

The situation is simpler when $X$ is of compact type. Then there are
no cycles, so $\IR_{\IJ^{\bfd}}=\IS^d$. Furthermore, $\mathfrak Y$ is a fiber
of $\IR_{\IJ^{\bfd}}$ over 
$\IJ^{\bfd}$. Thus the ``new compactified Jacobian'' is $\IJ^{\bfd}$,
which is of course nothing new in this case.  The general situation, where the ``new compactified Jacobian'' is
indeed a new construction, will be addressed in a future work. 
\end{remark}

\begin{remark}\rm There is a natural $\mathbf G_\mathbf m^V/\mathbf
  G_\mathbf m$-invariant relative torsion-free,
rank-one, degree-$d$ sheaf on $X\times \mathfrak Y/\mathfrak Y$ which
is given by the restriction of the sheaf on
$X\times\IS^d$ we called $\I$
in Section~\ref{embsheaves}. Being invariant, it is the pullback of a
relative torsion-free, rank-one sheaf on $X\times\mathfrak I/\mathfrak I$. It is this sheaf that
gives the structure of $\mathfrak I$ as a closed substack of $\IJ^d$.
\end{remark}

  
\section{Regeneration}\label{sec:regeneration}

Recall from Section~\ref{stackydeg} the definition of $\IR$ and of the
$Y^{a,b}_{\ell,\m}$. Recall from Section~\ref{embsheaves} 
the definition of the schemes $\IJ^{\bfb}$ and $\IR_{\IJ^{\bfd}}$ and the stacks
$\mathbf S^d$ and $\mathbf J^d$. Recall that given a point $s\in\IJ^{\bfb}$ corresponding 
to a collection of torsion-free, rank-one sheaves $(L_v\,;\, v\in V)$
with $L_v$ of degree $\bfb_v$, fixing trivializations
$L_v|_{N_e}\cong \k$ and $\O_{X_v}(N_e)|_{N_e}\cong \k$ for each $e\in
E$ and each $v\in e$, we obtain an isomorphism between the fiber of
$\IR^E_{\IJ^{\bfb}}/\IJ^{\bfb}$ over $s$ and $\IR$. The isomorphism is
uniquely defined once we specify that $\P^{\bfc+c}_{\IJ^{\bfb}}$ is taken to
  $\P_{c}$ for each $c\in C^1(G,\Z)$. For simplification, we assume
from now on that $\bfc=0$. 

In this section, we give a converse to Theorem~\ref{thm:main4}. 

\begin{thm}[Regeneration]\label{regeneration}
Let $\ell\:E\to\mathbb N$ be a length function, 
$\mathfrak m\in C^1(G,\Z)$ and let 
$$
a\:C^1(G,\Z)\to{\mathbf G }_{\mathbf m}(\k)\quad\text{and}\quad 
b\:H^1(G,\Z)\to {\mathbf G }_{\mathbf m}(\k)
$$
be characters. Let $Y_{\ell,\mathfrak m}^{a,b}$ be the 
corresponding subscheme of $\IR$. Let $\bfb\in C^0(G,\Z)$. For each $v\in V$, let
$K_v$ be a torsion-free, rank-one sheaf of degree $\bfb_v$ on
$X_v$. Let $s\in\IJ^{\bfb}$ be the corresponding point. 
Fix trivializations $K_v|_{N_e}\cong \k$ and $\O_{X_v}(N_e)|_{N_e}\cong \k$ for
each $e\in E$ and $v\in e$. Let $Y\subseteq
\IR_{\IJ^{\bfb}}$ be the closed subscheme of the fiber over $s$
corresponding to $Y_{\ell,\m}^{a,b}$ under the isomorphism of the fiber
with $\IR$ induced by the trivializations. Let $\mathfrak
I\subseteq\mathbf J^d$ be
the image of $Y$. Then there exist a
smoothing $\pi\:\X\to B$ of $X$ and an invertible sheaf $L_\eta$
on the generic fiber of $\pi$ such that the following two statements
hold:
\begin{enumerate}
\item The total space $\X$ is regular everywhere but possibly at
the nodes $N_e$ of $X$, where it has singularity degree $\ell_e$ for
each $e\in \E(G)$, and at the remaining nodes of $X$ where the sheaf
$\oplus K_v$ is not invertible, where it has singularity degree $2$.
\item $\mathfrak I$ is the
collection of stable limits of $L_\eta$.
\end{enumerate}
\end{thm}

The
whole section is devoted to the proof of the above theorem.

\smallskip

Let $G'$ denote the full dual graph of $X$, containing
  self loops. Consider a versal
deformation $V/M$ of $X$, as explained in Subsection~\ref{ges}. 
More precisely, $M$ is the spectrum of the power series ring $R$ over $\k$ in the
variables $t_e$ for $e\in E(G')$, and variables
$s_1,\dots,s_p$ for a
certain integer $p$, and we have an isomorphism of $R$-algebras
$$
\psi_e\:\widehat{\mathcal O}_{V,N_e}\to R[[z_e,w_e]]/(z_e w_e-t_e)
$$
for each $e\in E(G')$, with $z_e=0$ corresponding to the
component $X_v$ of $X$ and $w_e=0$ to $X_u$ if $e\in E(G)$ and
$e^{\mathfrak o}=uv$. The pullbacks of $z_e$ and $w_e$ under $\psi_e$
yield under restriction analytic
local parameters of $X_u$ and $X_v$ at $N_e$, respectively. We may
assume that the isomorphisms $\mathcal O_{X_v}(N_e)|_{N_e}\cong \k$ are
given by them, for all $e\in E(G)$ and $v\in V(G)$ with $v\in e$.

We may then consider the smoothing  $\pi\:\X\to B$ of $X$ induced by
the map $B\to M$ pulling back the $s_i$ to $t$, and pulling back  
$t_e$ to $a_{e^{\mathfrak o}}t^{\ell_e}$ for each $e\in E(G)$, to $t$
for each $e\in E(G')-E(G)$ if $K_v$ is invertible
at the corresponding node for $v\in e$, and to $t^2$ otherwise. Then
$\pi$ is as in the statement of the theorem. Also, as explained in
Subsection~\ref{ges}, the sheaves $\mathcal O_{\mathcal X^n}(f)|_X$ are
determined by the homomorphism $a$.

Let $\sigma\:\wt\X\to\X$ be the Cartier reduction of $\pi$. The
composition $\wt\pi:=\pi\sigma\:\wt\X\to B$ is a smoothing of
$X^{\ell}$. As in Section~\ref{stackydeg}, we will let $\wt\pi^n\:\wt{\mathcal
  X}^n\to B$ denote the Cartier reduction of the extension of
$\pi$ by the base change map $\mu_n\:B\to B$ for each $n\in\N$; it
is a smoothing of $X^{n\ell}$.

\smallskip

Recall from Theorem~\ref{unionorbits} that $Y_{\ell,\mathfrak
  m}^{a,b}$ is a union of orbits in
$\IR$ under the action of $\mathbf G_\mathbf m^V/\Gm$, the orbits of
the $p^n_f$ for $f\in C^0(G,\Z)$ and $n\in\N$, where $p^n_f$ is a point on
$\P_{\dl^{\m,n}_f}$ whose coordinates are given in terms of $a$ and
$b$, as made precise in that theorem. Now, for each $f\in C^0(G,\Z)$ and $n\in\N$, 
let $I^n_f$ be the torsion-free, rank-one sheaf on
$X$ parameterized by the point $t^n_f\in Y$ corresponding to $p^n_f$. For each $e=uv\in\E$, write
\begin{equation}\label{dfi}
f(v)-f(u)+n\mathfrak m_e=n\ell_e\delta_e^{\m,n}(f)+n\ell_e-i^n_e(f),
\end{equation}
where $i^n_e(f)$ is an integer satisfying $0<i^n_e(f)\leq n\ell_e$. 
We can write $I^n_f=\sigma^{n\ell}_*L^n_f$ for a certain almost 
invertible sheaf $L^n_f$ on $X^{n\ell}$, admissible for $\sigma^{n\ell}$, such
that $L^n_f$ has degree zero on each rational smooth curve $Z^{e,n}_i$ contracted
by $\sigma^{n\ell}$ onto $N_e$ for each $e=uv\in\mathbb E$, unless 
$\mathfrak d_f^{\mathfrak m,n}(e)\not\in\Z$ and $i=i^n_e(f)$.

The sheaf $L^n_f$ is not unique. Its associated divisor in $H^n$,
denoted henceforth by $D^n_f$, is. 
At any rate, with the choices we have already made, we have the
following claim.

\begin{claim}\label{cl1} Notation as above, for each $f,h\in C^0(G,\Z)$ 
and $n\in\N$, the sheaves 
\begin{equation}\label{Lnhf}
L^n_h\quad\text{and}\quad L^n_f\otimes\O_{\wt\X^n}\big(\,\widetilde
g\,\big)|_{X^{n\ell}}
\end{equation}
have the same degrees on the exceptional components of $X^{n\ell}$ and
the same restrictions to $X_v$ for each $v\in V(G)$, where $\widetilde
g$ is the canonical extension of $g:=h-f$ with respect to $D^n_f$. In
particular,
$$
D^n_h=D^n_f+\div(\wt g).
$$
\end{claim}

\begin{proof}
Indeed, both sheaves in \eqref{Lnhf} are
$\sigma^{n\ell}$-admissible. Furthermore, for each $e=uv\in\E$, the
sheaf on the right in \eqref{Lnhf} has degree 0 on $Z^{e,n}_i$ for each $i=1,\dots,n\ell_e-1$,
unless
$$
g(u)-g(v)+i^n_e(f)=n\ell_e\Big\lfloor\frac{g(u)-g(v)+i^n_e(f)}{n\ell_e}\Big\rfloor+i.
$$
Substituting for $i^n_e(f)$ from Equation~\eqref{dfi} we get
$$
g(u)-g(v)+f(u)-f(v)+n\ell_e\delta_e^{\m,n}(f)+n\ell_e-n\mathfrak
m_e=n\ell_e\Big\lfloor\frac{g(u)-g(v)+i^n_e(f)}{n\ell_e}\Big\rfloor+i,
$$
whence
$$
h(v)-h(u)+n\mathfrak
m_e=n\ell_e\Big(\delta_e^{\mathfrak m,n}(f)-
\Big\lfloor\frac{g(u)-g(v)+i^n_e(f)}{n\ell_e}\Big\rfloor\Big)+n\ell_e-i.
$$
Comparing with Equation~\eqref{dfi}, with $f$ replaced by $h$, it
follows that both sheaves in \eqref{Lnhf} have the same degree on
each exceptional component $Z^{e,n}_i$.

Now, given $v\in V(G)$ and $q\in C^0(G,\Z)$, it follows from
Propositions~\ref{cappc12} and \ref{ILv} that 
\begin{equation}\label{Lnh}
L^n_q|_{X_v}\cong(I^n_q)_v\cong K_v\otimes
\mathcal O_{X_v}\Big(\sum_{\substack{e\in\E\\ \te_e=v}}
\lfloor \mathfrak d_q^{\mathfrak m,n}(e)\rfloor N_e\Big),
\end{equation}
where $(I^n_q)_v$ is the restriction of $I^n_q$ to $X_v$ modulo
torsion. In addition, by Proposition~\ref{IfI0}, 
\begin{equation}\label{OXn}
\O_{\wt\X^n}\big(\,\widetilde
g\,\big)|_{X_v}\cong \mathcal O_{X_v}\Big(\sum_{\substack{e\in\E\\ \te_e=v}}
\lfloor \mathfrak d_g^{\mathfrak p,n}(e)\rfloor N_e 
+\sum_{\substack{e\in\E\\ \te_{e}=v\\ \mathfrak p_e<0}} N_e\Big),
\end{equation}
where $\mathfrak p\in C^1(G,\Z)$ satisfies 
\begin{equation}\label{pei}
\mathfrak p_e=\frac{1}{n}\sum_{i=1}^{n\ell_e-1}iD^n_f(z_i^{\ol e,n})
\end{equation}
for each $e\in E^{\mathfrak o}$. 

Since $D^n_f$ is the divisor
associated to $L^n_f$, we have that $n\mathfrak p_e=i^n_{\ol e}(f)$
for each $e\in E^{\mathfrak o}$ with $i^n_{\bar e}(f)\neq n\ell_e$, whence 
$n\mathfrak p_e=n\ell_e-i^n_{e}(f)$. The latter holds even if
$i^n_{\ol e}(f)=n\ell_e$, because then $\mathfrak p_{e}=0$ 
and $i^n_{e}(f)=n\ell_e$. Now, if $e\in\E-E^{\mathfrak o}$,
then $n\mathfrak p_{\ol e}=i^n_{e}(f)$,
whence $n\mathfrak p_e=-i^n_{e}(f)$, unless $i^n_e(f)=n\ell_e$, in
which case $n\mathfrak p_e=0$. To summarize,
\begin{equation}\label{pe0}
n\mathfrak p_{e}=\begin{cases}
n\ell_e-i^n_e(f)&\text{if }\mathfrak p_e\geq 0,\\
-i^n_{e}(f) &\text{if }\mathfrak p_e< 0.
\end{cases}
\end{equation}

Finally, for each $e\in\E$, 
\begin{align*}
n\ell_e\lfloor \mathfrak d_h^{\mathfrak m,n}(e)\rfloor
=&h(\he_e)-h(\te_e)+n\mathfrak m_e-n\ell_e+i^n_{e}(h)\\
=&h(\he_e)-h(\te_e) -f(\he_e)+f(\te_e)-i^n_{e}(f)+n\ell_e\lfloor \mathfrak d_f^{\mathfrak m,n}(e)\rfloor
+i^n_{e}(h)\\
=&g(\he_e)-g(\te_e) +n\mathfrak p_e-n\mathfrak p_e-i^n_{e}(f)
+n\ell_e\lfloor \mathfrak d_f^{\mathfrak m,n}(e)\rfloor
+i^n_{e}(h)\\
=&n\ell_e \lfloor \mathfrak d_g^{\mathfrak p,n}(e)\rfloor+\rho -n\mathfrak p_e-i^n_{e}(f)
+n\ell_e\lfloor \mathfrak d_f^{\mathfrak m,n}(e)\rfloor 
+i^n_{e}(h)
\end{align*}
for a certain integer $\rho$ satisfying $0\leq\rho<n\ell_e$. 
Since $0<\rho+i^n_{e}(h)< 2n\ell_e$, and Equation~\eqref{pe0}
yields that $n\mathfrak p_e+i^n_{e}(f)$ is a multiple of
$n\ell_e$, it follows that 
$\rho+i^n_{e}(h)=n\ell_e$ and 
\begin{equation}\label{ddd}
\lfloor \mathfrak d_h^{\mathfrak m,n}(e)\rfloor=
\begin{cases}
\lfloor \mathfrak d_g^{\mathfrak p,n}(e)\rfloor+
\lfloor \mathfrak d_f^{\mathfrak m,n}(e)\rfloor&\text{if
}\mathfrak p_e\geq 0,\\
\lfloor \mathfrak d_g^{\mathfrak p,n}(e)\rfloor+
\lfloor \mathfrak d_f^{\mathfrak m,n}(e)\rfloor+1&\text{if
}\mathfrak p_e< 0.
\end{cases}
\end{equation}
It follows now from Isomorphisms~\eqref{Lnh} for $q=h$ and $q=f$, and from
Equation~\eqref{OXn}, that the restrictions to
$X_v$ of both sheaves in \eqref{Lnhf} are equal, for each $v\in V(G)$,
finishing the proof of our claim.
\end{proof}

Notice as well that, since $\rho+i^n_{e}(h)=n\ell_e$, if
$\dl_h^{\m,n}(e)\in\Z$, then $i^n_e(h)=n\ell_e$ and thus $\rho=0$, implying
that also $\dl_g^{\mathfrak p,n}(e)\in\Z$. And if
$\dl_f^{\m,n}(e)\in\Z$ then $\mathfrak p_e=0$. Now, if both $I^n_f$ and $I^n_h$ are
invertible at $N_e$, then all of $\dl_f^{\m,n}(e),\dl_h^{\m,n}(e),
\dl_g^{\mathfrak p,n}(e)$ are integers, $\mathfrak p_e=0$ and it
follows from \eqref{ddd} that
\begin{equation}\label{dd0}
\dl_h^{\m,n}(e)=\dl_g^{\mathfrak p,n}(e)+\dl_h^{\m,n}(e).
\end{equation}
Furthermore, since we have already shown that the sheaves in \eqref{Lnhf} have the
same degrees on the exceptional components of $X^{n\ell}$, it follows
that also the pushforward $\sigma^{n\ell}_*\O_{\wt\X^n}(\widetilde g)|_{X^{n\ell}}$ 
is invertible at $N_e$.  Finally, the  gluing data of 
$I^n_f$ and $I^n_h$ at $N_e$ are 
$b_ea_e^{\dl^{\m,n}_f(e)}$ and $b_ea_e^{\dl^{\m,n}_h(e)}$, 
respectively, as given by the coordinates of $p^n_f$ and
$p^n_h$. As for the gluing data of
$\sigma^{n\ell}_*\O_{\wt\X^n}(\widetilde g)|_{X^{n\ell}}$ 
at $N_e$, since $\mathfrak
p_e=0$, it is the same as if we assumed that $\wt g$ is the canonical extension
of $g$ with respect to the divisor 0; it follows from
Subsection~\ref{ges} that the gluing data is 
$a_e^{\dl^{0,n}_g(e)}$. Clearly, $\dl^{0,n}_g(e)=\dl^{\mathfrak
  p,n}_g(e)$, and hence it follows from \eqref{dd0} that 
the sheaves in \eqref{Lnhf} have the same gluing data in a
neighborhood of $Z^{e,n}$. 

\smallskip

Since the chains $Z^{e,n}$ are curves of compact type, an invertible
sheaf on each of them is determined by its degrees on the
components. It follows 
that the sheaves in \eqref{Lnhf} coincide up to the chosen gluing
along the intersections of the chains $Z^{e,n}$ with the rest of
$X^{n\ell}$ for the $e\in\E$ such that either $\dl_h^{\m,n}(e)$ or
$\dl_f^{\m,n}(e)$ is not an integer.

In fact, we have not yet specified the gluing giving $L^n_f$ along $Z^{e,n}$ if
$I^n_f$ is not invertible at $N_e$. We do it now: For each oriented edge
$e=uv\in E^{\mathfrak o}$ such 
that $\mathfrak d^{\mathfrak m,n}_f(e)\not\in\Z$, give $L^n_f$ the
gluing along $Z^{e,n}$ such 
that the pushforwards under $\sigma^{n\ell}$ of 
\begin{equation}\label{Lnff}
L^n_{f_e}\quad\text{and}\quad
L^n_f\otimes\O_{\wt\X^n}\big(\,\widetilde
r\,\big)|_{X^{n\ell}}
\end{equation}
have the same gluing data at $N_e$, where $\widetilde r$ is the
canonical extension of $r:=(n\ell_e-i^n_e(f))\chi_{\indm u}$ with
respect to the divisor $D^n_f$, and $f_e:=f+r$. Notice 
that $I^n_{f_e}$ is invertible at $N_e$, and thus the gluing data for
$L^n_{f_e}$ along $Z^{e,n}$ is determined from that of $I^n_{f_e}$ at
$N_e$. 
From the above prescription,
the gluing data for $L^n_f$ is now fixed everywhere for each 
$f\in C^0(G,\Z)$ and $n\in\N$.

Finally, with the above choices, we can claim more than
Claim~\ref{cl1}:

\begin{claim}\label{cl2} Notation as above, for each $f,h\in C^0(G,\Z)$ and 
$n\in\N$, we have
$$
L^n_h\cong L^n_f\otimes\O_{\wt\X^n}\big(\,\widetilde g\,\big)|_{X^{n\ell}},
$$
where $\widetilde g$ is the canonical extension of $g:=h-f$ with respect to the divisor
$D^n_f$ in $H^n$ associated to $L^n_f$.
\end{claim}

\begin{proof} 
Indeed, given $e\in E^{\mathfrak o}$, notice first that, 
if $q\in C^0(G,\Z)$ is such that $I^n_q$ is invertible at
$N_e$, we have proved that the gluing data along $Z^{e,n}$ of 
$$
L^n_q\quad\text{and}\quad L^n_{f_e}\otimes
\O_{\wt\X^n}\big(\,\widetilde
g\,\big)|_{X^{n\ell}}
$$
coincide, where $g:=q-f_e$ and $\widetilde g$ is the
canonical extension of $g$ with respect to the divisor $D^n_{f_e}$ in
$H^n$ associated to $L^n_{f_e}$. It follows that the gluing data along
$Z^{e,n}$ of 
$$
L^n_q \quad\text{and}\quad L^n_f\otimes
\O_{\wt\X^n}\big(\,\widetilde g+\widetilde
r\,\big)|_{X^{n\ell}}
$$
coincide, where $\widetilde r$ is the
canonical extension of $r:=(n\ell_e-i^n_e(f))\chi_{\indm u}$ with
respect to the divisor $D^n_f$. But $D^n_{f_e}=D^n_f+\text{div}(\wt r)$ by Claim~\ref{cl1}. 
It thus follows from \cite{AE1},
Prop.~2.11, that $\widetilde r+\widetilde g$ is the canonical extension
of $r+g$ with respect to the divisor $D^n_f$. Notice that $r+g=q-f$.

Thus, given $q\in C^0(G,\Z)$ such that $\dl^{\m,n}_q(e)\in\Z$, the gluing
data along $Z^{e,n}$ of these three sheaves,
\begin{equation}\label{3shf}
L^n_q,\quad L^n_h\otimes
\O_{\wt\X^n}\big(\,\widetilde
p\,\big)|_{X^{n\ell}}\quad\text{and}\quad L^n_f\otimes
\O_{\wt\X^n}\big(\,\widetilde g\,\big)|_{X^{n\ell}},
\end{equation}
coincide, where $\widetilde p$ (resp.~$\widetilde g$) is the canonical extension
of $p:=q-h$ (resp.~$g:=q-f$) with respect to the divisor $D^n_h$
(resp.~$D^n_f$). But then so do the gluing data of 
$$
L^n_h\otimes
\O_{\wt\X^n}\big(\,\widetilde p+\widetilde
r\,\big)|_{X^{n\ell}}\quad\text{and}\quad
L^n_f\otimes
\O_{\wt\X^n}\big(\,\widetilde g+\widetilde r\,\big)|_{X^{n\ell}},
$$
where $\widetilde r$ is the canonical extension
of $r:=h-q$ with respect to the divisor $D^n_q$. Since all three
sheaves in \eqref{3shf} have the same associated divisor in $H^n$, it
follows again from \cite{AE1},
Prop~2.11, that $\widetilde p+\widetilde
r$ is the canonical extension of $p+r=0$ and that 
$\widetilde g+\widetilde r$ is the canonical extension
of $g+r=h-f$ with respect to $D^n_f$. Since $D^n_f$ is $G$-admissible,
$\widetilde p+\widetilde r=0$. So the sheaves above are
exactly those in \eqref{Lnhf}, and we have concluded that they have
the same gluing data along $Z^{e,n}$ for every $e\in E^{\mathfrak o}$,
whence our claim.
\end{proof}

We need two more claims to finish the proof of the theorem.

\begin{claim}\label{cl3} Notation as above, $L^n_0$ is the pullback of $L^1_0$
  under the map $X^{n\ell}\to X^\ell$.
\end{claim}

\begin{proof}
Notice that $p^n_0=p^1_0$ for every $n\in\Z$, since
$\dl^{\m,n}_0=\dl^{\m,1}_0$, and thus $I^n_0=I^1_0$. Thus the restrictions of $L^n_0$
and $L^1_0$ to the components $X_v$ for $v\in V(G)$ coincide. Also,
$i^n_e(0)=ni^1_e(0)$ for each $e\in\E$,  and thus $D^n_0$ is the
pullback of $D^1_0$. Thus the restriction of $L^n_0$ to $Z^{e,n}$
coincides with that of the pullback of $L^1_0$ for each $e\in\E$. 

Also, 
the gluing data of $L^n_0$
along $Z^{e,n}$ for each $e\in E^{\mathfrak o}$ such that
$\dl^{\m,1}_0(e)\in\Z$ is the gluing data of the pullback of
$L^1_0$. On the other hand, by construction, for each $e=uv\in E^{\mathfrak o}$ such that 
$\dl^{\m,1}_0(e)\not\in\Z$, the gluing data of $L^n_0$ along
$Z^{e,n}$ is such that the gluing data at $N_e$ of the pushforwards of 
$$
L^n_{g_n}\quad\text{and}\quad 
L^n_0\otimes\mathcal O_{\wt\X^n}\big(\,\widetilde
g_n\,\big)|_{X^{n\ell}}
$$
coincide, where $g_n:=(n\ell_e-i^n_e(0))\chi_{\indm u}$ and $\widetilde g_n$
is the canonical extension of $g_n$ with respect to $D^n_0$. Same for
$L^1_0$, for $n$ replaced by 1. Clearly,
since $i^n_e(0)=ni^1_e(0)$, we have that $g_n=ng_1$, 
whence $\mathcal O_{\wt\X^n}\big(\,\widetilde
g_n\,\big)$ is the pullback of $\mathcal O_{\wt\X}\big(\,\widetilde
g_1\,\big)$ to $\wt\X^n$. Finally,
the gluing data of the pushforward $I^n_{g_n}$ at $N_e$ is given by
$p^n_{g_n}$, whereas that of
$I^1_{g_1}$ is given by $p^1_{g_1}$. Since $g_n=ng_1$, we have that
$\dl^{\m,n}_{g_n}=\dl^{\m,1}_{g_1}$, and thus
$p^n_{g_n}=p^1_{g_1}$. It follows that the gluing data of $L^n_0$ and
of the pullback of $L^1_0$ coincide everywhere, proving the claim.
\end{proof}

\begin{claim}\label{cl4} Notation as above, there is an almost invertible sheaf $\L$ on
$\wt\X$ whose restriction to $X^\ell$ is $L^1_0$.
\end{claim}

\begin{proof} 
The sheaf $L^1_0$ fails to be
invertible only at the nodes of $X^\ell$ lying over nodes of $X$ other
than the $N_e$, at which
the corresponding $K_v$ is not invertible. Let $\widehat\X$ denote the
blowup of $\wt\X$ at these nodes and by $\widehat X$ the special
fiber of $\widehat\X/B$. The curve $\widehat X$ is obtained by
separating the branches of every node in the center of the blowup and
connecting them by a smooth rational curve. Then $L^1_0$ is the
pushforward of an invertible sheaf on $\widehat X$, having degree
one on the added rational curves, by \cite{EP16}, Prop.~5.5,
p.~77. Since this sheaf is invertible, since the relative Picard scheme of
$\widehat\X/B$ is smooth over $B$, and since $B$ is the spectrum of
$\k[[t]]$, it follows that it extends to an invertible sheaf on
$\widehat\X$. The pushforward of this extension under 
the blowup map $\widehat\X\to\wt X$ is the required sheaf $\L$.
\end{proof}

\begin{proof}[Proof of Theorem~\ref{regeneration}] We claim that 
$L_\eta:=\L|_{X_\eta}$ is the required
sheaf. Indeed, we need only prove
    that $L^n_h=\mathcal L^n(h)|_{X^{n\ell}}$ for each $h\in  C^0(G,\Z)$
    and $n\in\N$, where $\mathcal L^n$ is the pullback of $\L$ to
    $\wt\X^n$. But notice that $\mathcal L^n|_{X^{n\ell}}=L^n_0$ by 
Claims~\ref{cl3}~and~\ref{cl4}. And
$$
\mathcal L^n(h)=\mathcal L^n\otimes\O_{\wt\X^n}\big(\,\widetilde h\,\big)
$$
for each $h\in C^0(G,\Z)$, where $\wt h$ is the canonical extension
of $h$ with respect to the divisor $D^n_0$ on $H^n$ associated to
$L^n_0$.  Applying Claim~\ref{cl2} with $f=0$ finishes the proof.
\end{proof}

 \vspace{.7cm}

\subsection*{Acknowledgements.}  This project benefited very much from
the hospitality of the Mathematics Department at the \'Ecole Normale
Sup\'erieure (ENS) in Paris and the Instituto de Matem\'atica Pura e
Aplicada (IMPA) in Rio de Janeiro during mutual visits of both authors.
We thank the two institutes and their members for providing for those visits.
We are also specially grateful to the Brazilian-French Network in
Mathematics for providing support for a visit of E.E.~to ENS Paris
and a visit of O.A.~to IMPA. 
 
 \bibliographystyle{alpha}
\bibliography{bibliography}
\end{document}